\documentclass[11pt]{amsart}

\usepackage{a4wide}
\usepackage{bbm}
\usepackage{t1enc} 
\usepackage{epsf}
\usepackage{amssymb}
\usepackage{graphicx}
\usepackage{psfrag}

\newtheorem{theorem}{Theorem}
\newtheorem{lem}{Lemma}
\newtheorem{prop}{Proposition}
\newtheorem{coro}{Corollary}

\theoremstyle{definition}
\newtheorem{defi}{Definition}
\newtheorem{rem}{Remark}
\newtheorem{ex}{Example}

\newcommand{\kk}{\mathbbm{k}}
\newcommand{\zetan}{\zeta_{n}}
\newcommand{\zetam}{\zeta_{m}}

\newcommand{\Z}{\mathbbm{Z}}
\newcommand{\N}{\mathbbm{N}}
\newcommand{\Zn}{\mathbbm{Z}[\zeta_{n}]}

\newcommand{\Q}{\mathbbm{Q}}
\newcommand{\Qn}{\mathbbm{Q}(\zeta_{n})}

\newcommand{\R}{\mathbbm{R}}
\newcommand{\C}{\mathbbm{C}}
\newcommand{\OO}{\mathcal{O}}
\newcommand{\oo}{\thinspace\scriptstyle{\mathcal{O}}\displaystyle}

\begin{document}

\title[Discrete Tomography]{Uniqueness in Discrete Tomography\linebreak
  of Planar Model Sets}

\author[C. Huck]{Christian Huck}
\address{\hspace*{-1em} Fakult\"{a}t f\"{u}r Mathematik,
  Universit\"{a}t Bielefeld, Postfach 10 01 31, 33501 Bielefeld, Germany}
\email{huck@math.uni-bielefeld.de}
\urladdr{http://www.math.uni-bielefeld.de/baake/huck/}

\begin{abstract} 
The problem of determining finite subsets of characteristic planar model sets
(mathematical quasicrystals) $\varLambda$, called cyclotomic model
sets, 
by parallel $X$-rays is considered. Here, an $X$-ray in direction
$u$ of a finite
subset of the plane gives the number of points in the set on each line
parallel to $u$. For practical reasons, only $X$-rays in
$\varLambda$-directions, i.e., directions parallel to
non-zero elements of the difference set $\varLambda - \varLambda$, are
permitted. In particular, by combining methods from algebraic number theory and
convexity, it is shown that the convex subsets of a
cyclotomic model set $\varLambda$, i.e., finite
sets $C\subset \varLambda$ whose convex hulls contain no new points of
$\varLambda$, are determined, among all convex subsets of
$\varLambda$, by their $X$-rays in four prescribed  $\varLambda$-directions, whereas any set of three $\varLambda$-directions does not suffice for this purpose. We also study the interactive technique of successive
determination in the case of cyclotomic model sets, in which the
information from previous $X$-rays is used in deciding the direction
for the next $X$-ray. In particular, it is shown that the finite subsets of any
cyclotomic model set $\varLambda$ can be successively
determined by two $\varLambda$-directions. All results are illustrated by means of well-known examples, i.e., the cyclotomic model sets associated with the square tiling, the
triangle tiling, the tiling of Ammann-Beenker, the T\"ubingen triangle tiling and the shield tiling.
\end{abstract}

\maketitle

\section{Introduction}

{\em Discrete tomography} is concerned with the 
inverse problem of retrieving information about some discrete object from 
(generally noisy) information about its incidences with certain query sets.
A typical example is the {\em reconstruction} of a finite point set
from its line sums in a small number of directions. More precisely, a ({\em discrete parallel}\/) {\em
  X-ray} of a finite subset of
Euclidean $d$-space $\R^d$ in direction $u$ gives the
number of points in the set on each line in $\R^d$ parallel to
$u$. In the
 classical setting, motivated by crystals, the positions to be determined live on the square lattice $\Z^2$ or, more
generally, on arbitrary lattices $L$ in $\R^d$, where $d\geq 2$. (Here,
 a subset $S$ of $\R^{d}$ is said
  to {\em live on} a subgroup $G$ of $\R^{d}$ when its
  difference set $S-S:=\{s-s'\,|\,s,s'\in S\}$ is a subset of $G$.) 
In fact, many of the problems in discrete tomography
have been studied on $\Z^2$, the classical planar setting of
discrete tomography; see~\cite{HK}, \cite{Gr} and~\cite{GGP}. In the longer run, by
also having other structures than perfect crystals in mind, one has to
take into account wider classes of sets, or at least
significant deviations from the lattice structure. As an intermediate
step between periodic and random (or amorphous) Delone sets (defined below), we consider systems of {\em aperiodic
  order}, more precisely, of so-called {\em model sets} (or {\em
  mathematical quasicrystals}), which are commonly accepted as a good mathematical model for quasicrystalline structures in
nature~\cite{St}. It is interesting to know whether these sets behave like lattices as far as discrete tomography is concerned.

The main motivation for our interest in the discrete tomography of
model sets comes from the physical existence of quasicrystals that can
be described as aperiodic model sets together with the demand of materials
science to reconstruct three-dimensional (quasi)crystals or planar
layers of them from their images under quantitative {\em high
  resolution transmission electron microscopy} (HRTEM) in a small
number of directions. In fact, in~\cite{ks,sk} a technique is described, based
on HRTEM, which can effectively measure the number of atoms lying on
lines parallel to certain directions; it is called QUANTITEM
(quantitative analysis of the information coming from transmission
electron microscopy). At present, the measurement of the number
of atoms lying on a line can only be achieved for some crystals; see~\cite{ks,sk}. However, it is reasonable to expect that future developments in technology will improve this situation.  

Whether or not one has future applications in materials science in
mind, the starting point will always be a specific structure
model. This means that the specific type of the (quasi)crystal is
given (see~\cite{B,schw} for the concept), and one is confronted with
the $X$-ray data of an unknown finite subset of it. This is the
correct analogue of starting with translates of $\Z^{d}$ in the classical
setting. 

Here, our investigation of the
discrete tomography of systems with aperiodic order will be restricted to the study of two-dimensional systems. Fortunately, solving the
 problems of discrete tomography for two-dimensional systems with aperiodic
order also lies at the heart
of solving the corresponding problems in three dimensions. This is due
to the fact that model sets
 possess, as it is the case for lattices, a dimensional hierarchy, meaning that any model set in
 $d$ dimensions can be sliced into model sets of dimension
 $d-1$.

More precisely, we consider a well-known class of \emph{planar}
model sets which can be described in algebraic terms. It will be a
nice property of this class of planar model sets that it contains, as
a special case, arbitrary translates of the square lattice $\Z^2$,
which corresponds to the classical planar setting of discrete tomography. 

Let us become more specific. By using the Minkowski representation of algebraic number fields,
we introduce, for $n \notin \{1,2\}$, the corresponding class of {\em cyclotomic model sets} $\varLambda \subset \mathbbm{C}\cong \R^2$
which live on $\Zn$, where $\zeta_{n}$ is a
primitive $n$th root of unity in $\C$, e.g., $\zeta_{n}=e^{2\pi
    i/n}$. The $\Z$-module $\Zn$ is the ring of integers in the $n$th cyclotomic
field $\mathbbm{Q}(\zeta_{n})$, and, for $n \notin \{1,2,3,4,6\}$, when
viewed as a subset of the plane, is dense, whereas cyclotomic model
sets $\varLambda$ are {\em Delone sets}, i.e., they are uniformly
discrete and relatively dense. In fact, model sets are even {\em
  Meyer sets}, meaning that also $\varLambda-\varLambda$ is uniformly
discrete; see~\cite{Moody}. It turns out that, except the cyclotomic model sets living
on $\Zn$ with $n \in \{3,4,6\}$ (these are exactly the translations of
the square and the triangular lattice, respectively), cyclotomic model
sets $\varLambda$ are {\em aperiodic}, meaning that they have no
translational symmetries at all. Well-known examples are the planar model sets with $N$-fold
cyclic symmetry that are associated with the square tiling ($n=N=4$), the
triangle tiling ($2n=N=6$), the Ammann-Beenker tiling ($n=N=8$), the
T\"ubingen triangle tiling ($2n=N=10$) and the shield tiling
($n=N=12$). Note that orders $5,8,10$ and $12$ occur as standard
cyclic symmetries of genuine quasicrystals~\cite{St}, which are
basically stacks of planar aperiodic layers. 

In general, the {\em reconstruction problem} of discrete tomography can possess rather
different solutions. This leads to the investigation of the
corresponding {\em uniqueness problem}, i.e., the (unique) {\em
  determination} of finite subsets of a fixed cyclotomic model set
$\varLambda$ by $X$-rays in a
 small number of suitably prescribed directions. For practical reasons, only $X$-rays in
$\varLambda$-directions, i.e.,
directions parallel to non-zero elements of the difference set
$\varLambda-\varLambda$ of $\varLambda$,
are permitted.  
We say that a subset $\mathcal{E}$ of the set of all finite subsets of
a fixed cyclotomic model set $\varLambda$ is determined by the
$X$-rays in a finite set $U$ of directions if
different sets $F$ and $F'$ in $\mathcal{E}$ have different $X$-rays,
i.e., if there is a direction $u\in U$ such that the $X$-rays of $F$ and $F'$ in direction $u$ differ. 

Without the restriction to $\varLambda$-directions, one
can prove that the finite subsets of a fixed cyclotomic
model set $\varLambda$ can be determined by one $X$-ray (Proposition~\ref{irratdir}). In fact, it
turns out that any $X$-ray in a non-$\varLambda$-direction is suitable for this purpose. However, in
practice (HRTEM), $X$-rays in non-$\varLambda$-directions are meaningless
since the resolution coming from such $X$-rays would not be good
enough to allow a quantitative analysis -- neighbouring lines are not
sufficiently separated. 

Proposition~\ref{m+1} demonstrates that the finite sets $F$ of cardinality less
than or equal to some $k\in \mathbbm{N}$ in a fixed cyclotomic model
set $\varLambda$ are determined by any set of $k+1$ $X$-rays in
pairwise non-parallel $\varLambda$-directions. However, in practice, one is
interested in the determination of finite sets by $X$-rays in a small
number of directions since after about $3$ to $5$ images taken by
HRTEM, the object may be damaged or even destroyed by the radiation
energy. Observing that the typical atomic structures to be determined
comprise about $10^6$ to $10^9$ atoms, one realizes that the last
result is not practical at all. In fact, it is also shown that any fixed number of $X$-rays in
$\varLambda$-directions is insufficient to determine the entire class of
finite subsets of a fixed cyclotomic model set $\varLambda$ (Proposition~\ref{source}). 

In view of this, it is necessary to impose some restriction in order to obtain positive uniqueness results. 

As a first option, we restrict the set of finite subsets of a fixed
cyclotomic model set $\varLambda$ under consideration, i.e., we
consider for every $R>0$ the class of {\em bounded} subsets of
$\varLambda$ with diameter less than $R$. Since $\varLambda$ is
uniformly discrete, bounded subsets of
$\varLambda$ are finite. It is shown that, for all $R>0$ and for all
cyclotomic model sets $\varLambda$, there are two non-parallel
prescribed $\varLambda$-directions such that the bounded subsets of $\varLambda$ with
diameter less than $R$ are determined,
among all such subsets of $\varLambda$, by their $X$-rays in these
directions (Corollary~\ref{coromod}). In fact, in
Theorem~\ref{bounded}, it is shown that the corresponding result holds
for arbitrary Delone sets $\varLambda\subset\R^d$, $d\geq 2$, of 
{\em finite local complexity}, the latter meaning that the difference set $\varLambda-\varLambda$ is
 closed and discrete. Thereby, also the case of translates of lattices
 $L$ in $\R^{d}$, $d\geq 2$, is included (Corollary~\ref{bounded2}). Moreover, even the corresponding results for
orthogonal projections on orthogonal complements of $1$-dimensional (respectively, $(d-1)$-dimensional)
subspaces generated by elements of $\varLambda-\varLambda$ 
(resp., $L-L=L$) hold; the multiplicity information coming along with $X$-rays
is not needed. Unfortunately, these results are limited in practice
because, in general, one cannot make sure that all the directions
which are used yield images of high enough resolution.

As a second option, we restrict the set of finite subsets of a fixed
cyclotomic model set $\varLambda$ by
considering the class of {\em convex} subsets of
$\varLambda$. They are finite sets $C\subset \varLambda$ whose convex
hulls contain no new points of $\varLambda$, i.e., finite sets
$C\subset \varLambda$ with $C = \operatorname{conv}(C)\cap
\varLambda$. In Theorem~\ref{main1}, it is shown that there are four pairwise non-parallel
prescribed $\varLambda$-directions such that the convex subsets of any
cyclotomic model set $\varLambda$ living on $\Zn$ are determined,
among all convex subsets of $\varLambda$, by their $X$-rays in these
directions. It is further shown that three pairwise non-parallel $\varLambda$-directions never
suffice for this purpose, whereas {\em any} seven pairwise
non-parallel $\varLambda$-directions which satisfy a certain relation (see
below for details) have the property that the convex subsets of any
cyclotomic model set $\varLambda$ are determined,
among all convex subsets of $\varLambda$, by their $X$-rays in these
directions. It is also demonstrated that, in a sense, the number seven cannot be reduced to six in the latter result. These results heavily
depend on a massive restriction of the set of $\varLambda$-directions. We
shall also remove this restriction and present explicit results in the case of the
important subclass of cyclotomic
model sets $\varLambda$ with co-dimension $2$, i.e., cyclotomic
model sets living on $\Zn$, where $n\in\{5,8,10,12\}$. Here, by using
$p$-adic valuations, we are able to present a more suitable analysis which
shows that uniqueness will be provided by any set of four $\varLambda$-directions
whose slopes (suitably ordered) yield a cross ratio that does not map under the field norm of the field
extension $\Q(\zeta_{n}+\bar{\zeta}_{n})/\Q$ to a certain {\em finite
  set} of rational numbers (Theorem~\ref{main2}). (It will turn out
that these cross
ratios are indeed elements of the maximal real subfield
$\Q(\zeta_{n}+\bar{\zeta}_{n})$ of the $n$th cyclotomic field
$\Q(\zeta_{n})$.) It is also shown that,
by the same analysis, similar results can be obtained for arbitrary
cyclotomic model sets (Theorem~\ref{main3}). In the case of cyclotomic
model sets with co-dimension $2$, we shall be able to
present sets of four $\varLambda$-directions which
provide uniqueness and yield
images of rather high resolution in the quantitative HRTEM of the
corresponding (aperiodic) cyclotomic model sets, the latter making these results look promising.  

A major task in achieving the above results involves examining
so-called $U$-{\em polygons} in cyclotomic model sets $\varLambda$,
which exhibit a weak sort of regularity. In the context of
$U$-polygons, the question for affinely regular polygons in the plane
with all their vertices in a fixed cyclotomic model set $\varLambda$
arises rather naturally. Chrestenson~\cite{C} has shown that any
(planar) regular polygon whose vertices are contained in
$\mathbbm{Z}^{d}$ for some $d\geq 2$ must have $3,4$ or $6$
vertices. More generally, Gardner and Gritzmann~\cite{GG} have
characterized the {\em affinely} regular lattice polygons, i.e.,
images of non-degenerate regular polygons under a non-singular affine
transformation $\R^{2} \longrightarrow \R^{2}$ whose vertices are
contained in $\mathbbm{Z}^{2}$. It turned out that these are precisely
the affinely regular triangles, parallelograms and hexagons. Clearly,
their result remains valid when one replaces the square lattice
$\mathbbm{Z}^{2}$ by the triangular lattice $\mathbbm{Z}[\zeta_{3}]$,
since the latter is the image of the first under an invertible linear
transformation of the plane. Observing that both the square and
triangular lattice are examples of rings of integers $\Zn$ in
cyclotomic fields $\Qn$, one can ask for a generalization of the above
result to all these objects, viewed as $\Z$-modules in the
plane. Using standard results from algebra and algebraic number
theory, we provide a characterization in terms of a simple
divisibility condition (Theorem~\ref{th1}). Moreover, we show that our
characterization of affinely regular polygons remains valid when
restricted to the corresponding class of cyclotomic model sets
$\varLambda$ (Corollary~\ref{th1mod}).  

As a third option, we consider the interactive technique of {\em successive determination},
introduced by Edelsbrunner and Skiena~\cite{ES} for continuous $X$-rays, in the
case of aperiodic cyclotomic model sets, in which the information from previous
$X$-rays may be used in deciding the direction for the next $X$-ray. Here,
it is shown that the finite subsets of any cyclotomic model set
$\varLambda$ can be successively determined by
$X$-rays in two non-parallel $\varLambda$-directions. In fact, it is shown that, for $n \notin \{1,2\}$, the
finite subsets of the ring of cyclotomic integers $\Zn$ can be successively determined by
$X$-rays in two non-parallel $\Zn$-directions, i.e., directions
parallel to a non-zero element of $\Zn$ (Theorem~\ref{mainsucdet}). Again, the corresponding results do even hold for
orthogonal projections on orthogonal complements of $1$-dimensional subspaces generated by
elements of $\varLambda-\varLambda$; the multiplicity information coming alsong with $X$-rays
is again not needed. Unfortunately, these results are again limited in practice
because, in general, one cannot make sure that all the directions
which are used yield images of high enough resolution.

Previous studies have focussed on the 
`anchored' case that the underlying ground set is located in a 
linear space, i.e., in a space with a specified location of the origin.
The $X$-ray data is then taken with respect to this localization. Note that
the rotational orientation of a (quasi)crystalline probe in an electron microscope can
rather easily be done in the diffraction mode, prior to taking images
in the high-resolution mode, though, in general, a natural choice of a
translational origin is {\em not} possible. This problem has its origin in the practice of quantitative HRTEM since, in general, the
$X$-ray information does not allow us to locate the examined sets. We
believe that, in order to obtain applicable results, one has to deal
with this `non-anchored' case. 

In discussing inverse problems, it is
important to distinguish between determination and reconstruction. The
problem of reconstructing finite subsets of cyclotomic model sets $\varLambda$
given $X$-rays in {\em two} $\varLambda$-directions is treated
in~\cite{BG2}. There, it is shown that, for a large class of
cyclotomic model sets, this problem can be solved in
polynomial time in the real RAM-model of computation. For corresponding results in the lattice case, see~\cite{GGP}.     

Note that this work provides a generalization of well-known results,
obtained by Gardner and Gritzmann~\cite{GG} for the `anchored' case of
fixed translates of lattices $L$ in
$\R^{d}$, $d\geq 2$, to the class of cyclotomic model sets. Additionally,
we shall present results which deal with the `non-anchored' case described
above. Several of our proofs are only a slight modification of the
proofs from~\cite{GG}. The
key role in this text is played by Lemma~\ref{dilate}, which rests on
a result from the theory of {\em Pisot-Vijayaraghavan}
numbers~\cite{Sa}, and $p$-adic valuations. For more details, other settings, and the historical background of the above uniqueness problems, we would like to refer the reader to~\cite{GG} and references therein.

\section{Algebraic Background, Definitions and Notation}\label{found}

Natural numbers are always assumed to be positive, i.e.,
$\mathbbm{N}=\{1,2,3,\dots\}$. Throughout the text, we use the
convention that the symbol $\,\subset\,$ includes equality. We denote the norm in Euclidean $d$-space $\mathbbm{R}^{d}$
by $\Arrowvert \cdot \Arrowvert$. The unit sphere in
$\mathbbm{R}^{d}$ is denoted by $\mathbb{S}^{d-1}$, i.e.,
$$
\mathbb{S}^{d-1}=\left\{\left. x\in\R^d\,\right |\,\Arrowvert x \Arrowvert =1\right\}\,.
$$
Moreover, the elements of $\mathbb{S}^{d-1}$ are also called
{\em directions}. If $x\in\R$, then $\lfloor x \rfloor$
denotes the greatest integer less than or equal to $x$. For $r>0$ and $x\in\R^{d}$,
$B_{r}(x)$ is the open ball of radius $r$ about $x$. If
$k,l\in \mathbbm{N}$, then $(k,l)$
and $[k,l]$ denote their greatest common divisor
and least common multiple, respectively. There should be no
  confusion with the usual notation of open (resp., closed) and
  bounded real intervals. For a subset $S\subset
\mathbbm{R}^{d}$, $k\in \mathbbm{N}$ and $R>0$, we denote by
$\operatorname{card}(S)$, $\mathcal{F}(S)$, $\mathcal{F}_{\leq k}(S)$, $\mathcal{D}_{<R}(S)$,
$S^{\circ}$, $\overline{S}$, $\partial S$, $\langle S \rangle_{\Z}$, $\langle S \rangle_{\Q}$,
 $\langle S \rangle_{\R}$, $\operatorname{conv}(S)$, $\operatorname{diam}(S)$
and $\mathbbm{1}_{S}$ the cardinality, the set of finite subsets, the
set of finite subsets of $S$ having cardinality less than or equal to
$k$, the set of subsets of $S$
with diameter less than $R$, interior, closure,
boundary, $\Z$-linear hull, $\Q$-linear hull,
 $\R$-linear hull, convex hull, diameter and
characteristic function of $S$, respectively. The notation of the
closure operator should not be confused with the usual notation of
complex conjugation. The {\em
  dimension} of $S$ is the dimension of its affine hull
$\operatorname{aff}(S)$, and is denoted by
$\operatorname{dim}(S)$. Further, a
 linear subspace $T$ of $\mathbbm{R}^{d}$ is called an
$S${\em-subspace} if it is generated by elements of the difference set
$S-S:=\{s-s'\,|\,s,s'\in S\}$. A direction
$u\in\mathbb{S}^{d-1}$ is called an $S${\em-direction} if it is
parallel to a non-zero element of $S-S$.
 If $T$ is a linear subspace of $\mathbbm{R}^{d}$, we denote the  
orthogonal projection of an element $x$ of $\mathbbm{R}^{d}$ on $T$ by
$x|T$. Moreover, we denote the
 orthogonal projection of a subset $S$ of $\mathbbm{R}^{d}$ on $T$ by
$S|T$. The orthogonal complement of $T$ is denoted by $T^{\perp}$. If we use this notation with respect to
 other than the  
{\em canonical} inner product on $\mathbbm{R}^{d}$, we
shall say so explicitly. The
symmetric difference of two sets $A$ and $B$ is defined as
$A\,\triangle \, B:=(A\setminus B) \cup
(B\setminus A)$. A subset $S$ of $\R^{d}$ is said
  to {\em live on} a subgroup $G$ of $\R^{d}$ when its
  difference set $S-S$ is a subset of $G$. Obviously, this is
  equivalent to the existence of a suitable $t\in \R^{d}$ such
that $S\subset t+G$. Finally, the {\em centroid} of an element
$F\in\mathcal{F}(\R^d)$ is defined as
$(\sum_{f\in F}f)/\operatorname{card}(F)$.

\begin{defi}\label{xray..}
Let $d\in \mathbbm{N}$ and let $F\in \mathcal{F}(\mathbbm{R}^{d})$. Furthermore, let $u\in \mathbb{S}^{d-1}$ be a direction and let $\mathcal{L}_{u}$ be the set of lines in direction $u$ in $\mathbbm{R}^{d}$. Then, the ({\em discrete parallel}\/) {\em X-ray} of $F$ {\em in direction} $u$ is the function $X_{u}F: \mathcal{L}_{u} \longrightarrow \mathbbm{N}_{0}:=\mathbbm{N} \cup\{0\}$, defined by $$X_{u}F(\ell) := \operatorname{card}(F \cap \ell\,) =\sum_{x\in \ell} \mathbbm{1}_{F}(x)\,.$$ Moreover, the {\em support} of $X_{u}F$, i.e., the set of lines in $\mathcal{L}_{u}$ which pass through at least one point of $F$, is denoted by $\operatorname{supp}(X_{u}F)$. For $z\in \mathbbm{R}^{d}$, we denote by $\ell_{u}^{z}$ the element of $\mathcal{L}_{u}$ which passes through $z$.
\end{defi}

\begin{rem}
In the situation of Definition~\ref{xray..}, $\operatorname{supp}(X_{u}F)$ is finite and, moreover, the cardinality of $F$ is implicit in the $X$-ray, since one has $$\sum_{\ell \in \operatorname{supp}(X_{u}F)}X_{u}F(\ell) = \operatorname{card}(F)\,.$$ 
\end{rem}

\begin{lem}\label{cardinality}
Let $d\in \mathbbm{N}$ and let $u\in \mathbb{S}^{d-1}$ be a
direction. If $F, F'\in \mathcal{F}(\mathbbm{R}^{d})$, one has: 
\begin{itemize}
\item[(a)]
$X_{u}F=X_{u}F'$ implies $\operatorname{card}(F)=\operatorname{card}(F')$.
\item[(b)]
If $X_{u}F=X_{u}F'$, the centroids of $F$ and $F'$ lie on the same line parallel to $u$.
\end{itemize}
\end{lem}
\begin{proof}
See ~\cite[Lemma 5.1 and Lemma 5.4]{GG}.
\end{proof}

\begin{defi}
Let $d\geq 2$, let $\mathcal{E}\subset
\mathcal{F}(\mathbbm{R}^{d})$, and let $m\in\N$. Further, let $U$ be a
finite set of directions and let $\mathcal{T}$ be a
finite set of $1$-dimensional subspaces $T$ of $\mathbbm{R}^{d}$. 
\begin{itemize}
\item[(a)]
We say that $\mathcal{E}$ is {\em determined} by the $X$-rays in the directions of $U$ if, for all $F,F' \in \mathcal{E}$, one has
$$
(X_{u}F=X_{u}F'\;\,\forall u \in U) \;  \Longrightarrow\; F=F'\,.
$$
\item[(b)]
We say that $\mathcal{E}$ is {\em determined} by the (orthogonal)
projections on the othogonal complements $T^{\perp}$ of the subspaces
$T$ of $\mathcal{T}$ if, for all $F,F' \in \mathcal{E}$, one has
$$
(F|T^{\perp}=F'|T^{\perp}\;\,\forall T \in \mathcal{T}) \;  \Longrightarrow\; F=F'\,.
$$
\item[(c)]
We say that $\mathcal{E}$ is {\em successively determined} by the
$X$-rays in the directions of $U$, say $U=\{u_1,\dots,u_{m}\}$, if,
for a given $F\in \mathcal{E}$, these can be chosen inductively (i.e., the choice of $u_{j}$ depending on all $X_{u_{k}}F$ with \mbox{$k\in\{1,\dots,j-1\}$)} such that, for all $F' \in \mathcal{E}$, one has
$$
(X_{u}F'=X_{u}F\;\,\forall u \in U)\;  \Longrightarrow\; F'=F\,.
$$
\item[(d)]
We say that $\mathcal{E}$ is {\em successively determined} by the (orthogonal)
projections on the othogonal complements $T^{\perp}$ of the subspaces
$T$ of $\mathcal{T}$
, say $\mathcal{T}=\{T_1,\dots,T_{m}\}$, if,
for a given $F\in \mathcal{E}$, these can be chosen inductively (i.e., the choice of $T_{j}$ depending on all $F|T_{k}^{\perp}$ with $k\in\{1,\dots,j-1\}$) such that, for all $F' \in \mathcal{E}$, one has
$$
(F|T^{\perp}=F'|T^{\perp}\;\,\forall T \in \mathcal{T}) \;  \Longrightarrow\; F=F'\,.
$$
\item[(e)]
We say that $\mathcal{E}$ is {\em determined} (resp., {\em successively determined}) by $m$ $X$-rays if there is a set $U$ of $m$ pairwise non-parallel directions such that $\mathcal{E}$ is determined (resp., successively determined) by the $X$-rays in the directions of $U$. 
\item[(f)]
We say that $\mathcal{E}$ is {\em determined} (resp., {\em
  successively determined}) by $m$ (orthogonal) projections on
orthogonal complements of $1$-dimensional subspaces if there is a set $\mathcal{T}$ of $m$
pairwise non-parallel $1$-dimensional subspaces of $\R^d$ such that
$\mathcal{E}$ is determined (resp., successively determined) by the
(orthogonal) projections on the othogonal complements $T^{\perp}$ of the subspaces
$T$ of $\mathcal{T}$. 
\end{itemize}
\end{defi}

\begin{rem}
Let $\mathcal{E}\subset \mathcal{F}(\mathbbm{R}^{d})$. Note that if
$\mathcal{E}$ is determined by a set of $X$-rays (resp., projections), then
$\mathcal{E}$ is successively determined by the same $X$-rays (resp., projections).
\end{rem}

\begin{defi}
Consider the Euclidean plane, $\R^2$.
\begin{itemize}
\item[(a)]
A {\em linear transformation} (resp., {\em affine transformation}) $\Psi\!:\, \mathbbm{R}^{2} \longrightarrow \mathbbm{R}^{2}$ is given by $z \longmapsto Az$ (resp., $z \longmapsto Az+t$), where $A$ is a real $2\times 2$ matrix and $t\in \mathbbm{R}^{2}$. In both cases, $\Psi$ is called {\em singular} (resp., non-singular) when $\operatorname{det}(A)= 0$ (resp., $\operatorname{det}(A)\neq 0$). 
\item[(b)]
A {\em homothety} $h\!:\, \mathbbm{R}^{2} \longrightarrow
\mathbbm{R}^{2}$ is given by $z \longmapsto \lambda z + t$, where
$\lambda \in \R$ is positive and $t\in \mathbbm{R}^{2}$. A homothety
is called a {\em dilatation} if $t=0$. Further, we call a homothety {\em expansive} if $\lambda>1$.
\end{itemize}
\end{defi}

\begin{defi}\label{deficonvex} Let $S \subset \R^2$.
\begin{itemize}
\item[(a)]
A {\em convex polygon} is the convex hull of a finite set of points in $\R^2$. 
\item[(b)]
A {\em polygon in} $S$ is a convex polygon with all vertices in $S$. 
\item[(c)]
A finite subset $C$ of $S$ is called a {\em convex set in} $S$ if its
convex hull contains no new points of $S$, i.e., if $C =
\operatorname{conv}(C)\cap S$. Moreover, the set of all convex sets in
$S$ is denoted by $\mathcal{C}(S)$. 
\item[(d)]
A {\em regular polygon} is always assumed to be planar, non-degenerate and convex. An {\em affinely regular polygon} is a non-singular affine image of a regular polygon. In particular, it must have at least $3$ vertices.
\item[(e)]
Let $U\subset \mathbb{S}^{1}$ be a finite set of directions. We call a non-degenerate convex polygon $P$ a $U${\em -polygon} if it has the property that whenever $v$ is a vertex of $P$ and $u\in U$, the line $\ell_{u}^{v}$ in the plane in direction $u$ which passes through $v$ also meets another vertex $v'$ of $P$.
\end{itemize}    
\end{defi}

\begin{rem}
Note that $U$-polygons have an even number of vertices. Moreover, an affinely regular polygon with an even number of vertices is a $U$-polygon if and only if each direction of $U$ is parallel to one of its edges.     
\end{rem}

The following property is straight-forward.

\begin{lem}\label{homotu}
  Let $h\!:\,\mathbbm{R}^{2}  \longrightarrow \mathbbm{R}^{2}$ be a homothety and let $U\subset \mathbb{S}^{1}$ be a finite set of directions. Then, one has:
\begin{itemize}
\item[(a)]
If $P$ is a $U$-polygon, the image $h[P]$ of $P$ under $h$ is again a $U$-polygon.
\item[(b)]
If $F_{1}$ and $F_{2}$ are elements of $\mathcal{F}(\mathbbm{R}^2)$
with the same $X$-rays in the directions of $U$, the sets $h[F_{1}]$ and $h[F_{2}]$ also have the same $X$-rays in the directions of $U$.\hfill \qed
\end{itemize}
\end{lem}

For all $n \in \mathbbm{N}$, and $\zetan$ a fixed primitive $n$th root
of unity (e.g., $\zetan = e^{2\pi i/n}$), let $\Q(\zetan)$ be
the corresponding cyclotomic field. It is well known that
$\Q(\zetan+\bar{\zeta}_{n})$ is the maximal real subfield of
$\Q(\zetan)$; see~\cite{Wa}. Throughout this paper, we shall use the
notation $$\mathbbm{K}_{n}=\Qn,\;
\mathbbm{k}_{n}=\Q(\zetan+\bar{\zeta}_{n}),\; \OO_{n}=\Zn,\; \oo_{n}
=\Z[\zetan+\bar{\zeta}_{n}]\,,$$ while $\phi$ denotes Euler's phi-function (often also called
Euler's totient function), i.e., $$\phi(n) =
\operatorname{card}\left(\big\{k \in \mathbbm{N}\, |\,1 \leq k \leq n
  \textnormal{ and } (k,n)=1\big\}\right)\,.$$ Occasionally, we shall
identify $\C$ with $\R^{2}$. Moreover, we denote the set of primes by
$\mathbbm{P}$ and its subset of Sophie Germain prime numbers (i.e.,
primes $p$ for which the number $2p+1$ is prime as well) by $\mathbbm{P}_{SG}$. 

\begin{rem}
Sophie Germain prime numbers $p\in\mathbbm{P}_{SG}$ are the primes $p$ such that the equation
$\phi(n)=2p$
has solutions. It is not known whether there are infinitely many
Sophie Germain primes. The first few are
\begin{eqnarray*}&&\{2,3,5,11,23,29,41,53,83,89,113,131,173,\\
  &&\hphantom{\{}179,191,233,239,251,
           281,293,359,419,\dots\}\,,\end{eqnarray*} 
see entry A005384 of~\cite{Sl} for further details.
\end{rem}

\begin{lem}\label{Oo}
For $n \geq 3$, one has: 
\begin{itemize}
\item[(a)]
$\OO_{n}$ is an $\oo_{n}$-module of rank $2$. More precisely, $\OO_{n} =\, \oo_{n} +\, \oo_{n} \,\zetan$, and $\{1,\zetan\}$ is an $\oo_{n}$-basis of $\OO_{n}$.
\item[(b)] 
$\mathbbm{K}_{n}$ is a $\mathbbm{k}_{n}$-vector space of dimension $2$. More precisely, $\mathbbm{K}_{n} =\, \mathbbm{k}_{n} +\, \mathbbm{k}_{n} \,\zetan$, and $\{1,\zetan\}$ is a $\mathbbm{k}_{n}$-basis of $\mathbbm{K}_{n}$.
\end{itemize}
\end{lem}
\begin{proof}
First, we show (a). The linear independence of $\{1,\zetan\}$ over $\oo_{n}$ is clear (since by our assumption $n \geq 3$, $\{1,\zetan\}$ is even linearly independent over $\mathbbm{R}$). It suffices to prove that all non-negative integral powers $\zetan^m$ satisfy $\zetan^m= \alpha + \beta \zetan$ for suitable $\alpha, \beta \in \oo_{n}$. Using induction, we are done if we show  $\zetan^2= \alpha + \beta \zetan$ for suitable $\alpha, \beta \in \oo_{n}$. To this end, note $\bar{\zeta}_{n} = \zetan^{-1}$ and observe that $\zetan^2 = -1 + (\zetan + \zetan^{-1}) \zetan$. Claim (b) follows from the same argument. 
\end{proof}

\begin{rem}\label{dense}
Seen as a point set of $\mathbbm{R}^2$, $\OO_{n}$ has $N$-fold cyclic symmetry, where 
\begin{equation}\label{eq} N=N(n):= [n,2]\,.
\end{equation}
Except for the one-dimensional case ($n\in\{1,2\}$, where
$\OO_{1}=\OO_{2}=\Z$) and the crystallographic cases ($n\in\{3,6\}$
(triangular lattice $\OO_{3}=\OO_{6}$) and $n=4$ (square lattice $\OO_{4}$), see
Figure~\ref{fig:squaretri}), $\OO_{n}$ is dense in $\mathbbm{R}^2$;
see~\cite[Remark 1]{BG2}. 
\end{rem}

\begin{prop}[Gau\ss]\label{gau}
$[\mathbbm{K}_{n} : \mathbbm{Q}] = \phi(n)$, hence the set $\{1,\zetan,\zetan^2,
\dots ,\zetan^{\phi(n)-1}\}$ is a $\mathbbm{Q}$-basis of
$\mathbbm{K}_{n}$. The field extension $\mathbbm{K}_{n}/ \mathbbm{Q}$
is a Galois extension with Abelian Galois group $G(\mathbbm{K}_{n}/
\mathbbm{Q}) \cong (\Z / n\Z)^{\times}$,
where $a\, (\textnormal{mod}\, n)$ corresponds to the automorphism given by\/ $\zetan \longmapsto \zetan^{a}$. 
\end{prop}
\begin{proof}
See~\cite[Theorem 2.5]{Wa}.
\end{proof}

\begin{rem}
Note the identity $(\Z / n\Z)^{\times}=\{a\, (\textnormal{mod}\,
n)\,|\,(a,n)=1\}$ and see~\cite[Table 3]{BG} for examples of the explicit structure of
$G(\mathbbm{K}_{n}/ \mathbbm{Q})$. This reference also contains
further material and references on the role of $\OO_{n}$ in the
context of model sets.
\end{rem}

\begin{lem}[Degree formula]\label{degree}
Let $E/F/K$ be an extension of fields. Then $$[E:K] = [E:F] [F:K]\,.$$
\end{lem}
\begin{proof}
See \cite[Ch. V.1, Proposition 1.2]{La}.
\end{proof}

\begin{coro}\label{cr5}
If $n\geq 3$, one has $[\mathbbm{k}_{n} : \mathbbm{Q}] =
\phi(n)/2$. It follows that a $\mathbbm{Q}$-basis of
$\mathbbm{k}_{n}$ is given by
$\{1,(\zetan+\bar{\zeta}_{n}),(\zetan+\bar{\zeta}_{n})^2, \dots
,(\zetan+\bar{\zeta}_{n})^{\phi(n)/2-1}\}$. Moreover,
$\mathbbm{k}_{n} / \mathbbm{Q}$ is a Galois extension with Abelian
Galois group $G(\mathbbm{k}_{n} / \mathbbm{Q})\cong (\Z /
n\Z)^{\times}/\{\pm 1\, (\textnormal{mod}\, n)\}$.
\end{coro}
\begin{proof}
This is an immediate consequence of Lemma \ref{Oo}(b), Lemma
\ref{degree}, Proposition \ref{gau} and Galois theory.
\end{proof}

\begin{rem}\label{deg2}
The dimension statement of Corollary~\ref{cr5} also follows via Galois theory from Lemma~\ref{degree} in connection with the fact that $\mathbbm{k}_{n}$ is the fixed field of $\mathbbm{K}_{n}$ with respect to the subgroup $\{id, \bar{.}\}$ of $G(\mathbbm{K}_{n}/ \mathbbm{Q})$, where $\bar{.}$ denotes complex conjugation, i.e., the automorphism given by $\zetan \longmapsto \zetan^{-1}$ (recall that $\mathbbm{k}_{n}$ is the maximal real subfield of $\mathbbm{K}_{n}$). 
\end{rem}

\begin{lem}\label{cyclosec}
If $m,n \in \mathbbm{N}$, then $\mathbbm{K}_{m} \cap \mathbbm{K}_{n} = \mathbbm{K}_{(m,n)}$.
\end{lem}
\begin{proof}
The assertion follows from similar arguments as in the proof of the special case $(m,n)=1$; compare \cite[Ch. VI.3, Corollary 3.2]{La}. Here, one has to observe $\mathbbm{Q}(\zeta_{m},\zetan)=\mathbbm{K}_{m}\mathbbm{K}_{n}= \mathbbm{K}_{[m,n]}$ and then to employ the identity \begin{equation}\label{eq2}\phi(m)\phi(n)=\phi([m,n])\phi((m,n))\end{equation} instead of merely using the multiplicativity of the arithmetic function $\phi$. 
\end{proof}

\begin{lem}\label{incl}
Let $m,n \in \mathbbm{N}$. The following statements are equivalent:
\begin{itemize}
\item[(i)]
$\mathbbm{K}_{m} \subset \mathbbm{K}_{n}$.
\item[(ii)]
$m|n$, or $m\; \equiv \;2  \;(\operatorname{mod} 4)$ and $m|2n$.
\end{itemize}
\end{lem}
\begin{proof}
For direction (ii) $\Rightarrow$ (i), suppose first $m|n$. This
clearly implies the assertion. Secondly, if $m\; \equiv \;2
\;(\operatorname{mod} 4)$, say $m=2o$ for a suitable odd number $o$,
and $m|2n$, then $\mathbbm{K}_{o} \subset \mathbbm{K}_{n}$ (due to
$o|n$). However, Proposition~\ref{gau} shows that the inclusion of fields
$\mathbbm{K}_{o} \subset \mathbbm{K}_{2o}=\mathbbm{K}_{m}$ cannot be
proper since we have, by means of the multiplicativity of $\phi$,
$\phi(m)=\phi(2o)=\phi(o)$. This gives $\mathbbm{K}_{m}\subset \mathbbm{K}_{n}$. 

For direction (i) $\Rightarrow$ (ii), suppose the inclusion of fields $\mathbbm{K}_{m} \subset \mathbbm{K}_{n}$. Then, Lemma~\ref{cyclosec} implies $\mathbbm{K}_{m} = \mathbbm{K}_{(m,n)}$, hence \begin{equation}\label{equation}\phi(m)=\phi((m,n))\,;\end{equation} see Proposition~\ref{gau} again. Using the multiplicativity of $\phi$ together with $\phi(p^{j})\,=\,p^{j-1}\,(p-1)$ for $p\in\mathbbm{P}$ and $j\in\mathbbm{N}$, we see that, given the case $(m,n)< m$, (\ref{equation}) can only be fulfilled if $m\; \equiv \;2  \;(\operatorname{mod} 4)$ and $m|2n$. The remaining case $(m,n)= m$ is equivalent to the relation $m|n$.
\end{proof}

The following consequence is immediate.

\begin{coro}\label{unique}
Let $m,n \in \mathbbm{N}$. The following statements are equivalent:
\begin{itemize}
\item[(i)]
$\mathbbm{K}_{m} = \mathbbm{K}_{n}$.
\item[(ii)]
$m=n$, or $m$ is odd and $n=2m$, or $n$ is odd and $m=2n$.\qed
\end{itemize}
\end{coro}

\begin{rem}
Corollary~\ref{unique} implies that, for $m,n\; \not\equiv \;2
\;(\operatorname{mod} 4)$, one has the identity $\mathbbm{K}_{m} = \mathbbm{K}_{n}$ if and only if $m=n$.
\end{rem}

\begin{lem}\label{unique2}
Let $m,n \in \mathbbm{N}$ with $m,n \geq 3$. Then, one has:
\begin{itemize}
\item[(a)]
$\mathbbm{k}_{m}=\mathbbm{k}_{n} \;\Longleftrightarrow\; \mathbbm{K}_{m}=\mathbbm{K}_{n} \mbox{ or\, } m,n \in \{3,4,6\}.$
\item[(b)]
$\mathbbm{k}_{m}\subset\mathbbm{k}_{n} \;\Longleftrightarrow\; \mathbbm{K}_{m} \subset \mathbbm{K}_{n} \mbox{ or\, } m\in \{3,4,6\}.$
\end{itemize}
\end{lem}
\begin{proof}
For claim (a), let us suppose $\mathbbm{k}_{m}=\mathbbm{k}_{n}=:\mathbbm{k}$ first. Then, Proposition~\ref{gau} and Corollary~\ref{cr5} imply that $[\mathbbm{K}_{m}:\mathbbm{k}]=[\mathbbm{K}_{n}:\mathbbm{k}]=2$. Note that $\mathbbm{K}_{m} \cap \mathbbm{K}_{n} = \mathbbm{K}_{(m,n)}$ is a cyclotomic field containing $\mathbbm{k}$. It follows from Lemma~\ref{degree} that either $\mathbbm{K}_{m} \cap \mathbbm{K}_{n}=\mathbbm{K}_{(m,n)}=\mathbbm{K}_{m}=\mathbbm{K}_{n}$ or $\mathbbm{K}_{m} \cap \mathbbm{K}_{n}=\mathbbm{K}_{(m,n)}= \mathbbm{k}$ and hence $\mathbbm{k}_{m}=\mathbbm{k}_{n}= \mathbbm{k}= \mathbbm{Q}$, since the latter is the only real cyclotomic field. Now, this implies $m,n\in\{3,4,6\}$; see also Lemma~\ref{phin2p}(a) below. The other direction is obvious. Claim (b) follows immediately from the part (a).   
\end{proof}

\begin{lem}\label{phin2p}
Consider $\phi$ on $\{n\in \mathbbm{N}\, |\, n\; \not\equiv \;2  \;(\operatorname{mod} 4)\}$. Then, one has: 
\begin{itemize}
\item[(a)] $\phi(n)/2=1$ if and only if $n\in\{3,4\}$. These are the crystallographic cases of the plane.
\item[(b)] $\phi(n)/2 \in \mathbbm{P}$ if and only if $n \in \mathbbm{S}:=\{8,9,12\} \cup \{2p+1\, | \, p \in \mathbbm{P}_{SG}\}$.
\end{itemize}
\end{lem}
\begin{proof}
The equivalences follow from the multiplicativity of $\phi$ in conjunction with the identity $\phi(p^{j})\,=\,p^{j-1}\,(p-1)$ for $p\in\mathbbm{P}$ and $j\in\mathbbm{N}$.
\end{proof}

\begin{rem}
Let $n\; \not\equiv \;2  \;(\operatorname{mod} 4)$. By Corollary~\ref{cr5}, $\mathbbm{k}_{n} / \mathbbm{Q}$ is a Galois extension with Abelian Galois group $G(\mathbbm{k}_{n} / \mathbbm{Q})$ of order $\phi(n)/2$. Now, it follows from Lemma~\ref{phin2p} that $G(\mathbbm{k}_{n} / \mathbbm{Q})$ is trivial if and only if $n\in\{1,3,4\}$, and simple if and only if $n\in\mathbbm{S}$. 
\end{rem}

\begin{prop}\label{p1}
For $n\in \mathbbm{N}$, one has:
\begin{itemize}
\item[(a)]
$\OO_{n}$ is the ring of cyclotomic integers in $\mathbbm{K}_{n}$, and hence its maximal order.
\item[(b)]
$\oo_{n}$ is the ring of integers of $\mathbbm{k}_{n}$, and hence its maximal order.
\end{itemize}
\end{prop}
\begin{proof}
See \cite[Theorem 2.6 and Proposition 2.16]{Wa}.
\end{proof}

\begin{rem}\label{r1}
It follows from Proposition \ref{p1}(a) and Proposition \ref{gau}  that $\OO_{n}$ is a $\Z$-module of rank $\phi(n)$ with $\Z$-basis  $\{1,\zetan,\zetan^2, \dots ,\zetan^{\phi(n)-1}\}$. Likewise, Proposition \ref{p1}(b) and Corollary \ref{cr5} imply that $\oo_{n}$ is a $\Z$-module of rank $\phi(n)/2$ with $\Z$-basis $\{1,(\zetan+\bar{\zeta}_{n}),(\zetan+\bar{\zeta}_{n})^2, \dots ,(\zetan+\bar{\zeta}_{n})^{\phi(n)/2-1}\}$.  
\end{rem}

\begin{defi}
Let $\lambda$ be a real algebraic integer.
\begin{itemize}
\item[(a)]
We call $\lambda$ a {\em Pisot-Vijayaraghavan} number ({\em
  PV-number}) if $\lambda>1$ while all its (algebraic) conjugates have moduli strictly less than $1$.
\item[(b)]
We call $\lambda$ a {\em Pisot-Vijayaraghavan
  unit} ({\em PV-unit}) if $\lambda$ is a PV-number and is also a
unit, i.e., if $1/\lambda$ is an algebraic integer as well.
\end{itemize}
\end{defi}

\begin{lem}\label{pisot}
If $n\,\in\,\mathbbm{N}\setminus \{1,2\}$, then there is a PV-number of $($full$)$ degree $\phi(n)/2$ in $\oo_{n}$.
\end{lem}
\begin{proof}
This follows immediately from~\cite[Ch. I, Theorem 2]{Sa}.
\end{proof}

\begin{rem}\label{r2b}
For $n\,\in\,\mathbbm{N}\setminus \{1,2,3,4,6\}$, there is even a
PV-unit in $\oo_{n}$. This can be seen by considering a representation
of the maximal order $\oo_{n}$ of $\mathbbm{k}_{n}$
(cf. Proposition~\ref{p1}(b)) in logarithmic space together with the
fact that the units of $\oo_{n}$ form a (full) lattice in the
hyperplane given by $$x^{}_{1}+\dots +x_{\frac{\phi(n)}{2}}\, =\, 0\,;$$ see~\cite[Ch. 2, Sections 3 and 4]{Bo}. Clearly, there are  points of this lattice with $x_{1}>0$ and $x_{2},\dots,x_{\phi(n)/2}<0$, and the corresponding $\alpha \in \oo_{n}$ with respect to this representation are PV-units. Note that these PV-units in $\oo_{n}$ necessarily have (full) degree $\phi(n)/2$: if one of them had degree $(\phi(n)/2)/m$, its conjugates (under the effect of the Galois group $G(\mathbbm{k}_{n} / \mathbbm{Q})$, cf. Corollary~\ref{cr5}) would come in sets of $m$ equal conjugates and we could not have just one positive $x$ (this is essentially Dirichlet's unit theorem). We refer the reader also to~\cite[Lemma 8.1.5(b)]{Es}.
\end{rem}

\begin{rem}\label{norm}
Let $\mathbbm{K}/\mathbbm{k}$ be an extension of
algebraic number fields, say of degree $d:=[\mathbbm{K}:\mathbbm{k}]\in\N$. The
corresponding norm $N_{\mathbbm{K}/\mathbbm{k}}\!:\,
\mathbbm{K}\longrightarrow \mathbbm{k}$ is given by
$$
N_{\mathbbm{K}/\mathbbm{k}}(\kappa)=\prod_{j=1}^{d}\sigma_{j}(\kappa)\,,
$$
where the $\sigma_{j}$ are the $d$ distinct embeddings of
$\mathbbm{K}/\mathbbm{k}$ into $\C/\mathbbm{k}$; compare~\cite[Algebraic Supplement, Sec.2,
Corollary 1]{Bo}. In particular, for every $\kappa \in \mathbbm{k}$, one has $N_{\mathbbm{K}/\mathbbm{k}}(\kappa)=\kappa^d$.
Moreover, the norm is transitive in the following sense. If $\mathbbm{L}$ is any intermediate field of $\mathbbm{K}/\mathbbm{k}$
above, then one has
\begin{equation}\label{normtr}
N_{\mathbbm{K}/\mathbbm{k}}=N_{\mathbbm{L}/\mathbbm{k}}\circ N_{\mathbbm{K}/\mathbbm{L}}\,.
\end{equation}
\end{rem}

\begin{lem}\label{extension}
Let $\sigma\!:\, \mathbbm{K}\longrightarrow \mathbbm{K}'$ be an isomorphism of fields, let
$\mathbbm{E}$
be an algebraic extension of $\mathbbm{K}$, and let $\mathbbm{L}$ be an algebraically closed
extension of $\mathbbm{K}'$. Then, there exists a field homomorphism  $\sigma'\;:\;
\mathbbm{E}\longrightarrow \mathbbm{L}$ which extends $\sigma$. 
\end{lem}
\begin{proof}
See \cite[Ch. V.2, Theorem 2.8]{La}.
\end{proof}

From now on, for $n\in
\mathbbm{N}$, we always let $\zetan := e^{2\pi i/n}$, a
primitive $n$th root of unity. 

\section{A Cyclotomic Theorem}\label{sec9}

In this section, we need the following facts from the theory of
$p$-adic valuations; compare~\cite{Gou,ko}. 

Let $p\in \mathbbm{P}$. The $p$-adic valuation on $\Z$ is the
function $v_{p}$, defined by $v_{p}(0):=\infty$ and by the equation
$$
n=p^{v_{p}(n)}n'
$$
for $n\neq 0$, where $p$ does not divide $n'$; that is, $v_{p}(n)$ is
the exponent of the biggest power of $p$ that divides $n$. The function
$v_{p}$ is extended to $\Q$ by defining
$$
v_{p}\Big(\frac{a}{b}\Big):= v_{p}(a)-v_{p}(b)
$$
for $a,b\in\Z\setminus\{0\}$; see~\cite[p. 23]{Gou}. Note that $v_{p}$
is $\Z$-valued on $\Q\setminus\{0\}$. As in~\cite[Ch. 5]{Gou},
$v_{p}$ can further be extended to the algebraic closure $\Q_{p}^{\operatorname{alg}}$
of a field $\Q_{p}$, whose elements are called $p$-{\em adic numbers},
containing $\Q$. Note that $\Q_{p}^{\operatorname{alg}}$ contains the algebraic
closure $\Q^{\operatorname{alg}}$ of $\Q$ and hence all algebraic numbers. On
$\Q_{p}^{\operatorname{alg}}\setminus\{0\}$, $v_{p}$ takes values in $\Q$, and
satisfies
\begin{eqnarray}
v_{p}(-x)&=&v_{p}(x)\,,\label{minuspreserve}\\
v_{p}(xy)&=&v_{p}(x)+v_{p}(y)\,,\label{log1}\\
v_{p}\Big(\frac{x}{y}\Big)&=&v_{p}(x)-v_{p}(y)\label{log2}
\end{eqnarray} 
and
\begin{eqnarray}
v_{p}(x+y)&\geq& \operatorname{min}\{v_{p}(x),v_{p}(y)\}\,,
\end{eqnarray} 
compare also~\cite[p. 143]{Gou}. 

\begin{prop}\label{ppuou} 
Let $p\in \mathbbm{P}$ and let $r,s,t \in \mathbbm{N}$.  
If $r$ is not
a power of $p$ and $(r,s)=1$, one has 
\begin{equation}
v_{p}(1-\zeta_{r}^{s})=0\,.
\end{equation} 
Otherwise, if $(p,s)=1$, then
\begin{equation}
v_{p}(1-\zeta_{p^{t}}^{s})=\frac{1}{p^{t-1}(p-1)}\,.
\end{equation}
\end{prop}
\begin{proof}
See \cite[Proposition 3.6]{GG}.
\end{proof}

\begin{defi}
Let $k,m\in\mathbbm{N}$ and let $p\in\mathbbm{P}$. An $m$th root of unity $\zeta_{m}^{k}$ is called a $p${\em -power root
  of unity} if there is a $t\in\mathbbm{N}$ such that
$\frac{k}{m}=\frac{s}{p^t}$ for some $s\in\mathbbm{N}$
  with $(p,s)=1$. 
\end{defi}

\begin{rem}
Note that an $m$th root of unity $\zeta_{m}^{k}$ is a $p$-power root
  of unity if and only if it is a primitive $p^{t}$th root of unity
  for some $t\in\mathbbm{N}$.
\end{rem}

\begin{lem}\label{ppuou2}
Let $k,t \in\mathbbm{N}$. Further, let $j,m \in\mathbbm{N}$ with
$(j,m)=1$ and let
$p\in\mathbbm{P}$. Then, one has:
\begin{itemize}
\item[(a)]
$\zeta_{m}^{k}$ is a primitive $p^t$th root of unity
 if and only if $(\zeta_{m}^{j})^{k}$ is a
 primitive $p^t$th root of unity.
\item[(b)]
$\zeta_{m}^{k}$ is a $p$-power root of unity
 if and only if $(\zeta_{m}^{j})^{k}$ is a
$p$-power root of unity.
\end{itemize}
 \end{lem}
\begin{proof}
We first prove claim (a). Assume that $\frac{k}{m}=\frac{s}{p^t}$ for a suitable
$s\in\mathbbm{N}$ with $(p,s)=1$. In particular, it follows that $p|m$ and, since $(j,m)=1$,
one has $(p,j)=1$. Hence,
$\frac{jk}{m}=\frac{js}{p^t}$ and $(p,js)=1$. Conversely, assume $\frac{jk}{m}=\frac{s}{p^t}$ for a suitable
$s\in\mathbbm{N}$ with $(p,s)=1$. Since $(j,m)=1$, it follows that
$j|s$, say $jl=s$ for a suitable $l\in\mathbbm{N}$. Hence, $\frac{k}{m}=\frac{l}{p^t}$
and, moreover, $(p,l)=1$. Claim (b) follows immediately from the part (a).
\end{proof}

\begin{lem}\label{sigmapreserve0}
Let $m,k\in\mathbbm{N}$ and let $p\in\mathbbm{P}$. If $\sigma\in G(\mathbbm{K}_{m}/\Q)$, then
$$
v_{p}(1-\zeta_{m}^{k})=v_{p}\big(\sigma(1-\zeta_{m}^{k})\big)\,.
$$
\end{lem}
\begin{proof}
By Proposition~\ref{gau}, $\sigma$ is given by $\zeta_{m}\longmapsto
\zeta_{m}^{j}$, where $j\in\mathbbm{N}$ satisfies $(j,m)=1$. The
assertion follows immediately from Proposition~\ref{ppuou} in
conjunction with Lemma~\ref{ppuou2}.
\end{proof}

\begin{rem}
Note that Lemma~\ref{sigmapreserve0} is only one instance of the following more
general fact. It is well known that if $\mathbbm{K}$ is a normal algebraic number
field (i.e., a finite Galois extension of $\Q$) and if
$\kappa\in\mathbbm{K}$, then one has $v_p(\kappa)=v_p(\sigma(\kappa))$
for every $p\in\mathbbm{P}$ and every $\sigma\in G(\mathbbm{K}/\Q)$;
see~\cite[Ch. 3, Sec. 4, Prob. 7]{Bo} and the parenthetical clause of
the sentence following~\cite[Proposition 2.14]{Wa} in conjunction
with~\cite[Theorem 2.13]{Wa} and~\cite[Proposition 2.14]{Wa} itself.
\end{rem}

\begin{defi}
Let $m\geq 4$ be a natural number. We define 
$$
D'_{m}:=\{(k_1,k_2,k_3,k_4)\in \mathbbm{N}^4 \,|\,
k_1,k_2,k_3,k_4\leq m-1 \mbox{ and } k_1+k_2=k_3+k_4\}\,,
$$
together with its subset
$$
D_{m}:=\{(k_1,k_2,k_3,k_4)\in \mathbbm{N}^4 \,|\, k_3<k_1\leq
k_2<k_4\leq m-1 \mbox{ and } k_1+k_2=k_3+k_4\}\,,
$$
and define the function $f_{m}\,:\, D'_{m}\longrightarrow \C$ by
\begin{equation}\label{fmd}
f_{m}\big((k_1,k_2,k_3,k_4)\big):=\frac{(1-\zeta_{m}^{k_1})(1-\zeta_{m}^{k_2})}{(1-\zeta_{m}^{k_3})(1-\zeta_{m}^{k_4})}\,.
\end{equation}
\end{defi}

\begin{lem}\label{fmdg1} 
Let $m\geq 4$. The function
$f_{m}$ takes only positive values. Moreover, one has $f_{m}(d)>1$ for all $d\in
D_{m}$.
\end{lem}
\begin{proof}
See the proof of~\cite[Lemma 3.1]{GG}.
\end{proof}

\begin{coro}\label{fmdkm}
Let $m\geq 4$ and let $d\in
D'_{m}$. Then, one has $f_{m}(d)\in\mathbbm{k}_{m}$. 
\end{coro}
\begin{proof}
By Lemma~\ref{fmdg1}, $f_{m}(d)$ is real. The assertion follows
immediately from the fact that $f_{m}(d)\in\mathbbm{K}_{m}$
and the fact that $\mathbbm{k}_{m}$ is the maximal real subfield of
$\mathbbm{K}_{m}$.
\end{proof}

For our application to discrete tomography, we shall investigate the
set
\begin{equation}\label{eqn13}
\Big ( \bigcup_{m\geq 4} f_{m}[D_{m}]\Big )\cap \mathbbm{k}
\end{equation}
in the case of arbitrary real algebraic number fields $\mathbbm{k}$. Gardner and Gritzmann showed the
following result, which deals with the smallest among all algebraic
number fields, i.e., with $\mathbbm{k}=\Q$.

\begin{theorem}\label{intersectq}
$$
\Big ( \bigcup_{m\geq 4} f_{m}[D_{m}]\Big )\cap \Q=\Big\{\frac{4}{3},\frac{3}{2},2,3,4\Big\}=:N_1\,.
$$
Moreover, all solutions of $f_{m}(d)=q\in\Q$, where $m\geq 4$ and
$d:=(k_1,k_2,k_3,k_4)\in D_{m}$, are either given, up to multiplication of $m$
and $d$ by the same factor, by $m=12$ and one of the following
$$
\begin{array}{rlrl}
\textnormal{(i)}&d=(6,6,4,8),q=\frac{4}{3};&\textnormal{(ii)}&d=(6,6,2,10),q=4;\\
\textnormal{(iii)}&d=(4,8,3,9),q=\frac{3}{2};&\textnormal{(iv)}&d=(4,8,2,10),q=3;\\
\textnormal{(v)}&d=(4,4,2,6),q=\frac{3}{2};&\textnormal{(vi)}&d=(8,8,6,10),q=\frac{3}{2};\\
\textnormal{(vii)}&d=(4,4,1,7),q=3;&\textnormal{(viii)}&d=(8,8,5,11),q=3;\\
\textnormal{(ix)}&d=(3,9,2,10),q=2;&\textnormal{(x)}&d=(3,3,1,5),q=2;\\
\textnormal{(xi)}&d=(9,9,7,11),q=2;&&
\end{array}
$$
or by one of the following
$$
\begin{array}{rl}
\textnormal{(xii)}&d=(2k,s,k,k+s),q=2, \mbox{ where } s\geq 2, m=2s \mbox{
  and } 1\leq k\leq \frac{s}{2};\\
\textnormal{(xiii)}&d=(s,2k,k,k+s),q=2, \mbox{ where } s\geq 2, m=2s \mbox{
  and } \frac{s}{2}\leq k< s.
\end{array}
$$\end{theorem}
\begin{proof}
See \cite[Lemma 3.8, Lemma 3.9 and Theorem 3.10]{GG}.
\end{proof}

In this section, we shall not be able to present full analogues of
Theorem~\ref{intersectq} in the case of 
the set~\eqref{eqn13} for real algebraic number fields
$\mathbbm{k}\neq\Q$. (The study of this interesting problem is work in
progress.) Instead, we shall show that the elements of the
set~\eqref{eqn13} satisfy certain conditions.  

In fact, although we are only interested in the set~(\ref{eqn13}), we
shall investigate its superset 
\begin{equation}\label{eqn}
\Big ( \bigcup_{m\geq 4} f_{m}[D'_{m}]\Big )\cap \mathbbm{k}
\end{equation}
in the case of arbitrary real algebraic number fields $\mathbbm{k}$. Every (real) algebraic number field $\mathbbm{k}$ different from $\Q$
has some finite
degree $e:=[\mathbbm{k}:\Q]\geq 2$. In particular, we shall be
interested in the case $e=2$. In the
following, we shall denote certain degrees by $e$ (resp., $f$). We
wish to emphasize that this should not be confused with their usual
usage, i.e., where $e$ denotes the ramification index and $f$ denotes the residue class degree.

\begin{lem}\label{sigmapreserve}
Let $m\geq 4$, let $p\in \mathbbm{P}$, and let $d\in
D'_{m}$. If $\sigma\in G(\mathbbm{K}_{m}/\Q)$, then
$$
v_{p}\big(f_{m}(d)\big)=v_{p}\Big(\sigma\big(f_{m}(d)\big)\Big)\,.
$$
\end{lem}
\begin{proof}
The
assertion follows immediately from Lemma~\ref{sigmapreserve0} and the Equations~(\ref{log1})
and~(\ref{log2}).
\end{proof}

\begin{lem}\label{1dn}
Let $m\geq 4$ and let
$d\in D'_{m}$. Then, for any prime
factor $p\in\mathbbm{P}$ of the numerator of the field norm 
$N_{\Q(f_m(d))/\Q}(f_m(d))$, one has 
\begin{equation*}\label{vpn}
v_{p}\Big(N_{\Q(f_{m}(d))/\Q}\big(f_{m}(d)\big)\Big)=e\,v_{p}\big(f_{m}(d)\big)\in \mathbbm{N}\,,
\end{equation*}
where $e:=[\Q(f_m(d)):\Q]\in\N$ is the degree of $f_{m}(d)$ over $\Q$. 
\end{lem}
\begin{proof}
By Corollary~\ref{fmdkm}, one has 
$f_{m}(d)\in\mathbbm{k}_{m}\subset \mathbbm{K}_{m}$. It
follows the inclusion of fields 
\begin{equation}\label{incl10}
\Q\big(f_{m}(d)\big)\,\subset\,\mathbbm{k}_{m}\,\subset\,
\mathbbm{K}_{m}\,.
\end{equation} 
The norm 
$N_{\Q(f_{m}(d))/\Q}\!:\,\Q\big(f_{m}(d)\big)\longrightarrow \Q$ of the field
extension $\Q\big(f_{m}(d)\big)/\Q$ is given by
$$N_{\Q(f_{m}(d))/\Q}(q)=\prod_{j=1}^{e}\sigma_{j}(q)$$ for $q\in
\Q\big(f_{m}(d)\big)$, where
$\{\sigma_{j}\,|\,j\in\{1,\dots,e\}\}$ is the Galois group
$G\big(\Q\big(f_{m}(d)\big)/\Q\big)$; see Remark~\ref{norm} and note that
the field extension $\Q\big(f_{m}(d)\big)/\Q$ is indeed a Galois
extension. The latter follows immediately from Galois theory,
  since, by Proposition~\ref{gau}, the Galois extension
  $\mathbbm{K}_{m}/\Q$ has an Abelian Galois group; cf. Proposition~\ref{gau}. By
Relation~(\ref{incl10}), Lemma~\ref{extension}, and since $\mathbbm{K}_{m}/\Q$ is a Galois
extension,
 each field automorphism $\sigma_{j}\in G\big(\Q\big(f_{m}(d)\big)/\Q\big)$, $j\in\{1,\dots,e\}$, can be extended to
a field automorphism $\sigma'_{j}\in G(\mathbbm{K}_{m}/\Q)$. It
follows that
$$N_{\Q(f_{m}(d))/\Q}\big(f_{m}(d)\big)=\prod_{j=1}^{e}\sigma'_{j}\big(f_{m}(d)\big)\,.$$ Using the $p$-adic valuation $v_{p}$ in conjunction with
 Equation~(\ref{log1}) and Lemma~\ref{sigmapreserve}, one gets
\begin{equation*}
v_{p}\Big(N_{\Q(f_{m}(d))/\Q}\big(f_{m}(d)\big)\Big)=e\, v_{p}\big(f_{m}(d)\big)\in\mathbbm{N}\,,
\end{equation*}
which completes the proof.
\end{proof}

\begin{lem}\label{2dn}
Let $m\geq 4$ and let
$d\in D'_{m}$. Then, for any prime
factor $p\in\mathbbm{P}$ of the numerator of 
$N_{\Q(f_m(d))/\Q}(f_m(d))$, the expression $\zeta_{m}^{k_{j}}$ in~$(\ref{fmd})$ is
a $p$-power root of unity for at least one value of
$j\in\{1,2\}$.
\end{lem}
\begin{proof}
Let $e:=[\Q(f_m(d)):\Q]\in\N$ be the degree of $f_{m}(d)$ over $\Q$.  
By assumption, the numerator of $N_{\Q(f_{m}(d))/\Q}(f_{m}(d))$ has a
prime factor, say $p\in\mathbbm{P}$.
Using Lemma~\ref{1dn}, one gets
$$
v_{p}\big(f_{m}(d)\big)=v_{p}\left ( \frac{(1-\zeta_{m}^{k_1})(1-\zeta_{m}^{k_2})}{(1-\zeta_{m}^{k_3})(1-\zeta_{m}^{k_4})}
\right )\in \frac{1}{e}\mathbbm{N}\,.
$$ 
The assertion follows immeiately from Equation~(\ref{log1}), Equation~(\ref{log2})
and Proposition~\ref{ppuou}.
\end{proof}

The following result is a rather simple one. For convenience, we include a proof.

\begin{lem}\label{3dn}
Let $m\geq 4$ and let
$d:=(k_1,k_2,k_3,k_4)\in D'_{m}$. Then one has the equality of fields 
$\Q(f_{m}(d))=\Q(f_{m}(d'))$ and, moreover, 
$$
N_{\Q(f_{m}(d'))/\Q}\big(f_{m}(d')\big)=\Big(N_{\Q(f_{m}(d))/\Q}\big(f_{m}(d)\big)\Big)^{-1}\,,
$$
where $d':=(k_3,k_4,k_1,k_2)\in D'_{m}$. In particular, $f_{m}(d)$ and
$f_{m}(d')$ have the same degree over $\Q$.
\end{lem}
\begin{proof}
Clearly, one has $d':=(k_3,k_4,k_1,k_2)\in D'_{m}$. Futher, observing
 the identity $f_{m}(d')=1/f_{m}(d)$ (by Lemma~\ref{fmdg1}, one has
 $f_{m}(d)\neq 0$), one sees the equality of the fields
 $\Q(f_{m}(d))=\Q(f_{m}(d'))$. The identity
\begin{eqnarray*}\label{hom}
N_{\Q(f_{m}(d'))/\Q}\big(f_{m}(d')\big)&=&N_{\Q(f_{m}(d))/\Q}\big(f_{m}(d')\big)\\
&=&N_{\Q(f_{m}(d))/\Q}\Big(\big(f_{m}(d)\big)^{-1}\Big)\\
&=&\Big(N_{\Q(f_{m}(d))/\Q}\big(f_{m}(d)\big)\Big)^{-1}
\end{eqnarray*} 
completes the proof.
\end{proof}

The following definition will be useful.

\begin{defi}
Let $a\in\N$, let $S\subset \R\setminus\{0\}$, and let $\alpha\in\Z$. We set:
\begin{itemize}
\item[(a)]
 $a^{\downarrow}:=\{\frac{a}{b}\in\Q\,|\,b\in\N\mbox{ with } b<a \mbox{
    and } (a,b)=1\}$.
\smallskip
\item[(b)]
$[S]^{\alpha}:=\{x^{\alpha}\,|\,x\in S\}$.
\end{itemize}
\end{defi}

\subsection{The Case of Degree Two}

\begin{lem}\label{lem1}
Let $m\geq 4$ and let
$d:=(k_1,k_2,k_3,k_4)\in D'_{m}$. Suppose that
$f_{m}(d)$ is of degree two over $\Q$ and suppose that the absolute value of
$N_{\Q(f_{m}(d))/\Q}(f_{m}(d))$ is greater than $1$. Then, one has:
$$
N_{\Q(f_{m}(d))/\Q}\big(f_{m}(d)\big)\in
\bigcup_{a\in\{2,3,4,5,6,8,9,12,16\}}\pm \, a^{\downarrow}\,\,\,\,\,\,=:N \,.
$$
\end{lem}
\begin{proof}
By assumption, the
numerator of $N_{\Q(f_{m}(d))/\Q}(f_{m}(d))$ has a prime factor
$p\in\mathbbm{P}$. Further, by Lemma~\ref{1dn}, for every such prime factor
$p\in\mathbbm{P}$, one has 
\begin{equation}\label{restr0}
v_{p}\Big(N_{\Q(f_{m}(d))/\Q}\big(f_{m}(d)\big)\Big)= 2\,v_{p}\big(f_{m}(d)\big)=2\,v_{p}\left ( \frac{(1-\zeta_{m}^{k_1})(1-\zeta_{m}^{k_2})}{(1-\zeta_{m}^{k_3})(1-\zeta_{m}^{k_4})}
\right )\in \mathbbm{N}\,.
\end{equation}
Applying Equations~\eqref{log1} and
\eqref{log2}, Proposition~\ref{ppuou} and Lemma~\ref{2dn},
one sees that $v_{p}(f_{m}(d))$ is a sum of at most four terms of the
form $1/(p^{t'-1}(p-1))$ for various $t'\in\N$ with one or two
positive terms and at most two negative ones. Let $t$ be the smallest
$t'$ occuring in one of the positive terms. Then,
Relation~\eqref{restr0} particularly shows that   
$$
\frac{4}{p^{t-1}(p-1)}\geq 1
$$
or, equivalently,
\begin{equation}\label{restr1}
p^{t-1}(p-1)\leq 4\,.
\end{equation}
One can see easily that the only possibilities for~\eqref{restr1} are
$p=2$ and $t\in\{1,2,3\}$, or $p=3$ and $t=1$, or $p=5$ and
$t=1$. Moreover, using~\eqref{restr0} in conjunction with Equations~\eqref{log1} and
\eqref{log2} and Proposition~\ref{ppuou}, one can see the following. In the first
case ($p=2$ and $t\in\{1,2,3\}$), the only possibilities for powers of $2$ dividing the numerator
of $N_{\Q(f_{m}(d))/\Q}(f_{m}(d))$ are $2,4,8$ and $16$, with room for
a $3$ as well in the cases $2$ and $4$, i.e., in the case where the
$2$-adic value of 
$N_{\Q(f_{m}(d))/\Q}(f_{m}(d))$ is $1$ or $2$. Hence, one obtains   
$N_{\Q(f_{m}(d))/\Q}(f_{m}(d))\in
\cup_{a\in\{2,4,6,8,12,16\}}\pm \, a^{\downarrow}$. In the second 
case ($p=3$ and $t=1$), the only possibilities for powers of $3$ dividing the numerator
of $N_{\Q(f_{m}(d))/\Q}(f_{m}(d))$ are $3$ and $9$, with room for
a $2$ or $4$ as well in the case $3$, i.e., in the case where the
$3$-adic value of 
$N_{\Q(f_{m}(d))/\Q}(f_{m}(d))$ is $1$. Consequently, one obtains $N_{\Q(f_{m}(d))/\Q}(f_{m}(d))\in\cup_{a\in\{3,6,9,12\}}\pm \, a^{\downarrow}
$, 
whereas the last case ($p=5$ and
$t=1$) immediately implies that $N_{\Q(f_{m}(d))/\Q}(f_{m}(d))\in 5^{\downarrow}$.
\end{proof}

\begin{lem}\label{lem3}
Let $m\geq 4$ and let
$d:=(k_1,k_2,k_3,k_4)\in D'_{m}$. Suppose that
$f_{m}(d)$ is of degree two over $\Q$ and suppose that  the absolute
value of $N_{\Q(f_{m}(d))/\Q}(f_{m}(d))$ is smaller than $1$. Then, one has:
$$
N_{\Q(f_{m}(d))/\Q}(f_{m}(d))\in [N]^{-1}\,,
$$
with $N$ as defined in Lemma~$\ref{lem1}$.
\end{lem}
\begin{proof}
Note that, by Lemma~\ref{fmdg1}, one has $f_{m}(d)> 0$, which implies that the absolute value of $N_{\Q(f_{m}(d))/\Q}(f_{m}(d))$ is non-zero. The assertion now follows immediately from Lemma~\ref{3dn} and Lemma~\ref{lem1}.
\end{proof}

\begin{theorem}\label{intersectk8}
For any real quadratic algebraic number field $\kk$, one has:
$$
N_{\kk /\Q}\Big [ \big(\bigcup_{m\geq 4} f_{m}[D_{m}]\big )\cap 
 \kk \Big ] 
\subset\{-1,1\}\cup [N_1]^2\cup N\cup [N]^{-1}=:N_{2}\,,
$$
with $N_1$ as defined in Theorem~$\ref{intersectq}$ and $N$ as defined in Lemma~$\ref{lem1}$.
\end{theorem}
\begin{proof}
Let $f_{m}(d)\in (\cup_{m\geq 4} f_{m}[D_{m}])\cap 
 \kk $ for suitable $m\geq 4$ and  $d\in
D_{m}$. Note that, by Lemma~\ref{fmdg1}, one has $f_{m}(d)> 1$
and hence $f_{m}(d)\neq 0$. First, suppose that $f_{m}(d)$ is of degree two
 over $\Q$, i.e., a primitive element of the field extension
 $\kk/\Q$. Since $f_{m}(d)\neq 0$, the norm
 $N_{\Q(f_{m}(d))/\Q}(f_{m}(d))=N_{\kk/\Q}(f_{m}(d))$ is non-zero. Therefore,
 the absolute value of $N_{\kk/\Q}(f_{m}(d))$ is greater
 than zero.  The following is divided into the three possible
 cases. First, assume that the absolute value of
 $N_{\kk/\Q}(f_{m}(d))$ is greater than $1$. Then, Lemma~\ref{lem1}
 immediately implies that $N_{\kk/\Q}(f_{m}(d))\in N$. Secondly,
 suppose that the absolute value of 
$N_{\kk/\Q}(f_{m}(d))$ is smaller than $1$. Then, Lemma~\ref{lem3}
immediately implies that $N_{\kk/\Q}(f_{m}(d))\in
[N]^{-1}$. Finally, if the absolute value of $N_{\kk/\Q}(f_{m}(d))$
equals $1$, one obviously gets
$N_{\kk/\Q}(f_{m}(d))\in\{1,-1\}$. Now, suppose that
$f_{m}(d)\in\Q$. Then, by Theorem~\ref{intersectq}, one has
$f_{m}(d)\in N_1$ and hence
$N_{\kk/\Q}(f_{m}(d))\in [N_1]^2$; cf. Remark~\ref{norm}. This completes the proof.
\end{proof}

\subsection{The General Case}

\begin{theorem}\label{intersectk8general}
For any $e\in\N$, there is a finite set $N_{e}\subset \Q$
such that, for all real algebraic number fields $\kk$ having degree $e$, one has:
$$
N_{\mathbbm{k}/\Q}\Big [ \big(\bigcup_{m\geq 4} f_{m}[D_{m}]\big )\cap 
 \mathbbm{k}\Big ] 
\subset N_{e}\,.
$$
\end{theorem}
\begin{proof}
It suffices to show the corresponding assertion for the superset $D'_{m}$ of $D_{m}$. Let $e\in\N$ and let $\kk$ be a real algebraic number fields $\kk$ of degree $e$. Further, let $f_{m}(d)\in (\cup_{m\geq 4} f_{m}[D'_{m}])\cap 
 \mathbbm{k}$ for suitable $m\geq 4$ and  $d\in
D_{m}$. By Lemma~\ref{degree}, and since one has
$$
\Q\,\subset \,\Q\big(f_{m}(d)\big)\,\subset \,\mathbbm{k}\,,$$ $f_{m}(d)$ is of degree $f:=[\Q\big(f_{m}(d)\big):\Q]$ over $\Q$, where $f$
is a divisor of $e$. By Lemma~\ref{fmdg1}, one has $f_{m}(d)\neq
0$. So,  
 the norm $N_{\Q(f_{m}(d))/\Q}(f_{m}(d))$ is non-zero, hence its absolute value is greater
 than zero. Suppose that the absolute value of
 $N_{\Q(f_{m}(d))/\Q}(f_{m}(d))$ is greater than $1$. Then, the numerator of
$N_{\Q(f_{m}(d))/\Q}(f_{m}(d))$ has a prime factor, say
$p\in\mathbbm{P}$. Further, by Lemma~\ref{1dn}, for every such prime factor
$p\in\mathbbm{P}$, one has
\begin{equation}\label{restr17}
v_{p}\Big(N_{\Q(f_{m}(d))/\Q}\big(f_{m}(d)\big)\Big)= e\,v_{p}\big(f_{m}(d)\big)=e\,v_{p}\left ( \frac{(1-\zeta_{m}^{k_1})(1-\zeta_{m}^{k_2})}{(1-\zeta_{m}^{k_3})(1-\zeta_{m}^{k_4})}
\right )\in \mathbbm{N}\,.
\end{equation}
Applying Equations~\eqref{log1} and
\eqref{log2}, Proposition~\ref{ppuou} and Lemma~\ref{2dn},
one sees that $v_{p}(f_{m}(d))$ is a sum of at most four terms of the
form $1/(p^{t'-1}(p-1))$ for various $t'\in\N$ with one or two
positive terms and at most two negative ones. Let $t$ be the smallest
$t'$ occuring in one of the positive terms. Then,
Relation~\eqref{restr17} particularly shows that   
$$
\frac{2\,e}{p^{t-1}(p-1)}\geq 1
$$
or, equivalently,
\begin{equation}\label{restr171}
p^{t-1}(p-1)\leq 2\,e\,.
\end{equation}
Similar to the argumentation in Lemma~\ref{lem1} above, one can now see that, by
Relations~\eqref{restr17} and~\eqref{restr171} and the obvious fact that
$$
p^{t-1}(p-1)\longrightarrow \infty
$$
for fixed $p\in\mathbbm{P}$ as $t\longrightarrow\infty$ (resp., for
fixed $t\in\N$ as $\mathbbm{P}\owns p\longrightarrow\infty$), there is a finite set, say $N_{f}$,
such that $N_{\Q(f_{m}(d))/\Q}(f_{m}(d))\in N_{f}$. Moreover, analogous to the argumentation in
Lemma~\ref{lem3} above, one can see that, if the absolute value of 
$N_{\Q(f_{m}(d))/\Q}(f_{m}(d))$ is smaller than one, then one has
$$N_{\Q(f_{m}(d))/\Q}(f_{m}(d))\in[N_{f}]^{-1}\,,$$ whereas the
missing case, where the absolute value
of $N_{\Q(f_{m}(d))/\Q}(f_{m}(d))$ equals $1$, only gives $N_{\Q(f_{m}(d))/\Q}(f_{m}(d))\in
\{1,-1\}$. By the transitivity of the norm (cf. Remark~\ref{norm}), it
follows that $N_{\mathbbm{k}/\Q}(f_{m}(d))\in
\big[\{1,-1\}\cup N_{f}\cup[N_{f}]^{-1}\big]^{\frac{e}{f}}$. Setting
$$
N_{e}:=\bigcup_{f|e}\big[\{1,-1\}\cup
N_{f}\cup[N_{f}]^{-1}\big]^{\frac{e}{f}}\,,
$$
the assertion immediately follows, since the above analysis only depends on the fixed
 degree $e$. 
\end{proof}

\begin{defi}
We denote by $\Q_{>1}^{\mathbbm{P}_{\leq 2}}$ the subset of $\Q$
consisting of all rational
numbers greater than one, which have at most two prime factors in the numerator.\end{defi}

\begin{theorem}\label{intersectk8generalcoro}
For any real algebraic number field $\kk$, one has:
$$
N_{\mathbbm{k}/\Q}\Big [ \big(\bigcup_{m\geq 4} f_{m}[D_{m}]\big )\cap 
 \mathbbm{k}\Big ] 
\subset \pm\left(\{1\}\cup \big[\Q_{>1}^{\mathbbm{P}_{\leq 2}}\big]^{\pm 1}\right)\,.
$$
\end{theorem}
\begin{proof}
This follows from a careful analysis of the proof
of Theorem~\ref{intersectk8general}.
\end{proof}

\subsection{Applications}

For our later application to discrete tomography of aperiodic model sets, we formulate some
corollaries. We shall be interested in applications of
the above results to real algebraic number fields of the form 
$\mathbbm{k}_{n}$, in particular for $n\in\{5,8,10,12\}$, in view of
their practical relevance; see~\cite{St}. 

\begin{rem}
Note that $\mathbbm{k}_{10}=\mathbbm{k}_{5}=\Q(\sqrt{5})$, 
 $\mathbbm{k}_{8}=\Q(\sqrt{2})$ 
and $\mathbbm{k}_{12}=\Q(\sqrt{3})$. By
Corollary~\ref{cr5}, all of these fields are quadratic algebraic
 number fields. 
\end{rem}

\begin{coro}\label{intersectk510s}
For $n\in\{5,8,10,12\}$, one has:
$$
N_{\mathbbm{k}_{n}/\Q}\Big [ \big(\bigcup_{m\geq 4} f_{m}[D_{m}]\big )\cap 
 \mathbbm{k}_{n}\Big ] 
\subset N_{2}\,,$$
with $N_{2}$ as defined in
Theorem~$\ref{intersectk8}$.
\end{coro}
\begin{proof}
This is an immediate consequence of Theorem~\ref{intersectk8}.
\end{proof}

\begin{defi}
Let the map $\phi/2\!:\,\N\setminus\{1,2\} \longrightarrow
\N$ be given  by $n\longmapsto \phi(n)/2$.
\end{defi}

\begin{rem}
Recall that, for $n\in\N\setminus\{1,2\}$, $\phi/2(n)$ is the degree of
$\kk_{n}$ over $\Q$; cf. Corollary~\ref{cr5}.
\end{rem}

\begin{coro}\label{intersectk510sg}
For all $e\in\operatorname{Im}(\phi/2)$, there is a finite set $N_{e}\subset \Q$
such that, for all $n\in(\phi/2)^{-1}[\{e\}]$, one has:
$$
N_{\mathbbm{k}_{n}/\Q}\Big [ \big(\bigcup_{m\geq 4} f_{m}[D_{m}]\big )\cap 
 \mathbbm{k}_{n}\Big ] 
\subset N_{e}\,.$$
\end{coro}
\begin{proof}
This is an immediate consequence of Theorem~\ref{intersectk8general}.
\end{proof}

\section{Cyclotomic Model Sets}\label{deficyclo}

In the present paper, we shall study the discrete tomography of a special
class of planar model sets, the \emph{cyclotomic} model
sets, which can be described in algebraic terms and have an Euclidean
internal space (defined below).

By definition, {\em model sets} arise from so-called {\em cut and
  project schemes}, compare \cite{Moody}. Before we can present the
  cut and project scheme from which {\em cyclotomic} model sets arise,
  we need the following remark. In the following let $n\in \N
  \setminus \{1,2\}$ (this
means that $\phi(n)\neq 1$).

\begin{rem}\label{mink}
The elements of the Galois group $G(\mathbbm{K}_{n}/
\mathbbm{Q})$ come in pairs of complex conjugate automorphisms. Let the
set $\{\sigma_{1},\dots,\sigma_{\phi(n)/2}\}$ arise from
$G(\mathbbm{K}_{n}/ \mathbbm{Q})$ by choosing exactly one automorphism
from each such pair. Here, we always choose $\sigma_{1}$ as the
identity rather than the complex conjugation. Every such choice induces a map 
$$
\,\,.\,\widetilde{\hphantom{a}}^{_{n}}\,:\, \OO_{n}\longrightarrow
\R^{2}\times(\R^2)^{\frac{\phi(n)}{2}-1}\,,
$$
 given by
$$
z\longmapsto
\left(z,\big(\sigma_{2}(z),\dots,\sigma_{\frac{\phi(n)}{2}}(z)\big )\right)\,.
$$
With the understanding that for $\phi(n)=2$ (this means $n\in\{3,4,6\}$), the singleton $$(\R^2)^{\phi(n)/2-1}=(\R^2)^{0}$$
is the trivial (locally compact) Abelian group
$\{0\}$, each such choice induces, via projection on the second factor,
a map 
$$.^{\star_{n}}\! : \, \OO_{n} \longrightarrow
(\R^2)^{\frac{\phi(n)}{2}-1}\,,$$ i.e., a map given by $.^{\star_{n}}\!\equiv 0$, if
$n\in\{3,4,6\}$, and by $z\longmapsto
\big(\sigma_2(z),\dots,\sigma_{\phi(n)/2}(z)\big)$
otherwise.
Then, $[\OO_{n}]\widetilde{\hphantom{a}}^{_{n}}$ is a
Minkowski representation of the maximal order $\OO_{n}$ of
$\mathbbm{K}_{n}$; cf. Proposition~\ref{p1}(a) and see~\cite[Ch. 2,
Sec. 3]{Bo}. It follows that $[\OO_{n}]\widetilde{\hphantom{a}}^{_{n}}$ is a (full) lattice in
$\R^{2}\times(\R^2)^{\phi(n)/2-1}$, meaning that it is a
co-compact 
discrete subgroup of the Abelian group
$\R^{2}\times(\R^2)^{\phi(n)/2-1}$ (the quotient group is a $\phi(n)$-dimensional torus). Here,
 this is equivalent to the existence of $\phi(n)$ $\R$-linearly
 independent vectors in $\R^{2}\times(\R^2)^{\phi(n)/2-1}$
 whose $\Z$-linear hull equals
 $[\OO_{n}]\widetilde{\hphantom{a}}^{_{n}}$, compare~\cite[Ch. 2, Sec. 3
   and 4]{Bo}. In fact, the set
   $$\left\{1\widetilde{\hphantom{a}}^{_{n}},(\zetan)\widetilde{\hphantom{a}}^{_{n}},\dots,(\zetan^{\phi(n)-1})\widetilde{\hphantom{a}}^{_{n}}\right\}$$ has this property; cf. Proposition~\ref{p1} and Remark~\ref{r1}. Further, note that the image $[\OO_{n}]^{\star_{n}}$ is dense in $(\R^2)^{\phi(n)/2-1}$; see~\cite[Section 3.2]{BG2}. 
\end{rem}

The class of
 {\em cyclotomic} model sets arises from cut and project schemes of the
following form. Consider the following diagram (cut
and project scheme), where we follow Moody~\cite{Moody}, modified
in the spirit of the algebraic setting of Pleasants~\cite{PABP}.

\begin{equation} \label{cutproj1}
\renewcommand{\arraystretch}{1.5}
\begin{array}{ccccc}
& \pi & & \pi_{\textnormal{\tiny int}}^{} & \vspace*{-2.0ex} \\
\!\!\!\!\!\!\!\!\!\!\!\!\!\!\!\R^{2} & \longleftarrow & \;\;\;\;\;\;\;\;\;\;\R^{2}\times(\R^2)^{\frac{\phi(n)}{2}-1}  & \longrightarrow & \;\;\;(\R^2)^{\frac{\phi(n)}{2}-1} \vspace*{1.2ex}\\
\cup\mbox{\tiny\, dense}\;\; &&\;\;\;\;\;\,\,\cup\mbox{\tiny\, lattice}&&\cup\mbox{\tiny\, dense}\\
 & \mbox{\tiny 1--1} & &  & \vspace*{-2.0ex} \\
 \!\!\!\!\!\!\!\!\!\!\!\!\!\OO_{n}& \longleftrightarrow & \,\,\,[\OO_{n}]\widetilde{\hphantom{a}}^{_{n}}& \longrightarrow &\!\!\!\![\OO_{n}]^{\star_{n}} \\
\end{array}
\end{equation}
As described above, one has
$$
[\OO_{n}]\widetilde{\hphantom{a}}^{_{n}}=\Big\{\Big(z,\underbrace{\big(\sigma_{2}(z),\dots ,\sigma_{\frac{\phi(n)}{2}}(z)\big)}_{=z^{\star_{n}}}\Big)\, \Big| \, z\in \OO_{n}\Big\}\,.
$$

\begin{defi}\label{cyclodef}
Let $n\in\N\setminus\{1,2\}$. Given any subset
$W\subset(\R^2)^{\phi(n)/2-1}$ with $\varnothing\, \neq\,
W^{\circ}\subset W\subset \overline{W^{\circ}}$ and
$\overline{W^{\circ}}$ compact, a so-called {\em window}, and any $t\in\R^2$, we obtain a planar
{\em model set} $$\varLambda_{n}(t,W) := t+\varLambda_{n}(W)$$
relative to any choice of the set $\{\sigma_{j} \,|\,
j\in\{2,\dots,\phi(n)/2\}\}$ as described above by setting
$$\varLambda_{n}(W):=\{z\in \OO_{n}\,|\,z^{\star_{n}}\in W\}\,.$$
Further, $\R^{2}$ (resp.,
$(\R^2)^{\phi(n)/2-1}$) is called the {\em physical} (resp.,
{\em internal}\/) space, the window $W$ is referred to as the {\em window of}
$\varLambda_{n}(t,W)$ and the map ${.}^{\star_{n}}$ is the so-called {\em star} map. Moreover, the dimension $c$ of the internal space of
$\varLambda_{n}(t,W)\in \mathcal{M}(\OO_{n})$, i.e., $c=\phi(n)-2$, is called the {\em co-dimension} of
$\varLambda_{n}(t,W)$. We set 
$$\mathcal{M}(\mathcal{O}_{n}):=
\left\{\varLambda_{n}(t,W)  \, \left | 
\, t\in\R^2, W\subset(\R^2)^{\frac{\phi(n)}{2}-1} \mbox{ is a window} 
 \right. \right \}\,. $$
Then, the class $\mathcal{CM}$ of {\em cyclotomic} model sets is defined as 

$$
\mathcal{CM}:=\bigcup_{n\,\in\,\mathbbm{N}\setminus \{1,2\}}\mathcal{M}(\OO_{n})\,.
$$
\end{defi}

\begin{rem}\label{propcm}
We refer the reader to~\cite{Moody,PABP} for details and related general settings,
and to~\cite{BM} for general background. Note that the star map is a
homomorphism of Abelian groups. Further, the co-dimension of a
cyclotomic model set is always an even number. Moreover, it is zero
if and only if $n\in\{3,4,6\}$. In the following, we present some properties of
cyclotomic model sets. Setting
$\varLambda:=\varLambda_{n}(t,W)\subset
\R^2$, one has that $\varLambda$ is a {\em Delone set} in $\R^2$ 
(i.e., $\varLambda$ is both uniformly discrete and relatively dense)
and 
has {\em finite local complexity} (i.e.,
 $\varLambda-\varLambda$ is closed and discrete). In fact, $\varLambda$ is even a {\em
  Meyer set} (i.e., $\varLambda$ is a Delone set and $\varLambda-\varLambda$ is uniformly
discrete); compare \cite{Moody}. 

Further, $\varLambda$ is an {\em
  aperiodic} cyclotomic model set (i.e., $\varLambda$ has no translational
symmetries) if and only if the star map is a monomorphism (i.e., if it is
injective). In fact, the kernel of the star map is the group of
translation symmetries of $\varLambda$; compare \cite{Moody} again. It follows that $\varLambda$ is aperiodic if and only if
$n\notin\{3,4,6\}$, i.e., the translates of the square
(resp., triangular) lattice are the only cyclotomic model sets with translation symmetries. 

Moreover, if
$\varLambda$ is {\em regular} (i.e., if the boundary $\partial W$ has Lebesgue
measure $0$ in $(\R^2)^{\phi(n)/2-1}$), $\varLambda$ is {\em
  pure point diffractive} (cf.~\cite{Schl}); if $\varLambda$ is {\em
  generic} (i.e., if $[\OO_{n}]^{\star_{n}}\cap\, \partial W = \varnothing$),
$\varLambda$ is {\em repetitive}; see~\cite{Schl}. If $\varLambda$ is
regular, the frequency of repetition of finite
patches is well defined (cf. ~\cite{Schl2}) and, moreover, if $\varLambda$, for a given $n$, is both generic and regular, and, if the window $W$ has 
$m$-fold cyclic symmetry with $m$ a divisor of
$N(n)$ ($N(n)$ is the function from
(\ref{eq})) and all in a suitable representation of the cyclic group $\mathsf{C}_m$ of order $m$, then $\varLambda$ has 
$m$-fold cyclic symmetry
in the sense of symmetries of LI-classes,
meaning that a discrete structure has a certain symmetry if the
original and the transformed structure are locally indistinguishable
(LI); see~\cite{B} for details. Although not all of these properties are used below, they actually
show the similarity of aperiodic model sets with lattices -- except their lack of periods.
\end{rem}

\begin{defi}\label{corr}
Let $n\,\in\,\mathbbm{N}\setminus \{1,2,3,4,6\}$, let
$\lambda\in\oo_{n}$, and let $\OO_{n}
 \longrightarrow [\OO_{n}]\widetilde{\hphantom{a}}^{_{n}}$ be a fixed
Minkowski embedding of $\OO_{n}$. Further, let $.^{\star_{n}}$ be the induced star map;
 see Remark~\ref{mink}. 
\begin{itemize}
\item[(a)]
$m_{\lambda}^{_{(n)}}$ denotes the $\Z$-module endomorphism of
$\OO_{n}$, given by $z \longmapsto \lambda z$.\smallskip  
\item[(b)]
$(m_{\lambda}^{_{(n)}})\widetilde{\hphantom{a}}^{_{n}}$ denotes the $\Z$-module (lattice)
endomorphism of
$[\OO_{n}]\widetilde{\hphantom{a}}^{_{n}}$, given by $$(z,z^{\star_{n}}) \longmapsto (\lambda z,(\lambda
z)^{\star_{n}})\,.$$
\item[(c)] 
$(m_{\lambda}^{_{(n)}}){}^{\star_{n}}$ denotes the $\Z$-module endomorphism of $[\OO_{n}]^{\star_{n}}$, given by $$z^{\star_{n}} \longmapsto (\lambda z)^{\star_{n}}\,.$$
\end{itemize}
\end{defi}

\begin{rem}\label{rtilde}
Note that, for
 $n\notin\{3,4,6\}$, every star map $.^{\star_{n}}$ is a $\Z$-module monomorphism;
 cf. Remark~\ref{propcm}.
 Let $n\,\in\,\mathbbm{N}\setminus \{1,2,3,4,6\}$, let
$\lambda\in\oo_{n}$ (resp.,
$\lambda\in\oo_{n}^{\times}$), and let $.\,\widetilde{\hphantom{a}}^{_{n}}\! : \,\OO_{n}
 \longrightarrow [\OO_{n}]\widetilde{\hphantom{a}}^{_{n}}$ be a fixed
Minkowski embedding of $\OO_{n}$. Further, let $.^{\star_{n}}$ be the induced star map. Then, the
$\Z$-module endomorphism (resp., automorphism)
 $m_{\lambda}^{_{(n)}}$ of $\OO_{n}$ corresponds via the chosen
Minkowski embedding of
$\OO_{n}$ to the $\Z$-module (lattice)
endomorphism (resp., automorphism) $(m_{\lambda}^{_{(n)}})\widetilde{\hphantom{a}}^{_{n}}$ of
$[\OO_{n}]\widetilde{\hphantom{a}}^{_{n}}$. Further, $m_{\lambda}^{_{(n)}}$ corresponds via the
$\Z$-module automorphism $.^{\star_{n}}\!:\,\OO_{n}\longrightarrow [\OO_{n}]^{\star_{n}}$, given by the star
map, to the $\Z$-module endomorphism (resp.,
automorphism) $(m_{\lambda}^{_{(n)}}){}^{\star_{n}}$ of $[\OO_{n}]^{\star_{n}}$.
\end{rem}

\begin{defi}
For $n\,\in\,\mathbbm{N}\setminus \{1,2,3,4,6\}$, we denote by $\Arrowvert \cdot \Arrowvert_{\infty}$ the maximum norm
on ${(\mathbbm{R}^{2})}^{\phi(n)/2-1}$ with respect to the
Euclidean norm on all factors $\R^2$.
\end{defi}

\begin{lem}\label{r2}
Let $n\,\in\,\mathbbm{N}\setminus \{1,2,3,4,6\}$, let $.\,\widetilde{\hphantom{a}}^{_{n}}\! : \,\OO_{n}
 \longrightarrow [\OO_{n}]\widetilde{\hphantom{a}}^{_{n}}$ be a fixed
Minkowski embedding of $\OO_{n}$, and let $.^{\star_{n}}$ be the induced star map. Then, for any
PV-number $\lambda$ of $($full$)$ degree $\phi(n)/2$ in $\oo_{n}$, the $\Z$-module
endomorphism $(m_{\lambda}^{_{(n)}}){}^{\star_{n}}$ is contractive, i.e., there is a
$\xi \in (0,1)$ such that the inequality $$\Arrowvert (m_{\lambda}^{_{(n)}}){}^{\star_{n}}(z^{\star_{n}})\Arrowvert_{\infty} \leq \xi\, \Arrowvert z^{\star_{n}}\Arrowvert_{\infty}$$ holds for all $z\in \OO_{n}$.
\end{lem}
\begin{proof}
Since $\lambda$ is an algebraic integer of full degree,
the set $\{\sigma_{1}(\lambda),\dots,\sigma_{\phi(n)/2}(\lambda)\}$ equals the set of (algebraic)
conjugates of $\lambda$. To see this, note that the set of
\mbox{(co-)}restrictions $\{\sigma_{1}|_{\mathbbm{k}_{n}}^{\mathbbm{k}_{n}},\dots,\sigma_{\phi(n)/2}|_{\mathbbm{k}_{n}}^{\mathbbm{k}_{n}}\}$ equals the Galois group
$G(\mathbbm{k}_{n} / \mathbbm{Q})$; compare
Corollary~\ref{cr5}. Since $\lambda$ is a PV-number, the last
observation shows that $$\xi:=\operatorname{max}\big\{\lvert \sigma_{j}(\lambda)\rvert\,\big |
\, j\in\left\{2,\dots,\phi(n)/2\right\}\big\} \in (0,1)\,.$$ The assertion follows.
\end{proof}

The following result will play a key role.

\begin{lem}\label{dilate}
Let $n\,\in\,\mathbbm{N}\setminus \{1,2\}$ and let
$\varLambda_{n}(t,W)\in \mathcal{M}(\OO_{n})$ be a cyclotomic model
set. Then, for any finite set $F\subset t+\mathbbm{K}_{n}$, there is a homothety $h\!:\, \C \longrightarrow \C$ such that $h[F]\subset \varLambda_{n}(t,W)$. In particular, if $t=0$ and $0\in \, W^{\circ}$, there is even a dilatation $d\!:\, \C \longrightarrow \C$ such that $d[F]\subset \varLambda_{n}(0,W)$.  
\end{lem}
\begin{proof}
Without loss of generality, we may assume that $F$ is non-empty. We
consider the $\mathbbm{Q}$-coordinates of the elements of $F-t$ with
respect to the $\mathbbm{Q}$-basis $\{1,\zetan,\zetan^2, \dots
,\zetan^{\phi(n)-1}\}$ of $\mathbbm{K}_{n}$
(cf. Proposition~\ref{gau}) and let $l\in \mathbbm{N}$ be the least
common multiple of all their denominators. Then, by Remark~\ref{r1}, we
get $l(F-t) \subset \OO_{n}$. If $n\in\{3,4,6\}$, we are done by
 defining the homothety $h\!:\, \C \longrightarrow \C$ by
$$
z\longmapsto lz + (-lt+t)\,.
$$
Secondly, suppose that $n\notin\{3,4,6\}$. Let $.^{\star_{n}}$ be the star map that is used in
 the construction of $\varLambda_{n}(t,W)$. From
$W^{\circ}\neq \varnothing$ and the denseness of $[\OO_{n}]^{\star_{n}}$ in
$(\R^2)^{\phi(n)/2-1}$, there follows the existence of a
suitable $z_{0}\in \OO_{n}$ with $z_{0}{}^{\star_{n}}\in
\,W^{\circ}$. Consider the open neighbourhood $$V:= \,W^{\circ} - \, z_{0}{}^{\star_{n}}$$ of $0$ in
$(\R^2)^{\frac{\phi(n)}{2}-1}$. Next, we choose a PV-number $\lambda$
of degree $\phi(n)/2$ in $\oo_{n}$; see
Lemma~\ref{pisot}. Consider the $\Z$-module endomorphism $(m_{\lambda}^{_{(n)}}){}^{\star_{n}}\! : \, [\OO_{n}]^{\star_{n}} \longrightarrow
[\OO_{n}]^{\star_{n}}$, as defined in Definition~\ref{corr}(c). Since $\lambda$ is a
PV-number, Lemma~\ref{r2} shows that $(m_{\lambda}^{_{(n)}}){}^{\star_{n}}$ is
contractive (in the sense which was made precise in that lemma). Since all norms on
${(\mathbbm{R}^{2})}^{\phi(n)/2-1}$ are equivalent, it follows the existence of a suitable $k\in\mathbbm{N}$ such that $$\big((m_{\lambda}^{_{(n)}}){}^{\star_{n}}\big)^{k}\big[\,[l(F-t)]^{\star_{n}}\big]\subset V\,.$$ It follows that $\big\{(\lambda^{k} z + z_{0})^{\star_{n}}\, |\, z\in l(F-t)\big\}\subset \,W^{\circ}$ and further that $h[F]\subset \varLambda_{n}(t,W)$, where $h\!:\, \C \longrightarrow \C$ is the homothety given by $$z \longmapsto (l\lambda^{k}) z + (z_{0} - (l\lambda^{k})t + t)\,.$$ The additional statement follows immediately from the foregoing proof in connection with the trivial observation that $0\in \OO_{n}$ maps, under the star map $.^{\star_{n}}$, to $0\in (\R^2)^{\phi(n)/2-1}$.       
\end{proof}

\begin{rem}
As the general structure shows, the restriction in Lemma~\ref{dilate}
to cyclotomic model sets is by no means necessary. In fact, the
related 
general algebraic setting of~\cite{PABP} can be used to extend this result also to higher dimensions. There, the role of the
PV-numbers (or PV-units) will be taken by suitable
hyperbolic lattice endomorphisms (or automorphisms). 
\end{rem}

We shall also need the following result from Weyl's theory of
uniform distribution, which is concerned with an analytical aspect of regular cyclotomic model sets.

\begin{theorem}[Weyl]\label{weyl}
Let $n\,\in\,\mathbbm{N}\setminus \{1,2,3,4,6\}$ and let $\varLambda:=\varLambda_{n}(t,W)\in\mathcal{M}(\OO_{n}) $ be a
regular $($aperiodic$)$ cyclotomic model
set, say constructed by use of the star map $.^{\star_{n}}$. Then, for
all $a\in\R^2$, one has the identity
$$
\lim_{R\rightarrow \infty}\,\, \frac{1}{\operatorname{card}(\varLambda\cap
  B_{R}(a))}\sum_{x\in \varLambda\cap
  B_{R}(a)}x^{\star_{n}}=\frac{1}{\operatorname{vol}(W)}\int_{W}y\,{\rm d}\lambda(y)\,,
$$
where $\lambda$ denotes the Lebesgue measure on $(\R^2)^{\frac{\phi(n)}{2}-1}$.
\end{theorem}
\begin{proof}
This follows from the multi-dimensional version of~\cite[Theorem 3]{Moody}.
\end{proof}

\section{Affinely Regular Polygons in Cyclotomic Model Sets}

\subsection{Affinely Regular Polygons in Rings of Cyclotomic Integers}\label{char}
Gardner and Gritzmann~\cite[Theorem 4.1]{GG} have shown that there is an affinely regular $m$-gon in the square lattice $\Z^2=\Z[\zeta_{4}]=\OO_{4}$ if and only if $m \in \{3,4,6\}$. We start off with a generalization.    

\begin{theorem}\label{th1}
Let $m,n\in\mathbbm{N}$ with $m,n\geq 3$. The following statements are equivalent:
\begin{itemize}
\item[(i)]
There is an affinely regular $m$-gon in $\OO_{n}$.
\item[(ii)]
$\mathbbm{k}_{m}\subset\mathbbm{k}_{n}$.
\item[(iii)]
$m \in \{3,4,6\}$, or $\mathbbm{K}_{m}\subset \mathbbm{K}_{n}$.
\item[(iv)]
$m \in \{3,4,6\}$, or $m|n$, or $m=2d$ with $d$ an odd divisor of $n$.
\item[(v)]
$m \in \{3,4,6\}$, or $\OO_{m}\subset \OO_{n}$.
\item[(vi)]
$\oo_{m}\subset\oo_{n}$.
\end{itemize}
\end{theorem}
\begin{proof}
For (i) $\Rightarrow$ (ii), let $P$ be an affinely regular $m$-gon in
$\OO_{n}$. There is then a non-singular affine transformation $\psi
\!:\, \C \longrightarrow \C$ with $\psi[R_{m}] = P$, where $R_{m}$ is
the regular $m$-gon with vertices given in complex form by $1, \zetam,
\dots, \zetam^{m - 1}$. If $m\in\{3,4,6\}$, condition (ii)
 holds trivially. Suppose $6\neq m\geq 5$. The pairs $\{1,\zetam\}$, $\{\zetam^{-1},\zetam^{2}\}$ lie on parallel lines and so do their images under $\psi$. Therefore,
$$\frac{\Arrowvert \zetam^{2} - \zetam^{-1} \Arrowvert}{\Arrowvert \zetam - 1 \Arrowvert} = \frac{\Arrowvert \psi(\zetam^{2}) - \psi(\zetam^{-1}) \Arrowvert}{\Arrowvert \psi(\zetam) - \psi(1) \Arrowvert}\,.$$ Moreover, we get the relation $$(1 + \zetam + \bar{\zeta}_{m})^2 = (1 + \zetam + \zetam^{-1})^2 = \frac{\Arrowvert \zetam^{2} - \zetam^{-1} \Arrowvert^2}{\Arrowvert \zetam - 1 \Arrowvert^2} = \frac{\Arrowvert \psi(\zetam^{2}) - \psi(\zetam^{-1}) \Arrowvert^2}{\Arrowvert \psi(\zetam) - \psi(1) \Arrowvert^2} \in \mathbbm{k}_{n}\,.$$ The pairs $\{\zetam^{-1},\zetam\}$, $\{\zetam^{-2},\zetam^{2}\}$ also lie on parallel lines. An argument similar to that above yields $$ (\zetam + \bar{\zeta}_{m})^2 = (\zetam + \zetam^{-1})^2 = \frac{\Arrowvert \zetam^{2} - \zetam^{-2} \Arrowvert^2}{\Arrowvert \zetam - \zetam^{-1} \Arrowvert^2} \in \mathbbm{k}_{n}\,.$$ By subtracting these equations, we get $$2(\zetam + \bar{\zeta}_{m}) + 1 \in \mathbbm{k}_{n}\,,$$ and hence $\zetam + \bar{\zeta}_{m} \in \mathbbm{k}_{n}$, the latter being equivalent to the inclusion of the fields $ \mathbbm{k}_{m} \subset \mathbbm{k}_{n}$.

As an immediate consequence of Lemma~\ref{unique2}(b), we get (ii) $\Rightarrow$ (iii). Conditions (iii) and (iv) are equivalent by Lemma~\ref{incl}, and (iii) $\Rightarrow$ (v) $\Rightarrow$ (vi) follows from Proposition~\ref{p1}.    

For (vi) $\Rightarrow$ (i), let $R_{m}$ again be the regular $m$-gon as defined in the step (i) $\Rightarrow$ (ii). Since $m,n \geq 3$, the sets $\{1,\zetam\}$ and $\{1,\zetan\}$ are $\R$-bases of $\C$. Hence, we can define an $\R$-linear map $L_{m}^{n}\! :\, \C \longrightarrow \C$ as the linear extension of $1 \longmapsto 1$ and $\zetam \longmapsto \zetan$. Obviously, $L_{m}^{n}$ is non-singular. Then, using $\oo_{m}\subset\oo_{n}$, one can see, by means of Lemma~\ref{Oo}(a), that the vertices of $L_{m}^{n}[R_{m}]$, i.e., $L_{m}^{n}(1), L_{m}^{n}(\zetam), \dots, L_{m}^{n}(\zetam^{m - 1})$, lie in $\OO_{n}$ (in fact, $L_{m}^{n}$ maps the whole $\oo_{m}$-module $\OO_{m}$ into the $\oo_{n}$-module $\OO_{n}$), whence $L_{m}^{n}[R_{m}]$ is a polygon in $\OO_{n}$. This shows that there is an affinely regular $m$-gon in $\OO_{n}$.              
\end{proof}

\begin{rem}\label{r3}
Alternatively, there is the following direct argument for (v) $\Rightarrow$ (i) in Theorem~\ref{th1}. If $m\in\{3,4,6\}$, then $L_{m}^{n}[R_{m}]$ from the proof of (vi) $\Rightarrow$ (i) above is an affinely regular $m$-gon in $\OO_{n}$. Otherwise, one can simply choose the regular $m$-gon $R_{m}$ (formed by the $m$th roots of unity) itself. 
\end{rem}

Although the equivalence (i) $\Leftrightarrow$ (iv) in Theorem~\ref{th1} is a satisfactory characterization, the following corollary deals with the two cases where condition (ii) can be used more effectively.

\begin{coro}\label{cor2}
Let $m\in \mathbbm{N}$ with $m\geq 3$. Consider $\phi$ on $\{n\in \mathbbm{N}\, |\, n\; \not\equiv \;2  \;(\operatorname{mod} 4)\}$. Then, one has:
\begin{itemize}
\item[(a)] If $\phi(n)/2=1$ $($i.e., $n\in\{3,4\}$; see
  Lemma~$\ref{phin2p}({\rm a}))$, there is an affinely regular $m$-gon in $\OO_{n}$ if and only if $m \in \{3,4,6\}$, i.e., the affinely regular polygons in $\OO_{n}$ in this case are exactly the affinely regular triangles, parallelograms and hexagons.  
\item[(b)] If $\phi(n)/2\in\mathbbm{P}$ $($i.e.,
  $n\in\mathbbm{S}$; see Lemma~$\ref{phin2p}({\rm b}))$, there is an affinely regular $m$-gon in $\OO_{n}$ if and only if 
$$\left\{
\begin{array}{ll}
m \in \{3,4,6,n\}, & \mbox{if $n=8$ or $n=12$,}\\
m \in \{3,4,6,n,2n\}, & \mbox{otherwise.}
\end{array}\right.
$$ 
\end{itemize}
\end{coro}
\begin{proof}
If $\phi(n)/2=1$, condition (ii) of Theorem~\ref{th1}
specializes to $\mathbbm{k}_{m} = \Q$. This is equivalent to
$\phi(m)=2$ by Corollary \ref{cr5}, which means $m\in\{3,4,6\}$. This
proves the part~(a). 

If $\phi(n)/2 \in \mathbbm{P}$, we have $[\mathbbm{k}_{n}:\Q]=
\phi(n)/2 \in \mathbbm{P}$. Hence, we get, by means of
condition (ii) of Theorem~\ref{th1} and Lemma~\ref{degree}, either
$\mathbbm{k}_{m}=\Q$ or $\mathbbm{k}_{m}=\mathbbm{k}_{n}$. The former
case implies $m\in\{3,4,6\}$ as in the proof of the part (a), while the proof
 follows from Lemma~\ref{unique2}(a) in conjunction with Corollary~\ref{unique} in the latter case.
\end{proof}

\subsection{Application to Cyclotomic Model Sets}

\begin{coro}\label{th1mod}
Let $m\in \mathbbm{N}$ with $m\geq 3$. Further, let $n\,\in\,\mathbbm{N}\setminus \{1,2\}$ and let $\varLambda_{n}(t,W)\in \mathcal{M}(\OO_{n})$ be a cyclotomic model set. The following statements are equivalent:
\begin{itemize}
\item[(i)]
There is an affinely regular $m$-gon in $\varLambda_{n}(t,W)$.
\item[(ii)]
$\mathbbm{k}_{m}\subset\mathbbm{k}_{n}$.
\item[(iii)]
$m \in \{3,4,6\}$, or $\mathbbm{K}_{m}\subset \mathbbm{K}_{n}$.
\item[(iv)]
$m \in \{3,4,6\}$, or $m|n$, or $m=2d$ with $d$ an odd divisor of $n$.
\item[(v)]
$m \in \{3,4,6\}$, or $\OO_{m}\subset \OO_{n}$.
\item[(vi)]
$\oo_{m}\subset\oo_{n}$.
\end{itemize}
\end{coro}
\begin{proof}
The assertion is a consequence of Theorem~\ref{th1} and Lemma~\ref{dilate}, since homotheties are non-singular affine transformations.
\end{proof}

Similarly, we employ Corollary~\ref{cor2}(b) to obtain

\begin{coro}\label{cor3}
Let $m\in \mathbbm{N}$ with $m\geq 3$. Further, let $n\in\mathbbm{S}$
and $\varLambda_{n}(t,W)\in \mathcal{M}(\OO_{n})$ be a $($aperiodic$)$ cyclotomic model set. Then, there is an affinely regular $m$-gon in $\varLambda_{n}(t,W)$ if and only if 
$$\left\{
\begin{array}{ll}
m \in \{3,4,6,n\}, & \mbox{if $n=8$ or $n=12$,}\\
m \in \{3,4,6,n,2n\}, & \mbox{otherwise.}\hspace*{10.735em}\qed\hspace*{-10.735em}
\end{array}\right. 
$$
\end{coro}

\subsection{Examples and Comments}\label{sec4}

In the following, let $R_{m}$ and $L_{m}^{n}[R_{m}]$ be as
in Remark~\ref{r3}. The first two examples are of the form
$\varLambda_{n}(0,W)\in \mathcal{M}(\OO_{n})$ for $n\in \{3,4\}$ and
hence, necessarily, $W=\{0\}$.

\smallskip
\begin{itemize}

\begin{figure}
\begin{minipage}[t]{0,48\textwidth}
\epsfxsize=\textwidth\epsfbox{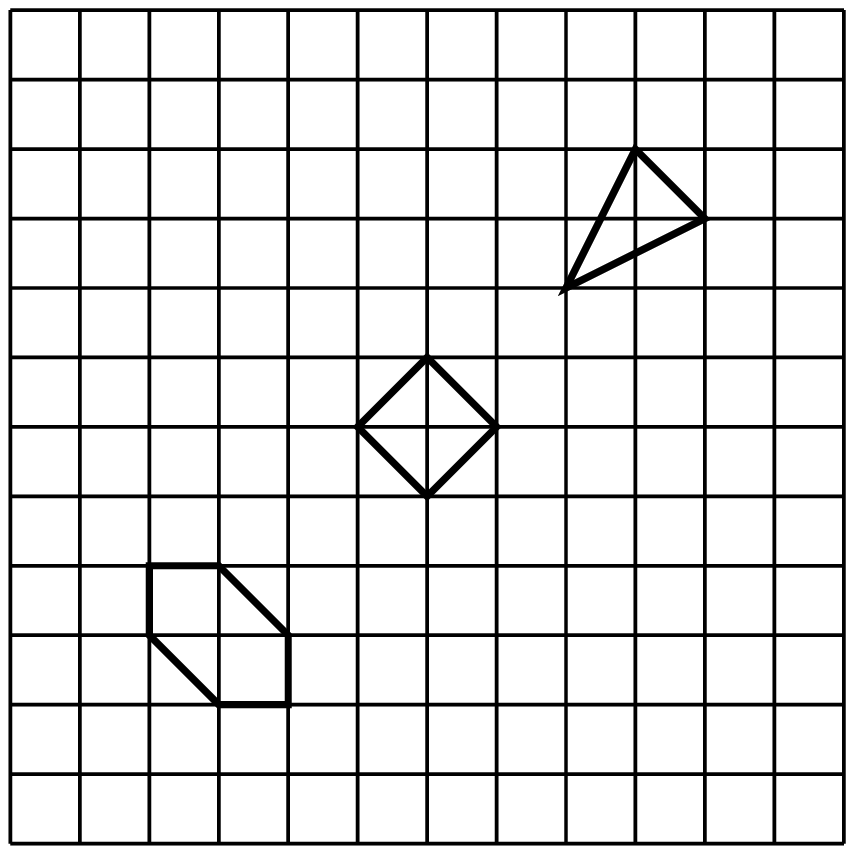}
\end{minipage}
\hfill
\begin{minipage}[t]{0,48\textwidth}
\vspace*{-7.0cm}
\epsfxsize=\textwidth\epsfbox{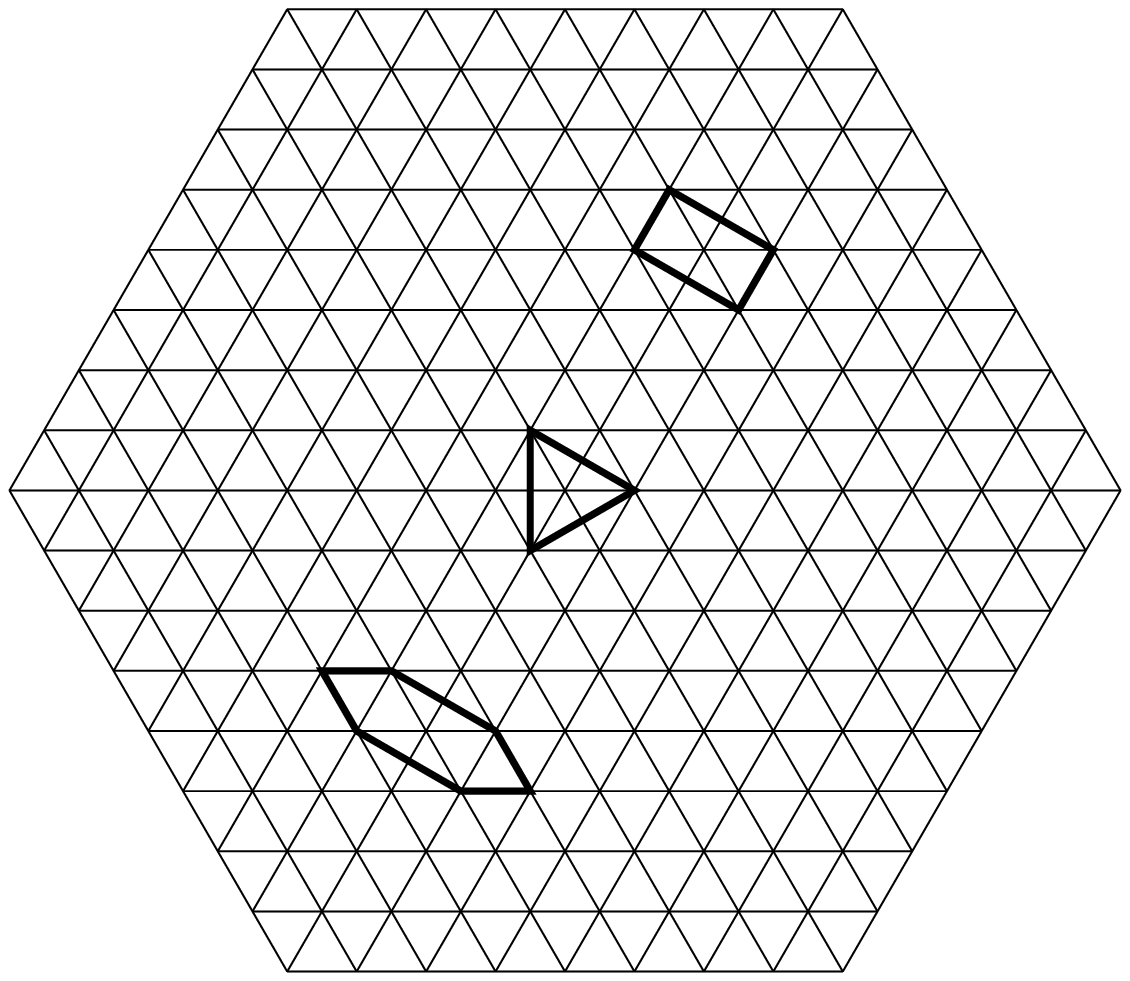}
\end{minipage}
\caption{On the left, a central patch of the square tiling with vertex
  set $\varLambda_{\text{SQ}}$ with $L_{4}^{4}[R_{4}]$ and translates of
  $L_{3}^{4}[R_{3}]$ and $L_{6}^{4}[R_{6}]$ as described in example
  {\rm (SQ)} is shown. On
  the right, a central patch of the triangular tiling with vertex set $\varLambda_{\text{TRI}}$ with $L_{3}^{3}[R_{3}]$ and translates of  $L_{4}^{3}[R_{4}]$ and $L_{6}^{3}[R_{6}]$ as described in example
  {\rm (TRI)} is shown.}
\label{fig:squaretri}
\end{figure}

\item[(SQ)]
The planar generic regular periodic cyclotomic model set with $4$-fold
cyclic symmetry associated with the well-known square tiling is the
square lattice, which can be described in algebraic terms as $\varLambda_{\text{SQ}} :=\varLambda_{4}(0,W)=\Z[i] = \OO_{4}$. By
Corollary~\ref{cor2}(a), there is an affinely regular $m$-gon in
$\OO_{4}$ if and only if $m \in \{3,4,6\}$. Moreover, Remark~\ref{r3} shows that
the affinely regular $3$-, $4$- and $6$-gons $L_{3}^{4}[R_{3}]$,
$L_{4}^{4}[R_{4}]$ and $L_{6}^{4}[R_{6}]$ are polygons in $\OO_{4}$;
see Figure~\ref{fig:squaretri}.
\item[(TRI)]
The planar generic regular periodic cyclotomic model set with $6$-fold
cyclic symmetry associated with the well-known triangle tiling is the
triangle lattice (also
commonly known as the hexagonal lattice), which can be described in algebraic terms as $\varLambda_{\text{TRI}} :=\varLambda_{3}(0,W)=\OO_{3}$. By Corollary~\ref{cor2}(a),
there is an affinely regular $m$-gon in $\OO_{3}$ if and only if $m \in
\{3,4,6\}$. Moreover, Remark~\ref{r3} shows that the affinely regular
$3$-, $4$- and $6$-gons $L_{3}^{3}[R_{3}]$, $L_{4}^{3}[R_{4}]$ and
$L_{6}^{3}[R_{6}]$ are polygons in $\OO_{3}$; see
Figure~\ref{fig:squaretri}.
\end{itemize}

All further examples below are aperiodic cyclotomic model sets of the
form $\varLambda_{n}(0,W)\in \mathcal{M}(\OO_{n})$ for suitable $n\in
\mathbbm{S}$ and satisfy $0\, \in\, W^{\circ}$. An analysis of the
proof of  Lemma~\ref{dilate} in connection with Corollary~\ref{cor3}
and Remark~\ref{r3} shows the existence of an expansive dilatation $d\!:\, \C \longrightarrow \C$, given by multiplication by a suitable non-negative integral power of a PV-number of $\oo_{n}$, such that $d[L_{m}^{n}[R_{m}]] \subset \varLambda_{n}(0,W)$, if $m\in\{3,4,6\}$ (respectively, $d[R_{m}] \subset \varLambda_{n}(0,W)$, otherwise), with $n,m$ in accordance with the statement of Corollary~\ref{cor3}. Let us demonstrate this in some more detail.

\smallskip
\begin{itemize}
\item[(AB)]
The planar generic regular aperiodic cyclotomic model set with $8$-fold cyclic symmetry associated with the well-known Ammann-Beenker tiling~\cite{am,bj,ga} can be described in algebraic terms as
$$\varLambda_{\text{AB}} := \{z \in \OO_{8}\, | \,z^{\star_{8}} \in O\}\,,$$
where the star map $.^{\star_{8}}$ is the Galois automorphism
(cf. Proposition \ref{gau}) in $G(\mathbbm{K}_{8}/ \mathbbm{Q})$,
defined by $\zeta_{8} \longmapsto \zeta_{8}^3$, and the window $O$ is
the regular octagon centred at the origin, with vertices in the
directions that arise from the $8$th roots of unity by a rotation
through $\pi/8$, and of unit edge length. This construction also gives a tiling with squares and rhombi, both having edge length $1$; see Figure~\ref{fig:ab}.

By Corollary~\ref{cor3}, there is an affinely regular $m$-gon in
$\varLambda_{\text{AB}}$  if and only if $m \in \{3,4,6,8\}$. See
Figure~\ref{fig:ab} for all assertions below. The affinely regular
\mbox{$3$-,} $4$- and $8$-gons $L_{3}^{8}[R_{3}]$, $L_{4}^{8}[R_{4}]$ and
$R_{8}$ are polygons in $\varLambda_{\text{AB}}$. For an affinely
regular $6$-gon with its vertices in $\varLambda_{\text{AB}}$,
consider the expansive dilatation of $L_{6}^{8}[R_{6}]$, which is given by
multiplication by the PV-unit $1+\sqrt{2}$ in $\oo_{8}$ (the
fundamental unit in $\oo_{8}$), i.e., the convex $6$-gon with vertices
$1+\sqrt{2}$, $(1+\sqrt{2})\zeta_{8}$, $(1+\sqrt{2})(-1 +\zeta_{8})$,
 $-(1+\sqrt{2})$, $-(1+\sqrt{2})\zeta_{8}$, $(1+\sqrt{2})(1-\zeta_{8})$. 
Here, also the PV-number $2+\sqrt{2}$ of (full) degree $2$ would be suitable.

\begin{figure}
\centerline{\epsfysize=0.55\textwidth\epsfbox{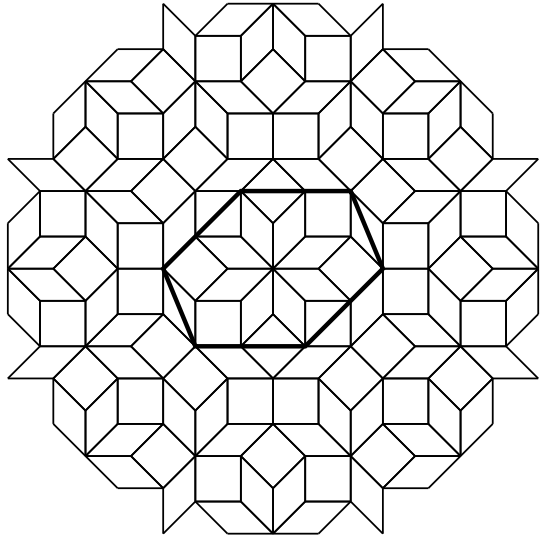}}
\caption{A central patch of the eightfold symmetric Ammann-Beenker tiling with vertex set $\varLambda_{\rm AB}$. In the tiling, the affinely regular $6$-gon as described in example {\rm (AB)} is shown, the other solutions being rather obvious.}
\label{fig:ab}
\end{figure}

\item[(TTT)]
The planar regular aperiodic cyclotomic model set with $10$-fold cyclic symmetry associated with the T\"ubingen triangle tiling~\cite{bk1,bk2} can be described in algebraic terms as
$$\varLambda_{\text{TTT}}^{u} := \{z \in \OO_{5}\, | \,z^{\star_{5}} \in u+W\}\,,$$
where the star map $.^{\star_{5}}$ is the Galois automorphism
(cf. Proposition \ref{gau}) in $G(\mathbbm{K}_{5}/ \mathbbm{Q})$,
defined by $\zeta_{5} \longmapsto \zeta_{5}^2$, and the window $W$ is
the regular decagon centred at the origin, with vertices in the
directions that arise from the $10$th roots of unity by a rotation
through $\pi/10$, and of edge length
$\tau/\sqrt{\tau+2}$, where $\tau$ denotes the golden ratio, i.e., $\tau=(1+\sqrt{5})/2$. Furthermore, $u$ is an element of $\R^{2}$. $\varLambda_{\text{TTT}}^{0}$ is not generic, while generic examples are obtained by shifting the window, i.e., $\varLambda_{\text{TTT}}^{u}$ is generic for almost all $u\in \R^{2}$. Generic $\varLambda_{\text{TTT}}^{u}$ always give a triangle tiling with long (short) edges of length $1$ ($1/\tau$). See Figure~\ref{fig:ttt} for a generic example which we call $\varLambda_{\text{TTT}}$; different generic choices of $u$ result in locally indistinguishable T\"ubingen triangle tilings.

Now, again by Corollary~\ref{cor3}, there is an affinely regular
$m$-gon in $\varLambda_{\text{TTT}}$  if and only if $m \in
\{3,4,5,6,10\}$. See Figure~\ref{fig:ttt} for all assertions
below. The affinely regular $4$-, $5$- and $10$-gons
$L_{4}^{5}[R_{4}]$, $R_{5}$ and $R_{10}$ are polygons in
$\varLambda_{\text{TTT}}$. For an affinely regular polygon in
$\varLambda_{\text{TTT}}$ with $3$ vertices, consider the expansive dilatation
of $L_{3}^{5}[R_{3}]$, which is given by multiplication by the PV-unit
$\tau^{2}$ in $\oo_{5}$, i.e., the convex $3$-gon with vertices
$\tau^{2}, \tau^{2}\zeta_{5}, \tau^{2}(-1 - \zeta_{5})$. Note the
identity $\tau^{2}=2+\frac{1}{\tau}$. Note further that also the
PV-unit $1+\frac{1}{\tau}=\tau$ would be suitable. For an affinely
regular polygon in $\varLambda_{\text{TTT}}$ with $6$ vertices,
consider the expansive dilatation of $L_{6}^{5}[R_{6}]$, which is again given by multiplication with the PV-unit $\tau^{2}=2+1/\tau$ in $\oo_{5}$, i.e., the convex $6$-gon with vertices $\tau^{2}$, $\tau^{2}\zeta_{5}$, $\tau^{2}(-1+\zeta_{5})$, $-\tau^{2}$, $-\tau^{2}\zeta_{5}$, $\tau^{2}(1-\zeta_{5})$. Again, also the PV-unit $1+1/\tau=\tau$ would be suitable.  

\begin{figure}
\centerline{\epsfxsize=0.55\textwidth\epsfbox{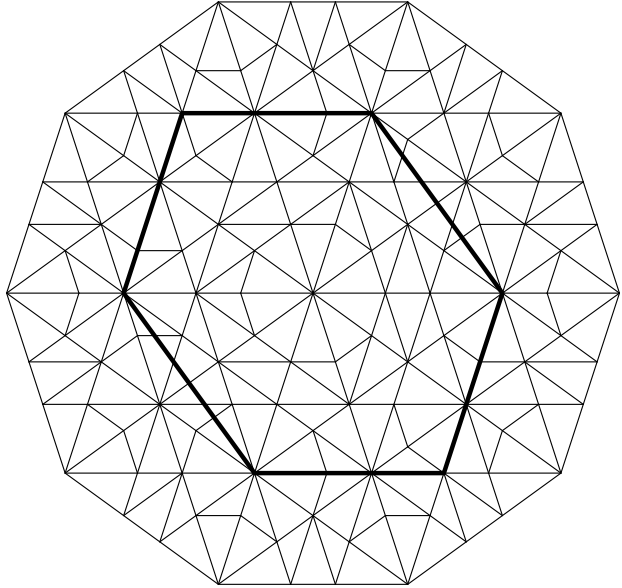}}
\caption{A central patch of the tenfold symmetric T\"ubingen triangle tiling. Therein, the affinely regular $6$-gon as described in example {\rm (TTT)} is marked, the other solutions being rather obvious.}
\label{fig:ttt}
\end{figure}

\item[(S)]
The planar regular aperiodic cyclotomic model set with $12$-fold cyclic symmetry associated with the shield tiling~\cite{ga} can be described in algebraic terms as
$$\varLambda_{\text{S}}^{u} := \{z \in \OO_{12}\, | \,z^{\star_{12}} \in u+W\}\,,$$
where the star map $.^{\star_{12}}$ is the Galois automorphism
(cf. Proposition~\ref{gau}) in $G(\mathbbm{K}_{12}/ \mathbbm{Q})$,
defined by $\zeta_{12} \longmapsto \zeta_{12}^5$, and the window $W$
is the regular dodecagon centred at the origin, with vertices in the
directions that arise from the $12$th roots of unity by a rotation
through $\pi/12$, and of edge length $1$. Again, $u$ is an
element of $\R^{2}$. $\varLambda_{\text{S}}^{0}$ is not generic, while
generic examples are obtained by shifting the window, i.e.,
$\varLambda_{\text{S}}^{u}$ is generic for almost all $u\in
\R^{2}$. The shortest distance between points in a generic
$\varLambda_{\text{S}}^{u}$ is $(\sqrt{3}-1)/\sqrt{2}$. Joining
such points by edges results in a shield tiling, i.e., a tiling with
triangles, squares and so-called shields as tiles, all having edge length $(\sqrt{3}-~1)/\sqrt{2}$. See Figure~\ref{fig:s} for a generic example which we call $\varLambda_{\text{S}}$; different generic choices of $u$ result in locally indistinguishable shield tilings.

Once more, by Corollary~\ref{cor3}, there is an affinely regular
$m$-gon in $\varLambda_{\text{S}}$ if and only if $m \in \{3,4,6,12\}$. See
Figure~\ref{fig:s} for all assertions below. It is immediate that
 $L_{3}^{12}[R_{3}],
L_{4}^{12}[R_{4}]$ and $R_{12}$ are polygons in $\varLambda_{\text{S}}$. For an affinely regular polygon in
$\varLambda_{\text{S}}$ with $6$ vertices, consider the expansive
dilatation of $L_{6}^{12}[R_{6}]$ which is given by multiplication
with the PV-number $1+\sqrt{3}$ of (full) degree $2$ in $\oo_{12}$,
i.e., the convex $6$-gon with vertices $1+\sqrt{3}$,
$(1+\sqrt{3})\zeta_{12}$, $(1+\sqrt{3})(-1 +\zeta_{12})$,
$-(1+\sqrt{3})$, $-(1+\sqrt{3})\zeta_{12}$,
$(1+\sqrt{3})(1-\zeta_{12})$. Alternatively, simply consider the
affinely regular polygon $R_{6}$ in
$\varLambda_{\text{S}}$.

\begin{figure}
\centerline{\epsfxsize=0.55\textwidth\epsfbox{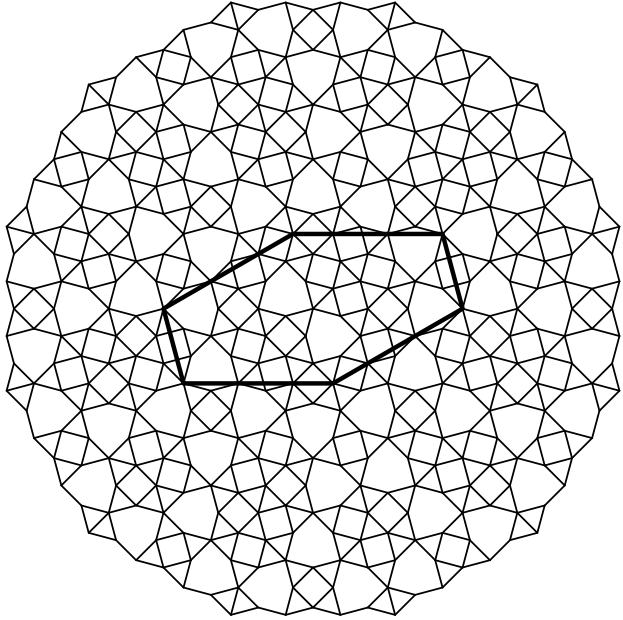}}
\caption{A central patch of the twelvefold symmetric shield
  tiling. Therein, one of the affinely regular $6$-gons as described in example {\rm (S)} is shown, the other solutions being rather obvious.}
\label{fig:s}
\end{figure} 

\end{itemize}

\section{U-Polygons in Cyclotomic Model Sets}\label{upol}

\begin{defi}\label{crossr}
Let $(t_1,t_2,t_3,t_4)$ be an ordered tuple of four distinct
elements of $\mathbbm{R}\cup\{\infty\}$. Then, its {\em cross ratio}
$(t_1,t_2;t_3,t_4)$ is defined by
$$
(t_1,t_2;t_3,t_4) := \frac{(t_3 - t_1)(t_4 - t_2)}{(t_3 - t_2)(t_4 - t_1)}\,,
$$
with the usual conventions if one of the $t_i$ equals $\infty$, so $(t_1,t_2;t_3,t_4)\in \mathbbm{R}$.
\end{defi}

\begin{defi}\label{u4nq}
Let $z\in \R^2\setminus\{0\}$, say $z=(x_{z},y_{z})^{T}$, and let $u\in
\mathbb{S}^{1}$ be a direction.
\begin{itemize}
\item[(a)]
We denote
by $s_{z}$ the slope of $z$, i.e., $s_{z}=y_{z}/x_{z}\in
\mathbbm{R}\cup\{\infty\}$.
\item[(b)]
We denote
by $u_{z}$ the direction $z/\Arrowvert z\Arrowvert\in\mathbb{S}^1$.
\item[(c)]
The {\em angle between} $u$ {\em and the positive real axis} is understood to be
the unique angle $\theta\in [0,\pi)$ having the property that a
rotation of $1\in\C$ by $\theta$ in counter-clockwise order is a
direction parallel to $u$.
\item[(d)]
For $n\,\in\,\mathbbm{N}\setminus \{1,2\}$, we denote by
$\mathcal{U}_{4,\mathbbm{Q}}^{n}$ the collection of all sets $U\subset
\mathbb{S}^{1}$ of four or more pairwise non-parallel
$\OO_{n}$-directions having the property that the cross ratio of
slopes of any four directions of $U$ (arranged in arbitrary order) is an element of $\mathbbm{Q}$.
\end{itemize}
\end{defi}

\begin{lem}\label{lodirections}
Let $u\in\mathbb{S}^1$ be an $\OO_{n}$-direction. Then, one has
$
s_u\in \mathbbm{k}_{[n,4]}\cup
\{\infty\}
$.
\end{lem}
\begin{proof}
Let $u\in
\mathbb{S}^{1}$ be an $\OO_{n}$-direction, say parallel to
$o\in\OO_{n}\setminus \{0\}$. Then, one has
\begin{equation}\label{loeq}
s_{u}=s_{o}=\frac{\frac{o-\bar{o}}{2i}}{\frac{o+\bar{o}}{2}}=-i\,\,\frac{o-\bar{o}}{o+\bar{o}}\in \left(\mathbbm{Q}(\zetan,i)\cap\R\right)\cup
\{\infty\} \,,
\end{equation}
because $\bar{o}\in\OO_{n}\setminus\{0\}$. The assertion follows from the equality
$
\mathbbm{Q}(\zetan,i)=\mathbbm{Q}(\zetan,\zeta_{4})=\mathbbm{K}_{n}\mathbbm{K}_{4}= \mathbbm{K}_{[n,4]}
$
, which implies that $\mathbbm{Q}(\zetan,i)\cap
\R=\mathbbm{k}_{[n,4]}$, since $\mathbbm{k}_{[n,4]}$ is the maximal
real subfield of $\mathbbm{K}_{[n,4]}$. 
\end{proof}

In view of Lemma~\ref{lodirections}, it is clear that the cross ratio of slopes of four pairwise non-parallel
$\OO_{n}$-directions is an element of the algebraic number field
$\mathbbm{k}_{[n,4]}$. Astonishingly, one even has the following result.

\begin{lem}\label{crkn4}
The
cross ratio of slopes of four pairwise non-parallel
$\OO_{n}$-directions is an element of the real algebraic number field
$\mathbbm{k}_{n}$.
\end{lem}
\begin{proof}
One can easily see from Equation~(\ref{loeq}) that, in the cross ratio
 of slopes of four pairwise non-parallel
$\OO_{n}$-directions, the appearing terms of the form $-i$ can be
cancelled out (even if one of the slopes equals $\infty$), hence the cross ratio
 is an element of $\Q(\zeta_{n})\cap\R=\kk_{n}$.
\end{proof}

For $n\in\{3,4,6\}$, the restriction imposed by the definition of
 the set $\mathcal{U}_{4,\mathbbm{Q}}^{n}$ is always fulfilled.

\begin{lem}\label{u4nq346}
If $n\in\{3,4,6\}$, then $\mathcal{U}_{4,\mathbbm{Q}}^{n}$ is
precisely the set of all sets $U\subset
\mathbb{S}^{1}$ of four or more pairwise non-parallel
$\OO_{n}$-directions.
\end{lem}
\begin{proof}
Observing that one has $\kk_{n}=\Q$ for $n\in\{3,4,6\}$, the
assertion follows immediately from Lemma~\ref{crkn4}.
\end{proof}

\begin{lem}\label{dirlem}
Let $n\,\in\,\mathbbm{N}\setminus \{1,2\}$ and let
$\varLambda_{n}(t,W)\in \mathcal{M}(\OO_{n})$ be a cyclotomic model
set. Then, one has: 
\begin{itemize}
\item[(a)]
The set of $\OO_{n}$-directions is precisely
the set of $\varLambda_{n}(t,W)$-directions.    
\item[(b)]
The set of $1$-dimensional $\OO_{n}$-subspaces is precisely
the set of $1$-dimensional $\varLambda_{n}(t,W)$-subspaces.    
\end{itemize}
\end{lem}
\begin{proof}
Let us start with (a). Since one has $\varLambda_{n}(t,W) - \varLambda_{n}(t,W)\subset \OO_{n}$, every $\varLambda_{n}(t,W)$-direction
is an $\OO_{n}$-direction. For the converse, let $u\in
\mathbb{S}^{1}$ be an $\OO_{n}$-direction, say parallel to
$o\in\OO_{n}\setminus \{0\}$. By Lemma~\ref{dilate}, there is a homothety $h\!:\, \C \longrightarrow
\C$ such that $h[\{0,o\}]\subset \varLambda_{n}(t,W)$. It follows that
$h(o)-h(0)\in (\varLambda_{n}(t,W) -
\varLambda_{n}(t,W))\setminus\{0\}$. Since $h(o)-h(0)$ is parallel to
$o$, the assertion follows. Part (b) follows from similar arguments.
\end{proof}

\begin{rem}
Lemma~\ref{dirlem} shows
that the notion of
$\OO_{n}$-directions in the context of cyclotomic model sets is a natural extension of the notion of {\em
  lattice directions} (i.e., $\Z^{d}$-directions, $d\geq 2$)
in~\cite{GG}.
\end{rem}

The following standard result is usually stated in the framework of projective
geometry; compare~\cite[Ch. XI, in particular, Corollary 96.11]{Borsuk}. For convenience, we give a
reformulation and also include a proof. 

\begin{lem}\label{crossratio}
Let $z_{j}\in \R^2$, $j\in \{1,\dots,4\}$, be four pairwise
non-parallel elements of the Euclidean plane with slopes $s_{z_{j}}\in \mathbbm{R}\cup\{\infty\}$. Furthermore, let $\psi$ be a non-singular linear transformation of the plane. Then, one has
$$
(s_{z_{1}},s_{z_{2}};s_{z_{3}},s_{z_{4}}) = (s_{\psi(z_1)},s_{\psi(z_2)};s_{\psi(z_3)},s_{\psi(z_4)})\,.
$$  
\end{lem}
\begin{proof}
Let $z_j=(
x_j,
y_j)^{T}$, $j\in \{1,\dots,4\}$. Then, one has
\begin{equation}\label{crpsi}
(s_{z_{1}},s_{z_{2}};s_{z_{3}},s_{z_{4}}) = \frac{(\frac{y_3}{x_3}-\frac{y_1}{x_1})(\frac{y_4}{x_4}-\frac{y_2}{x_2})}{(\frac{y_3}{x_3}-\frac{y_2}{x_2})(\frac{y_4}{x_4}-\frac{y_1}{x_1})}=\frac{\operatorname{det}\left(\begin{smallmatrix}
x_1&x_3\\
y_1&y_3
\end{smallmatrix}\right)\operatorname{det}\left(\begin{smallmatrix}
x_2&x_4\\
y_2&y_4
\end{smallmatrix}\right)}{\operatorname{det}\left(\begin{smallmatrix}
x_2&x_3\\
y_2&y_3
\end{smallmatrix}\right)\operatorname{det}\left(\begin{smallmatrix}
x_1&x_4\\
y_1&y_4
\end{smallmatrix}\right)}\,.
\end{equation}
The map $\psi\!:\,\R^2 \longrightarrow \R^2$ is given by
$z\longmapsto Az$, where $A$ is a real $2\times 2$ matrix with
non-zero determinant. The assertion follows immediately from
Equation~(\ref{crpsi}) in conjunction with the multiplication theorem for determinants.
\end{proof}

\begin{rem}\label{r4}
Let $z_{j}\in \OO_{4}=\mathbbm{Z}^2$, $j\in \{1,\dots,4\}$, be four pairwise
non-parallel elements of the square lattice. Clearly, the cross ratio
$(s_{z_{1}},s_{z_{2}};s_{z_{3}},s_{z_{4}})$
is a rational number, say $q$. Let $n\,\in\,\mathbbm{N}\setminus
\{1,2\}$ and consider the non-singular linear transformation
$L_{4}^{n}$ of the Euclidean plane as defined in the proof of
Theorem~\ref{th1} (proof of direction (vi) $\Rightarrow$ (i)). By Lemma~\ref{crossratio},
one has
$$q=(s_{L_{4}^{n}(z_1)},s_{L_{4}^{n}(z_2)};s_{L_{4}^{n}(z_3)},s_{L_{4}^{n}(z_4)})\,.$$
Hence, since $L_{4}^{n}$ maps $\mathbbm{Z}^2$ into $\OO_{n}$ (see the
proof of Theorem~\ref{th1} ((vi) $\Rightarrow$ (i)) or
Lemma~\ref{Oo}(a)), we see that $\mathcal{U}_{4,\mathbbm{Q}}^{n}$ from
Definition~\ref{u4nq}(d) is not empty. In fact, one has
$$
V:=\left\{u_{L_{4}^{n}(z_j)}\,\Big|\, j\in \{1,\dots,4\}\right
\}\in \mathcal{U}_{4,\mathbbm{Q}}^{n}\,.
$$
We shall make use
of $L_{4}^{n}$ in conjunction with Lemma~\ref{dilate} in order to
transfer results obtained by Gardner and Gritzmann~\cite{GG} for the
square lattice to the class of cyclotomic model sets.
\end{rem}

The following result was proved using Darboux's theorem~\cite{D} on
midpoint polygons; see~\cite{GM} or~\cite[Ch. 1]{G} and
compare~\cite[Lemma 4.3.6]{GG2}.

\begin{prop}\label{uaffine}
If $U\subset \mathbb{S}^{1}$ is a finite set of directions, there exists a $U$-polygon if and only if there is an affinely regular polygon such that each direction of $\,U$ is parallel to one of its edges.
\end{prop}
\begin{proof}
See \cite[Proposition 4.2]{GG}.
\end{proof}

\begin{rem}
Note that a $U$-polygon need not be affinely regular, even if it is a
$U$-polygon in a cyclotomic model set.
\end{rem}

\begin{ex}\label{ex612}
 Consider
$\varLambda_{n}(t,W)\in \mathcal{M}(\OO_{n})$, where
$n\,\in\,\{3,4,5,8,12\}$. Thereby, our standard examples
$\varLambda_{\text{SQ}}$, $\varLambda_{\text{TRI}}$, $\varLambda_{\text{TTT}}$, $\varLambda_{\text{AB}}$ and $\varLambda_{\text{S}}$ as introduced in Section~\ref{sec4} are included. If $n\,\in\,\{3,4,5,8\}$, let $U$ consist of the six pairwise
non-parallel $\OO_{n}$-directions $u_1$, $u_{2+\zetan}$,
$u_{1+\zetan}$, $u_{1+2\zetan}$, $u_{\zetan}$ and $u_{-1+\zetan}$. Let
$P$ be the non-degenerate convex dodecagon with vertices at
$3+\zetan$, $3+2\zetan$, $2+3\zetan$, $1+3\zetan$, $-1+2\zetan$,
$-2+\zetan$, and the reflections of these points in the origin $0$;
compare~\cite[Example 4.3 and Figure 1]{GG} and see
Remark~\ref{r4}. Then, one easily verifies that $P$ is a $U$-polygon
in $\OO_{n}$. By Lemma~\ref{dilate}, there is a homothety $h\!:\, \C
\longrightarrow \C$ such that $P':=
h[P]$ is a polygon in $\varLambda_{n}(t,W)$. Since $P'$ is a
$U$-polygon (see Lemma~\ref{homotu}(a)), $P'$ is a $U$-polygon
in $\varLambda_{n}(t,W)$. However, by Corollary~\ref{cor2}(a) and
Corollary~\ref{cor3}, $P'$ is not affinely regular. If $n=12$, let $U$ consist of the four pairwise non-parallel
$\OO_{12}$-directions $u_{1-\zeta_{12}}$, $u_{1}$, $u_{1+\zeta_{12}}$
and $u_{\zeta_{12}}$. Let $P$ be the non-degenerate convex
octagon with vertices at $1$, $\zeta_{12}$, $-1+\zeta_{12}$, $-2$,
$-2-\zeta_{12}$, $-1-2\zeta_{12}$, $-2\zeta_{12}$ and $1-\zeta_{12}$;
compare~\cite[Ch. 4, Figure 4.3]{HK} and see Remark~\ref{r4}. Then,
one easily verifies that $P$ is a $U$-polygon in $\OO_{12}$. By
Lemma~\ref{dilate}, there is a homothety $h\!:\, \C \longrightarrow
\C$ such that $P':= h[P]$ is a
polygon in $\varLambda_{12}(t,W)$. Since $P'$ is a $U$-polygon
(see Lemma~\ref{homotu}(a)), $P'$ is a $U$-polygon in
$\varLambda_{12}(t,W)$. However, by Corollary~\ref{cor3}, $P'$ is not affinely regular.
\end{ex}

The proof of the following result is a modified version of that of~\cite[Lemma 4.4]{GG}; compare Remark~\ref{r4}.

\begin{lem}\label{uleq3}
Let $n\,\in\,\mathbbm{N}\setminus \{1,2\}$ and let
$\varLambda_{n}(t,W)\in \mathcal{M}(\OO_{n})$ be a cyclotomic model set. If $U\subset \mathbb{S}^{1}$ is any set of three or fewer pairwise non-parallel $\OO_{n}$-directions, then there exists a $\,U$-polygon in $\varLambda_{n}(t,W)$.
\end{lem}
\begin{proof}
We begin with the case $\operatorname{card}(U)=3$. For
$i\in\{1,2,3\}$, let $\alpha_{i},\beta_{i} \in \oo_{n}$ so that
$\alpha_{i}+\beta_{i}\zetan \in \OO_{n}$, $i\in\{1,2,3\}$, are
parallel to the directions of $U$; cf. Lemma~\ref{Oo}(a). Without loss of generality, we may assume that $\alpha_{1}>\alpha_{2}>\alpha_{3}$, and that either $\beta_{1}=\beta_{2}=\beta_{3}>0$, or $\beta_{1}=0$, $\alpha_{1}>0$ and $\beta_{2}=\beta_{3}>0$. Set
$$
h:=\alpha_{2}\beta_{3}-\alpha_{3}\beta_{2}, \, k:= \alpha_{1}\beta_{3}-\alpha_{3}\beta_{1}, \, l:=\alpha_{1}\beta_{2}-\alpha_{2}\beta_{1}\,. 
$$
It turns out that $h,k,l>0$, and one can verify that the points
$0,h\alpha_{1}+h\beta_{1}\zetan$,
$(h\alpha_{1}+k\alpha_{2})+(h\beta_{1}+k\beta_{2})\zetan,(h\alpha_{1}+k\alpha_{2}+l\alpha_{3})+(h\beta_{1}+k\beta_{2}+l\beta_{3})\zetan$, 
$(k\alpha_{2}+l\alpha_{3})+(k\beta_{2}+l\beta_{3})\zetan,l\alpha_{3}+l\beta_{3}\zetan$ are the vertices of a $U$-polygon, say $P$, in $\OO_{n}$. By Lemma~\ref{dilate}, there is a  homothety $h\!:\, \C \longrightarrow \C$ such that $P':= h[P]$ is a polygon in $\varLambda_{n}(t,W)$. Since $P'$ is a $U$-polygon (see Lemma~\ref{homotu}(a)), $P'$ is a $U$-polygon in $\varLambda_{n}(t,W)$. The remaining cases $\operatorname{card}(U)<3$ follow from similar arguments.
\end{proof}

\begin{rem}
A geometric way of seeing the proof of Lemma~\ref{uleq3} is that one
first constructs a triangle in $\mathbbm{K}_n$ having sides parallel to the
given directions of $U$. If two of the vertices are chosen in
$\mathbbm{K}_n$, then the third is automatically in $\mathbbm{K}_n$. Now, fit six of these triangles together in the obvious way
to make an
affinely regular hexagon in $\mathbbm{K}_n$. The latter is then a $U$-polygon, say $P$,
in $\mathbbm{K}_n$. By Lemma~\ref{dilate}, there is a  homothety $h\!:\, \C \longrightarrow \C$ such that $P':= h[P]$ is a polygon in $\varLambda_{n}(t,W)$. Since $P'$ is a $U$-polygon (see Lemma~\ref{homotu}(a)), $P'$ is a $U$-polygon in $\varLambda_{n}(t,W)$. 
\end{rem}

The proof of the following result is step by step analogous to that of~\cite[Theorem
4.5]{GG}, wherefore we need not repeat it here. The most important tools for this proof are
Proposition~\ref{uaffine} and $p$-adic valuations. 

\begin{theorem}\label{finitesetcr}
Let $n\,\in\,\mathbbm{N}\setminus \{1,2\}$, let $U\in \mathcal{U}_{4,\mathbbm{Q}}^{n}$, and suppose the existence of a $U$-polygon. Then, one has $\operatorname{card}(U)\leq 6$, and the cross ratio of slopes of any four directions of $U$, arranged in order of increasing angle with the positive real axis, is an element of the set
$N_{1}$ $($as defined in
Theorem~$\ref{intersectq}$$)$.\qed 
\end{theorem}

Due to the roughness of our analysis in Section~\ref{sec9}, we shall not
be able in this section to prove a full analogue of
Theorem~\ref{finitesetcr} in the case of arbitrary sets of four or
more pairwise non-parallel $\OO_{n}$-directions. In fact, it is
an open problem whether there is always a cardinality statement as above.

Suppose the existence of a $U$-polygon in a cyclotomic model
set. Then, the set $U$ consists of
$\OO_{n}$-directions, where $n\in\N\setminus\{1,2\}$. The proof of the
following result is a modified version of the first part of the proof of ~\cite[Theorem 4.5]{GG}.

\begin{theorem}\label{finitesetncr0}
Let $n\,\in\,\mathbbm{N}\setminus \{1,2\}$, let $U\subset \mathbb{S}^1$ be a set of four or more pairwise
non-parallel $\OO_{n}$-directions, and suppose the existence of a
$U$-polygon. Then, the
cross ratio of slopes of any four directions of $U$, arranged in order
of increasing angle with the positive real axis, is an element of the set
$$
\Big ( \bigcup_{m\geq 4} f_{m}[D_{m}]\Big )\cap \mathbbm{k}_{n}
\,.$$
\end{theorem}
\begin{proof}
Let $U$ be as in the assertion. By Proposition~\ref{uaffine}, $U$ consists
of directions parallel to the edges of an affinely regular
polygon. Hence, there is a non-singular linear transformation $\psi$ of the
plane with the following property: 
 If one sets
$$
V:= \left\{\left. u_{\psi(u')}\, \right |\, u'\in U \right \} \subset \mathbb{S}^{1}\,,
$$
then $V$ is contained in a set of directions that are equally spaced
in $\mathbb{S}^{1}$, i.e., the angle between each pair of adjacent
directions is the same. Since the directions of $U$ are pairwise
non-parallel, there is an $m\in\N$ with $m\geq 4$ such that
each direction of $V$ is parallel to a direction of the form $e^{h\pi i/m}$, where
$h\in\N_{0}$ satisfies $h\leq m-1$. Let $u'_{j}$, $1\leq
j\leq 4$, be four
directions of $U$, arranged in order
of increasing angle with the positive real axis. By
Lemma~\ref{crkn4}, the cross ratio of the slopes of
these $\OO_{n}$-directions, say
$q:=(s_{u'_{1}},s_{u'_{2}};s_{u'_{3}},s_{u'_{4}})$, is an element of the
real algebraic
number field
$\mathbbm{k}_{n}$. One can see by Lemma~\ref{crossratio} together
with the fact that every
non-singular linear transformation of the plane either
preserves or inverts orientation that we may assume, without
loss of generality, that each direction $u_{\psi(u'_{j})}\in V$ is parallel to a direction of the form $e^{h_{j}\pi i/m}$, where $h_j \in\N_{0}$,
$1\leq j \leq 4$, and, moreover, $h_1<h_2<h_3<h_4\leq m-1$. Using
Lemma~\ref{crossratio} again, one gets
$$
q=(s_{\psi(u'_1)},s_{\psi(u'_2)};s_{\psi(u'_3)},s_{\psi(u'_4)})=\frac{(\tan (\frac{h_3 \pi}{m})- \tan (\frac{h_1 \pi}{m}))(\tan (\frac{h_4 \pi}{m})- \tan (\frac{h_2 \pi}{m}))}{(\tan (\frac{h_3 \pi}{m})- \tan (\frac{h_2 \pi}{m}))(\tan (\frac{h_4 \pi}{m})- \tan (\frac{h_1 \pi}{m}))}\,.
$$
Manipulating the right-hand side, one obtains
$$
q=\frac{\sin(\frac{(h_3-h_1)\pi}{m})\sin(\frac{(h_4-h_2)\pi}{m})}{\sin(\frac{(h_3-h_2)\pi}{m})\sin(\frac{(h_4-h_1)\pi}{m})}\,.
$$
Setting $k_1:=h_3-h_1$, $k_2:=h_4-h_2$, $k_3:=h_3-h_2$ and
$k_4:=h_4-h_1$, one gets $1\leq k_3<k_1, k_2<k_4\leq m-1$ and
$k_1+k_2=k_3+k_4$.

Using $\sin(\theta)=\frac{-e^{-i\theta}(1-e^{2i\theta})}{2i}$, one
obtains
$$
\mathbbm{k}_{n}\owns q=\frac{(1-\zeta_{m}^{k_1})(1-\zeta_{m}^{k_2})}{(1-\zeta_{m}^{k_3})(1-\zeta_{m}^{k_4})}=f_m(d)\,,
$$
with $d:=(k_1,k_2,k_3,k_4)$, as in~(\ref{fmd}). Then, $d\in D_m$ if its
first two coordinates are interchanged, if necessary, to ensure that
$k_1\leq k_2$; note that this operation does not change the value of
$f_m(d)$. This completes the proof.
\end{proof}

Suppose that there is a $U$-polygon in an (aperiodic) cyclotomic model
set with co-dimension two. Then, the set $U$ consists of
$\OO_{n}$-directions, where $n\in\{5,8,10,12\}$. 

\begin{theorem}\label{finitesetncr}
Let $n\in\{5,8,10,12\}$, let $U\subset \mathbb{S}^1$ be a set of four or more pairwise
non-parallel $\OO_{n}$-directions, and suppose the existence of a
$U$-polygon. Then, the
cross ratio of slopes of any four directions of $U$, arranged in order
of increasing angle with the positive real axis, maps under the norm
$N_{\mathbbm{k}_{n}/\Q}$ to the set $N_{2}$ as defined in
Theorem~$\ref{intersectk8}$.
\end{theorem}
\begin{proof}
This is an immediate consequence of Theorem~\ref{finitesetncr0} in
conjunction with Corollary~\ref{intersectk510s}.
\end{proof}

For the general case $n\in\N\setminus\{1,2\}$, one obtains the
following results.

\begin{theorem}\label{finitesetncrgeneral}
For all $e\in\operatorname{Im}(\phi/2)$, there is a finite set $N_{e}\subset \Q$
such that, for all $n\in(\phi/2)^{-1}[\{e\}]$, and all sets $U\subset \mathbb{S}^1$ of four or more pairwise
non-parallel $\OO_{n}$-directions, one has the following:

If there
exists a
$U$-polygon, then the
cross ratio of slopes of any four directions of $U$, arranged in order
of increasing angle with the positive real axis, maps under the norm
$N_{\mathbbm{k}_{n}/\Q}$ to $N_{e}$.
\end{theorem}
\begin{proof}
This is an immediate consequence of Theorem~\ref{finitesetncr0} in
conjunction with Corollary~\ref{intersectk510sg}.
\end{proof}

\begin{theorem}\label{finitesetncrgeneral2}
For all $n\in\N\setminus\{1,2\}$ and all sets $U\subset \mathbb{S}^1$ of four or more pairwise
non-parallel $\OO_{n}$-directions, one has the following:

If there
exists a
$U$-polygon, then the
cross ratio of slopes of any four directions of $U$, arranged in order
of increasing angle with the positive real axis, maps under the norm
$N_{\mathbbm{k}_{n}/\Q}$ to the set
$$\pm\left(\{1\}\cup \big[\Q_{>1}^{\mathbbm{P}_{\leq 2}}\big]^{\pm 1}\right)\,.$$
\end{theorem}
\begin{proof}
This is an immediate consequence of Theorem~\ref{finitesetncr0} in
conjunction with Theorem~\ref{intersectk8generalcoro}.
\end{proof}

\section{Discrete Tomography of Cyclotomic Model Sets}\label{secset}

\subsection{Setting}

For the discrete tomography of aperiodic cyclotomic model sets,
one additional difficulty, in comparison to the crystallographic case,
stems from the fact that it is not sufficient to consider one pattern
and its translates to define the setting. Hence, to define the
analogue of a specific crystal, one has to add all infinite patterns that
emerge as limits of sequences of translates defined in the local
topology (LT). Here, two patterns are $\varepsilon$-close if, after a
translation by a distance of at most $\varepsilon$, they agree on a ball of
radius $1/\varepsilon$ around the origin. If the starting
pattern $P$ is crystallographic, no new patterns are added; but if $P$
is a generic
aperiodic cyclotomic model set,
one ends up with uncountably many different patterns, even up to
translations! Nevertheless, all of them are locally indistinguishable (LI).
This means that every {\em finite} patch in $\varLambda$ also appears in 
any of the other elements of the LI-class and vice versa; see~\cite{B} for details.

\begin{rem}\label{LI}
The entire LI-class of a regular, generic cyclotomic model set $\varLambda(W)$
can be shown to consist of all sets $t+\varLambda(\tau+W)$, with
$t \in \mathbbm{R}^{2}$ and $\tau\in(\R^2)^{\frac{\phi(n)}{2}-1}$ such that $\partial(\tau+ W)\cap
[\OO_{n}]^{\star_{n}}=\varnothing$ (i.e., $\tau$ is in a generic position), 
and all patterns obtained as limits of
sequences $t+\varLambda(\tau_{n}+W)$, with all $\tau_{n}$
in a generic position; see~\cite{B,Schl}. Each such limit is then a {\em subset} of some
$t+\varLambda(\tau+W)$, as $\overline{W}$ was assumed compact, but $\tau$ might not be in a generic
position. In view of this complication, we must make sure that we
deal with finite subsets of {\em generic} cyclotomic model sets.  This
restriction to the generic case is the proper analogue of the
restriction to {\em perfect} lattices and their translates in the
classical setting.
\end{rem}

\begin{defi}\label{W-sets}
Let $n\in\N\setminus\{1,2\}$, let $W\subset
(\R^2)^{\phi(n)/2-1}$ be a window $($cf. Definition~$\ref{cyclodef}$$)$, 
and let a star map $.^{\star_{n}}$ be given, i.e., a map $.^{\star_{n}}\! : \, \mathcal{O}_{n}\longrightarrow
(\R^2)^{\phi(n)/2-1}$, given by $z\longmapsto 0$, if
$n\in\{3,4,6\}$, and given by $z\longmapsto
(\sigma_2(z),\dots,\sigma_{\phi(n)/2}(z))$ otherwise $($as
described in Definition~$\ref{cyclodef}$$)$. We let
$\mathcal{M}_{g}(\mathcal{O}_{n})$ be the set of generic elements of $\mathcal{M}(\mathcal{O}_{n})$ and define
$$
W^{\star_{n}}_{\mathcal{M}_{g}(\mathcal{O}_{n})}:=\left\{\varLambda_{n}^{\star_{n}}(t,\tau
+ W)\,\left |\,\begin{array}{l} t\in\R^{2},\tau \in
    (\R^2)^{\frac{\phi(n)}{2}-1}\mbox{ and}\\
    \left[\OO_{n}\right]^{\star_{n}} \cap \, \partial (\tau+W)
    =\varnothing \end{array}\right. \right \}
\subset \mathcal{M}_{g}(\mathcal{O}_{n})\,,
$$ 
 where the elements $\varLambda_{n}^{\star_{n}}(t,\tau+ W)$ of this set are understood to be defined by use of the
 above star map $.^{\star_{n}}$, i.e.,
 $$\varLambda_{n}^{\star_{n}}(t,\tau+
 W)=t+\left\{z\in\OO_{n}\,\left |\,z^{\star_{n}}\in\tau+W \right. \right\}\,.$$ 
\end{defi}

\begin{rem}\label{remwsets}
Let $n\in\N\setminus\{1,2\}$. Note that, if $W\subset
(\R^2)^{\phi(n)/2-1}$ is a window, then every translate of it, i.e., 
$\tau+W$, is a window as well. Note further that for $n=4$ (resp.,
$n\in\{3,6\}$) the set
$W^{\star_{n}}_{\mathcal{M}_{g}(\mathcal{O}_{n})}$ simply consists of
all translates of the square lattice $\OO_{4}$ (resp., triangular
lattice $\OO_{3}=\OO_{6}$).
\end{rem}

The setting for the uniqueness problem of discrete tomography of cyclotomic model sets
 looks as follows. Let $n\,\in\,\mathbbm{N}\setminus \{1,2\}$, let $W\subset
 (\R^2)^{\phi(n)/2-1}$ be a window, 
 and let a star map $.^{\star_{n}}$ be given (as described in
 Definition~\ref{cyclodef}). Then, we are interested in the
 (successive) determination of the set $$\bigcup_{\varLambda\in
   W^{\star_{n}}_{\mathcal{M}_{g}(\mathcal{O}_{n})}}\mathcal{F}(\varLambda)$$ or suitable subsets thereof by the $X$-rays in a small number of $\OO_{n}$-directions. 

\subsection{Grids}\label{grids}

It is advantageous to narrow down the set of possible positions in a
way that matches the algebraic setting at hand.

\begin{defi}\label{defgrid}
Let $U\subset\mathbb{S}^1$
be a finite set of pairwise
non-parallel directions. We define the {\em complete grid}
$G_{\OO_{n}}^{U}$ of $\OO_{n}$ with respect to $U$ as
$$G_{\OO_{n}}^{U} := \bigcap_{u\in U}\,\,\Big( \bigcup_{\ell \in
  \mathcal{L}_{\OO_{n}}^{u}} \ell\Big)\,,$$
where $\mathcal{L}_{\OO_{n}}^{u}$ denotes the set of lines in
direction $u$ in $\R^2$, which pass through at least one element of $\OO_{n}$.
 Further, for a
finite subset   $F$ of $\R^2$, we define the \emph{grid} of $F$ with respect to the $X$-rays in the
directions of $U$ as
$$
G_{F}^{U}\,\,:=\,\,\bigcap_{u\in U}\,\,\left( \bigcup_{\ell \in
  \mathrm{supp}(X_{u}F)} \ell\right)\,.
$$
\end{defi}

\begin{rem}
Note that, in the situation of Definition~\ref{defgrid}, one has the
inclusion \begin{equation}\label{triveq}\OO_{n}\subset G_{\OO_{n}}^{U}\,.\end{equation}
\end{rem}

The following property is immediate.

\begin{lem}\label{fgrid}
If $U\subset\mathbb{S}^1$
is a finite set of pairwise
non-parallel directions, then for all finite subsets $F,F'$ of
$\R^2$ one has
$$
(X_{u}F=X_{u}F'\;\forall u \in U) \;  \Longrightarrow\; F,F'\,\,\subset\,\,G_{F}^{U}\,\,=\,\,G_{F'}^{U}\,.\hspace*{5.675em}\qed\hspace*{-5.675em}
$$
\end{lem}

Naturally, we are interested in the following special situation. Let $n\in\mathbbm{N}\setminus\{1,2\}$, let $U\subset\mathbb{S}^1$
be a finite set of pairwise
non-parallel $\OO_{n}$-directions, and let $F$ be a
finite subset of $t+\OO_{n}$, where $t\in\R^2$. Clearly, the grid
$G_{F}^{U}$ of $F$ with respect to the $X$-rays in the
directions of $U$ may be a proper superset of $F$. In fact, it may
even contain points that lie in a {\em different} 
translate of $\OO_{n}$ than $F$ itself. This problem was studied
in~\cite[Fig. 5, Sec. 4.2 and Sec. 5.1]{BG2}. Here, we are interested
in the case 
where the latter phenomenon cannot occur. The main result of this
section (Theorem~\ref{oneequiv} below) will show that this
goal can be reached by allowing only a special kind of sets $U$ of at
least two $\OO_{n}$-directions. First, we need the following results.

\begin{prop}\label{gridprop}
Let $n\in\mathbbm{N}\setminus\{1,2\}$ and let $o,o'$ be two
non-parallel elements of $\mathcal{O}_{n}$. Then, the complete grid
 $G_{\OO_{n}}^{\{u_{o},u_{o'}\}}$ of $\OO_{n}$ with respect to the
directions $u_{o}$ and $u_{o'}$ $($as defined in
Definition~$\ref{u4nq}$$({\rm b})$$)$ satisfies
$$
\mathcal{O}_{n}\,\,\subset\,\,G_{\OO_{n}}^{\{u_{o},u_{o'}\}}\,\,\subset\,\,
   \big\langle \left\{o/(\alpha_{o} \beta_{o'} - \beta_{o} \alpha_{o'}),
       o'/(\alpha_{o} \beta_{o'} - \beta_{o} \alpha_{o'}) \right\}\big
   \rangle_{\thinspace\scriptstyle{\mathcal{O}}\displaystyle_{n}}\,\,\subset\,\,\mathbbm{K}_{n}\,\,\subset\,\,\C
\,, 
$$
where the elements $\alpha_{o},\alpha_{o'}, 
\beta_{o}, \beta_{o'} \in
\thinspace\scriptstyle{\mathcal{O}}\displaystyle_{n}$ are determined
by $o = \alpha_{o} + \beta_{o} \zeta_{n}$ and $o' = \alpha_{o'} +
\beta_{o'} \zeta_{n}$ $($cf. Lemma~$\ref{Oo}$$({\rm a})$$)$, and $\langle
\,.\,\rangle_{\thinspace\scriptstyle{\mathcal{O}}\displaystyle_{n}}$
denotes the $\thinspace\scriptstyle{\mathcal{O}}\displaystyle_{n}$-linear hull.
\end{prop}
\begin{proof}
See~\cite[Proposition 5 and Remark 14]{BG2}.
\end{proof}

\begin{rem}\label{xyzneq0}
The linear independence of $\{o, o'\}$ and $\{1, \zeta_{n}\}$ over
$\R$ implies that $\alpha_{o} \beta_{o'} - \beta_{o} \alpha_{o'} \neq
0$.
\end{rem}

\begin{defi}
For elements $o\in\OO_{n}\setminus\{0\}$ and $o'\in\OO_{n}$, we say that $o$ {\em divides} $o'$ and
write $o|o'$ if $\frac{o'}{o}\in\OO_{n}$. 
\end{defi}

\begin{defi}\label{rib}
Let $\mathbbm{K}/\mathbbm{k}$ be an extension of
algebraic number fields, say of degree
$d:=[\mathbbm{K}:\mathbbm{k}]\in\N$. Further, let $\OO_{\mathbbm{K}}$
(resp., $\oo_{\mathbbm{k}}$) be the ring of integers of $\mathbbm{K}$
(resp., $\mathbbm{k}$). Then, a subset $\{o_1,\dots,o_d\}\subset
\OO_{\mathbbm{K}}$ is called a {\em relative integral basis} of
$\mathbbm{K}/\mathbbm{k}$ if it is an $\oo_{\mathbbm{k}}$-basis of the
$\oo_{\mathbbm{k}}$-module $\OO_{\mathbbm{K}}$, i.e., if every element
$o\in \OO_{\mathbbm{K}}$ is uniquely expressible as an
$\oo_{\mathbbm{k}}$-linear combination of $\{o_1,\dots,o_d\}$.
\end{defi}

\begin{rem}\label{oo'rem}
By Lemma~\ref{Oo}(b), for $n\in\mathbbm{N}\setminus\{1,2\}$, one has  
$[\mathbbm{K}_{n}:\mathbbm{k}_{n}]=2$. Moreover, by Proposition~\ref{p1}, one has the
identities $\OO_{\mathbbm{K}_{n}}=\OO_{n}$
and $\oo_{\mathbbm{k}_{n}}=\oo_{n}$. Let $o,o'$ be two
non-parallel elements of $\mathcal{O}_{n}$. Then, $\{o,o'\}$ is a
relative integral basis of $\mathbbm{K}_{n}/\mathbbm{k}_{n}$ if and
only if the $\oo_{n}$-linear hull $\langle\{o,o'\}\rangle_{\oo_{n}}$
of $\{o,o'\}$ equals $\OO_{n}$. 
\end{rem}

\begin{prop}\label{proponeequiv}
Let $n\in\mathbbm{N}\setminus\{1,2\}$ and let $o,o'$ be two
non-parallel elements of $\mathcal{O}_{n}$, say
$o=\alpha_{o}+\beta_{o}\zeta_{n}$ and
$o'=\alpha_{o'}+\beta_{o'}\zeta_{n}$ for uniquely determined
$\alpha_{o},\alpha_{o'},\beta_{o},\beta_{o'}\in\oo_{n}$ 
$($cf. Lemma $\ref{Oo}$$({\rm a})$$)$. The following statements are equivalent:
\begin{itemize}
\item[(i)]
$\alpha_o\beta_{o'}-\beta_{o}\alpha_{o'}\in\oo_{n}^{\times}$.
\item[(ii)]
$\alpha_o\beta_{o'}-\beta_{o}\alpha_{o'}\big|o$ and $\alpha_o\beta_{o'}-\beta_{o}\alpha_{o'}\big|o'$.
\item[(iii)]
$\{o,o'\}$ is a relative integral basis of 
$\mathbbm{K}_{n}/\mathbbm{k}_{n}$.
\end{itemize}
Moreover, each of the above conditions \rm{(i)}-\rm{(iii)} implies the equation
$$
G_{\OO_{n}}^{\{u_o,u_{o'}\}}=\OO_{n}\,.
$$ 
\end{prop}
\begin{proof}
Direction (i) $\Rightarrow$ (ii) is immediate. Next, we show direction (ii)
$\Rightarrow$ (i). Suppose that
$\alpha_o\beta_{o'}-\beta_{o}\alpha_{o'}\big|o$ and
$\alpha_o\beta_{o'}-\beta_{o}\alpha_{o'}\big|o'$, say
$(\alpha_o\beta_{o'}-\beta_{o}\alpha_{o'})(\gamma+\delta \zeta_{n})=o$
and $(\alpha_o\beta_{o'}-\beta_{o}\alpha_{o'})(\gamma'+\delta'
\zeta_{n})=o'$ for suitable $\gamma,\gamma',\delta,\delta'\in\oo_{n}$. Using the $\R$-linear
independence of $\{1,\zeta_{n}\}$, the latter implies the equations 
\begin{eqnarray*}(\alpha_o\beta_{o'}-\beta_{o}\alpha_{o'})\gamma&=&\alpha_{o}\,,\\
(\alpha_o\beta_{o'}-\beta_{o}\alpha_{o'})\delta&=&\beta_{o}\,,\\
(\alpha_o\beta_{o'}-\beta_{o}\alpha_{o'})\gamma'&=&\alpha_{o'}\,,\\
(\alpha_o\beta_{o'}-\beta_{o}\alpha_{o'})\delta'&=&\beta_{o'}\,.\\
\end{eqnarray*} 
Consequently, one has
$$
\alpha_o\beta_{o'}-\beta_{o}\alpha_{o'}=(\alpha_o\beta_{o'}-\beta_{o}\alpha_{o'})^2
\,(\gamma\delta'-\delta\gamma')
\,.
$$
Further, dividing by $\alpha_o\beta_{o'}-\beta_{o}\alpha_{o'}$
(cf. Remark~\ref{xyzneq0}), one obtains the equation
$$
1=(\alpha_o\beta_{o'}-\beta_{o}\alpha_{o'})\,(\gamma\delta'-\delta\gamma')
\,,
$$
and the assertion follows. For direction (i)
$\Rightarrow$ (iii), first observe that Proposition~\ref{gridprop} particularly shows that  
$$
\mathcal{O}_{n}\,\,\subset\,\,
   \big\langle \left\{o/(\alpha_{o} \beta_{o'} - \beta_{o} \alpha_{o'}),
       o'/(\alpha_{o} \beta_{o'} - \beta_{o} \alpha_{o'}) \right\}\big
   \rangle_{\thinspace\scriptstyle{\mathcal{O}}\displaystyle_{n}}\,\,\subset\,\,\OO_{n}
\,,
$$
hence
$$\left\langle \left\{o/(\alpha_{o} \beta_{o'} - \beta_{o} \alpha_{o'}),
       o'/(\alpha_{o} \beta_{o'} - \beta_{o} \alpha_{o'}) \right\}\right
   \rangle_{\thinspace\scriptstyle{\mathcal{O}}\displaystyle_{n}}\,\,=\,\,\langle\{o,o'\}\rangle_{\oo_{n}}=\OO_{n}
\,.
$$
By Remark~\ref{oo'rem}, the assertion follows. Finally, let us prove direction (iii) $\Rightarrow$ (i). Here, by Remark~\ref{oo'rem}, there are 
$\gamma,\gamma',\delta,\delta'\in\oo_{n}$ such that $\gamma o+\delta
o'=1$ and $\gamma' o+\delta' o'=\zeta_{n}$. Using the $\R$-linear
independence of $\{1,\zeta_{n}\}$, the latter implies that $\gamma
\alpha_{o}+\delta \alpha_{o'}=1$, $\gamma
\beta_{o}+\delta \beta_{o'}=0$, $\gamma'
\alpha_{o}+\delta' \alpha_{o'}=0$ and $\gamma'
\beta_{o}+\delta' \beta_{o'}=1$. Hence, one has
\begin{eqnarray*}
1&=&(\gamma
\alpha_{o}+\delta \alpha_{o'})(\gamma'
\beta_{o}+\delta' \beta_{o'})-(\gamma
\beta_{o}+\delta \beta_{o'})(\gamma'
\alpha_{o}+\delta' \alpha_{o'})\\ 
&=&(\alpha_o\beta_{o'}-\beta_{o}\alpha_{o'})(\gamma\delta'-\gamma'\delta)\,,
\end{eqnarray*}
and the assertion follows. The additional statement follows immediately from
Proposition~\ref{gridprop}, which in this case shows that  
$$
\mathcal{O}_{n}\,\,\subset\,\,G_{\OO_{n}}^{\{u_{o},u_{o'}\}}\,\,\subset\,\,
   \big\langle \left\{o/(\alpha_{o} \beta_{o'} - \beta_{o} \alpha_{o'}),
       o'/(\alpha_{o} \beta_{o'} - \beta_{o} \alpha_{o'}) \right\}\big
   \rangle_{\thinspace\scriptstyle{\mathcal{O}}\displaystyle_{n}}\,\,\subset\,\,\OO_{n}
\,.
$$
This completes the proof.
\end{proof}

\begin{theorem}\label{oneequiv}
Let $n\in\mathbbm{N}\setminus\{1,2\}$, let $O\subset
\OO_{n}\setminus\{0\}$ be a finite set of $m\geq 2$ pairwise
non-parallel elements. Suppose the existence of two different
elements $o,o'\in O$ satisfying one of the equivalent conditions {\rm
  (i)}-{\rm (iii)} 
of Proposition~$\ref{proponeequiv}$. Then, setting $U_{O}:=\{u_{o}\,|\,o\in
O\}\subset \mathbb{S}^1$, for all $t\in\R^2$ and all finite subsets $F$ of $t+\OO_{n}$, one
has the inclusion $G_{F}^{U_{O}}\subset t+\OO_{n}$.
\end{theorem}
\begin{proof}
Let $t\in\R^2$ and let $F$ be a finite subset of $t+\OO_{n}$. It
follows that $F-t$ is a finite subset of $\OO_{n}$. Clearly, one has
the inclusion
$
G_{F-t}^{U_{O}}\subset\,\,G_{\OO_{n}}^{U_{O}}\subset
G_{\OO_{n}}^{\{u_{o},u_{o'}\}}
$. Hence, by Proposition~\ref{gridprop} in conjunction with the
additional statement of Proposition~\ref{proponeequiv}, one further obtains
$$
G_{F-t}^{U_{O}}\,\subset\,G_{\OO_{n}}^{U_{O}}\,\subset\,G_{\OO_{n}}^{\{u_{o},u_{o'}\}}\,=\,\OO_{n}
\,.$$ The equality $G_{F-t}^{U_{O}}=G_{F}^{U_{O}}-t$ completes the proof.
\end{proof}

\subsection{Practical relevance}

The above setting for the discrete tomography of cyclotomic model sets
is motivated by the practice of quantitative HRTEM. This is due to the
fact that, because of the $N(n)$-fold cyclic symmetry of genuine
planar (quasi)crystals (with $N(n)$ being the function from
(\ref{eq})), 
the determination of the rotational orientation of a
(quasi)crystalline probe in an electron microscope can 
rather easily be done in the diffraction mode, prior to taking images
in the high-resolution mode, though, in general, a natural choice of a
translational origin is {\em not} possible. Therefore, in order to
prove practically relevant and rigorous results, one has to deal
with the `non-anchored' case of the whole LI-class of a regular, generic cyclotomic model set
$\varLambda$, rather than dealing with the `anchored' case of a fixed
such $\varLambda$. It will turn out in the following that the treatment of the
`non-anchored' case is often feasible. Moreover, in
Section~\ref{convdet}, we shall even be able to provide positive
uniqueness results for which we can make sure that all the
$\OO_{n}$-directions used correspond to dense lines in the
corresponding cyclotomic model sets, the latter meaning that the
resolution coming from these directions is rather high. Hence, we
believe that these results look promising in view of real applications.  

\section{Simple Uniqueness Results}

In this section, we present some uniqueness results which are meant to motivate the next three sections.

\begin{prop}\label{irratdir}
Let $n\,\in\,\mathbbm{N}\setminus \{1,2\}$ and let
$\varLambda_{n}(t,W)\in \mathcal{M}(\OO_{n})$ be a cyclotomic model set. Further, let $u \in \mathbb{S}^{1}$ be a non-$\OO_{n}$-direction. Then, $\mathcal{F}(\varLambda_{n}(t,W))$ is determined by the single $X$-ray in direction $u$.
\end{prop}
\begin{proof}
This follows immediately from the fact that a line in the plane in a
non-$\OO_{n}$-direction passes through at most one point of
$t+\OO_{n}$. 
\end{proof}

\begin{rem}
Let $n\,\in\,\mathbbm{N}\setminus \{1,2\}$. By Lemma~\ref{lodirections}, the slope of
any $\OO_{n}$-direction is an element of the set $\mathbbm{k}_{[n,4]}\cup
\{\infty\}$. It follows that any direction having any other slope is a
non-$\OO_{n}$-direction.  Of course, there are
non-$\OO_{n}$-directions, even non-$\OO_{n}$-directions having algebraic slopes. 
\end{rem}

The following result represents a fundamental source of difficulties in discrete tomography. There exist several versions; compare~\cite[Theorem 4.3.1]{HK} and~\cite[Lemma 2.3.2]{G}. 

\begin{prop}\label{source}
Let $n\,\in\,\mathbbm{N}\setminus \{1,2\}$ and let
$\varLambda_{n}(t,W)\in \mathcal{M}(\OO_{n})$ be a cyclotomic model set. Further, let $U\subset \mathbb{S}^{1}$ be an arbitrary, but fixed finite set of pairwise non-parallel $\OO_{n}$-directions. Then, $\mathcal{F}(\varLambda_{n}(t,W))$ is not determined by the $X$-rays in the directions of $U$.
\end{prop}
\begin{proof}
We argue by induction on $\operatorname{card}(U)$. The case
$\operatorname{card}(U)=0$ means $U=\varnothing$ and is obvious. Suppose the assertion to be
true whenever $\operatorname{card}(U)=k\in \mathbbm{N}_{0}$ and let
$\operatorname{card}(U)=k+1$. By induction hypothesis, there are
different elements $F$ and $F'$ of $\mathcal{F}(\varLambda_{n}(t,W))$
with the same $X$-rays in the directions of $U'$, where $U'\subset U$
satisfies $\operatorname{card}(U')=k$. Let $u$ be the remaining
direction of $U$. Choose a non-zero element $o\in \OO_{n}$ parallel to
$u$ such that $o+(F\cup F')$ and $F\cup F'$ are disjoint. Then,
$F'':= F\cup (o+F')$ and
$F''':= F'\cup (o+F)$ are different
elements of $\mathcal{F}(t+\OO_{n})$ with the same $X$-rays in the
directions of $U$. By Lemma~\ref{dilate}, there is a homothety $h\!:\,
\C \longrightarrow \C$ such that $h[F''\cup F''']=h[F'']\cup
h[F''']\subset \varLambda_{n}(t,W)$. It follows that $h[F'']$ and
$h[F''']$ are different elements of $\mathcal{F}(\varLambda_{n}(t,W))$
with the same $X$-rays in the directions of $U$; see Lemma~\ref{homotu}(b). 
\end{proof}

\begin{figure}
\centerline{\epsfysize=0.55\textwidth\epsfbox{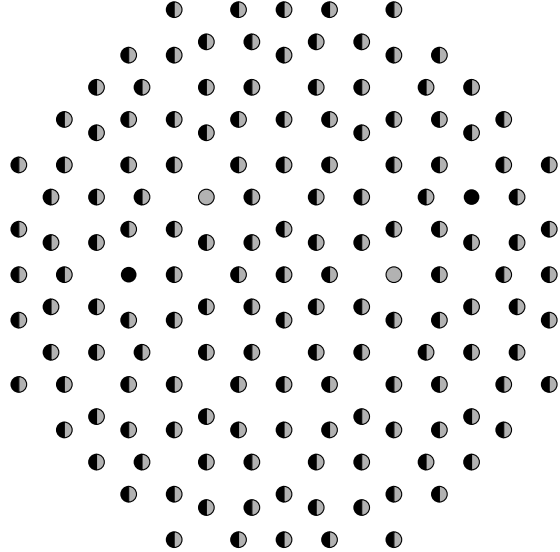}}
\caption{Two contiguous subsets of $\varLambda_{\rm AB}$ with the same
  $X$-rays in the $\OO_{8}$-directions $1$ and $\zeta_{8}$.}
\label{fig:contig}
\end{figure}

\begin{rem}\label{source2}
An analysis of the proof of Proposition~\ref{source} shows that, for
any finite set $U\subset \mathbb{S}^{1}$ of $k$ pairwise non-parallel
$\OO_{n}$-directions, there are disjoint elements  $F$ and $F'$ of
$\mathcal{F}(\varLambda_{n}(t,W))$ with
$\operatorname{card}(F)=\operatorname{card}(F')=2^{(k-1)}$ and with the same
$X$-rays in the directions of $U$. Consider any convex set $C$ in
$\R^2$ which contains $F$ and $F'$ from above. Then, the subsets
$F_1:=(C\cap\varLambda_{n}(t,W))\setminus F$ and
$F_2:=(C\cap\varLambda_{n}(t,W))\setminus F'$ of
$\mathcal{F}(\varLambda_{n}(t,W))$ also have the same
$X$-rays in the directions of $U$. Whereas the points in $F$ and $F'$
are widely dispersed over a region, those in $F_1$ and $F_2$ are
contiguous in a way similar to atoms in a quasicrystal. This procedure
is illustrated in Figure~\ref{fig:contig} in the case of the
 aperiodic cyclotomic model set $\varLambda_{\rm AB}$ associated with
the Amman-Beenker tiling as described in
  example {\rm (AB)} in Section~\ref{sec4}; compare~\cite[Remark 4.3.2]{GG2}. 
\end{rem}

The proof of the following result is the same as that of~\cite[Theorem
4.3.3]{HK}. Originally, the proof is due to R\'enyi;
see~\cite{Re}. For clarity, we prefer to repeat the details here, in a
slightly modified way.

\medskip
\begin{prop}\label{m+1}
Let $n\,\in\,\mathbbm{N}\setminus \{1,2\}$ and let
$\varLambda_{n}(t,W)\in \mathcal{M}(\OO_{n})$ be a cyclotomic model set. Further, let $U\subset \mathbb{S}^{1}$ be any set of $k+1$ pairwise non-parallel $\OO_{n}$-directions where $k\in \mathbbm{N}_{0}$. Then, $\mathcal{F}_{\leq k}(\varLambda_{n}(t,W))$ is determined by the $X$-rays in the directions of $U$.
\end{prop}
\begin{proof}
Let $F,F'\in \mathcal{F}_{\leq k}(\varLambda_{n}(t,W))$ have the same $X$-rays in the directions of $U$. Then, one has $\operatorname{card}(F)=\operatorname{card}(F')$ by Lemma~\ref{cardinality}(a) and $$
F,F'\,\,\subset\, \, G_{F}^{U}$$ by Lemma~\ref{fgrid}. But we have $G_{F}^{U}=F$ since the existence of a point in $G_{F}^{U}\setminus F$ implies the existence of at least $\operatorname{card}(U)\geq k+1$ points in $F$, a contradiction. It follows that $F=F'$.
\end{proof}

\begin{rem}\label{exdodecagon}
Let $n\,\in\,\mathbbm{N}\setminus \{1,2\}$ and let
$\varLambda_{n}(t,W)\in \mathcal{M}(\OO_{n})$ be a cyclotomic model set. Remark~\ref{source2} and Proposition~\ref{m+1} show
that $\mathcal{F}_{\leq k}(\varLambda_{n}(t,W))$ can be determined by
the $X$-rays in
any set of $k+1$ pairwise non-parallel $\OO_{n}$-directions but not by
$1+\lfloor\log_{2}k\rfloor$ pairwise non-parallel $X$-rays in
$\OO_{n}$-directions. 
\end{rem}

\section{Determination of  Subsets with Bounded Diameter by $X$-Rays or
  Projections -- Delone Sets of Finite Local Complexity with
  Applications to Meyer Sets and Model Sets}

Let $d\in\N$ and let $R>0$. Recall that a Delone set
$\varLambda\subset \R^{d}$ is uniformly discrete (i.e., there is a
radius $r>0$ such that every ball of the form $B_{r}(x)$, where $x\in\R^{d}$,
contains at most one point of $\varLambda$) and relatively
dense (i.e., there is a
radius $R>0$ such that every ball of the form $B_{R}(x)$, where $x\in\R^{d}$,
contains at least one point of $\varLambda$). The uniform
discreteness of Delone sets $\varLambda$ immediately implies the inclusion $\mathcal{D}_{<
  R}(\varLambda)\subset \mathcal{F}(\varLambda)$.

\begin{lem}\label{dirdense}
Let $d\geq 2$ and let $\varLambda\subset \R^d$ be relatively dense. Then,
the set of $\varLambda$-directions is dense in $\mathbb{S}^{d-1}$. 
\end{lem}
\begin{proof}
We may assume, without loss of generality, that $0\in\varLambda$. Let $u\in\mathbb{S}^{d-1}$ and let
$B_{\varepsilon}(u)\cap \,\mathbb{S}^{d-1}$ be an arbitrary open
$\varepsilon$-neighbourhood of $u$ in $\mathbb{S}^{d-1}$. Without
restriction, let $\varepsilon < 1$. Then,
$B_{\varepsilon}(u)\cap \,\mathbb{S}^{d-1}$ is a $(d-1)$-dimensional open ball, i.e.,
homeomorphic to the open ball $B_{1}(0):=\{x\in\R^{d-1}\,|\,\Arrowvert
x\Arrowvert < 1\}$. Consider the smallest convex cone $C$ in $\R^d$ with apex $0$ and
 containing the set $B_{\varepsilon}(u)\cap \,\mathbb{S}^{d-1}$, i.e.,
$$
C:=\left\{\left. \sum_{j=1}^{n}\lambda_j x_j \,\right | \,n\in\N,
    \R\owns \lambda_1,\dots,\lambda_n \geq 0, x_1,\dots,x_n\in
    B_{\varepsilon}(u)\cap \,\mathbb{S}^{d-1} \right\}\,. 
$$  
Since $\varLambda$ is relatively dense, there is a radius $R>0$ such that every open ball $B_{R}(z)$, where
$z\in\R^d$, contains at least one element of $\varLambda$. Clearly, the
interior $C^{\circ}=C\setminus \{0\}$ of the convex 
cone $C$ contains open balls of arbitrary large radius, hence points
of $\varLambda$. This completes the proof. 
\end{proof}

\begin{rem}
Let $d\in\N$. Note that a subset $\varLambda\subset\R^{d}$ has finite
local complexity (i.e., the difference set $\varLambda-\varLambda$ is
closed and discrete) if and only if for every $r>0$ there are, up to translation, only finitely many
point sets (called {\em patches of diameter} $r$) of the form
$\varLambda\cap B_{r}(x)$, where $x\in\R^d$.
\end{rem}

\begin{theorem}\label{bounded}
Let $d\geq 2$, let $R>0$, and let
$\varLambda\subset\R^{d}$ be a Delone set of finite local
complexity. Then, one has:
\begin{itemize}
\item[(a)]
The set $\mathcal{D}_{<R}(\varLambda)$ is determined by two $X$-rays in $\varLambda$-directions.
\item[(b)]
The set $\mathcal{D}_{<R}(\varLambda)$ is determined by
two projections on orthogonal complements of $1$-dimensional $\varLambda$-subspaces.
\end{itemize}
\end{theorem}
\begin{proof}
Let us first prove part (a). Since $\varLambda$ has finite
local complexity, there are only finitely many possible $\varLambda$-directions having the
property that there may be more than one point of a set $F\in \mathcal{D}_{<R}(\varLambda)$ on a line in this
direction. We denote the finite set of all these $\varLambda$-directions by
$U$. Let $u\in\mathbb{S}^{d-1}$ be an arbitrary
$\varLambda$-direction. For every $F\in
\mathcal{D}_{<R}(\varLambda)$, Lemma~\ref{fgrid} shows that
$$
F\,\,\subset\,\, G_{F}^{\{u\}}\,\,\cap\,\, \varLambda\,.
$$
Choose $u''\in\mathbb{S}^{d-1}\cap(\R u)^{\perp}$ and note
that $(G_{F}^{\{u\}}\cap\varLambda)|(\R u)^{\perp}$ is a finite set with diameter $$D_{F}^{u}:=\operatorname{diam}\big((G_{F}^{\{u\}}\cap\varLambda)|(\R u)^{\perp}\big)<R\,.$$ Since $\varLambda$ has finite
local complexity, the set of diameters $$\left\{D_{F}^{u}\,\left
  |\,F\in
\mathcal{D}_{<R}(\varLambda)\right. \right\}$$ is finite. Set
$$D:=\operatorname{max}\big(\{D_{F}^{u}\, |\,F\in
\mathcal{D}_{<R}(\varLambda)\}\big)<R\,.$$
Note that there is an $\varepsilon_{0}\in\R$ with $0<\varepsilon_{0}<1$ such that every element of
the set $B_{\varepsilon_{0}}(u'')\, \cap\, \mathbb{S}^{d-1}$ is a direction having the property that on each line in
this direction there are no two points of any set $G_{F}^{\{u\}}\cap\varLambda$, $F\in
\mathcal{D}_{<R}(\varLambda)$, on that line with a distance 
$\geq R$.
Since the set of $\varLambda$-directions is dense in
$\mathbb{S}^{d-1}$ by Lemma~\ref{dirdense} (by assumption, $\varLambda$ is
a Delone set and hence relatively dense), and by the finiteness of the
set $U$, this observation shows that one can choose a
$\varLambda$-direction non-parallel to $u$, say $u'$, such that $u'\notin U$,
 and with the property that on each line in
this direction there are no two points of any set $G_{F}^{\{u\}}\cap\varLambda$, $F\in
\mathcal{D}_{<R}(\varLambda)$, on that line with a distance
greater or equal to $R$. We claim that $\mathcal{D}_{<
  R}(\varLambda)$ is determined by the $X$-rays in the
directions $u$ and $u'$. To see this, let $F,F'\in\mathcal{D}_{<
  R}(\varLambda)$ satisfy $X_{u}F=X_{u}F'$. Then, by
Lemma~\ref{fgrid}, one has
$F,F'\subset G_{F}^{\{u\}}\cap\varLambda$. In order to show that the identity
$X_{u'}F=X_{u'}F'$ implies the equality $F=F'$, we shall even
prove that each line in direction $u'$ meets at most one point of
$G_{F}^{\{u\}}\cap\varLambda$. Assume the existence of a line $\ell_{u'}$ in
direction $u'$, and assume the existence of two distinct points $g$ and
$g'$ in $\ell_{u'}\cap (G_{F}^{\{u\}}\cap\varLambda)$. By construction, the
distance of $g$ and
$g'$ is less than $R$. Hence, one has $\{g,g'\}\in \mathcal{D}_{<
  R}(\varLambda)$, and further $u'\in U$, a contradiction. Part (b) follows immediately from an analysis of the proof of part (a).
\end{proof}

Let us note some implications of Theorem~\ref{bounded}.

\begin{coro}\label{meyerb}
Let $d\geq 2$, let $R>0$, and let
$\varLambda\subset\R^{2}$ be a Meyer set. Then, one has:
\begin{itemize}
\item[(a)]
The set $\mathcal{D}_{<R}(\varLambda)$ is determined by two $X$-rays in $\varLambda$-directions.
\item[(b)]
The set $\mathcal{D}_{<R}(\varLambda)$ is determined by
two projections on orthogonal complements of $1$-dimensional $\varLambda$-subspaces.
\end{itemize}
\end{coro}
\begin{proof}
This follows immediately from Theorem~\ref{bounded}, since Meyer sets
$\varLambda$ are, by definition, Delone sets having the property that
$\varLambda-\varLambda$ is uniformly discrete. Clearly, the latter
property implies the finiteness of local complexity.
\end{proof}

\begin{coro}\label{modelb}
Let $d\geq 2$, let $R>0$, and let
$\varLambda\subset\R^{d}$ be an arbitrary model set. Then, one has:
\begin{itemize}
\item[(a)]
The set $\mathcal{D}_{<R}(\varLambda)$ is determined by two $X$-rays in $\varLambda$-directions.
\item[(b)]
The set $\mathcal{D}_{<R}(\varLambda)$ is determined by
two projections on orthogonal complements of $1$-dimensional $\varLambda$-subspaces.
\end{itemize}
\end{coro}
\begin{proof}
This follows immediately from Corollary~\ref{meyerb}, since every
planar model set is a Meyer set; see~\cite{Moody}.
\end{proof}

\begin{coro}\label{coromod}
Let $n\,\in\,\mathbbm{N}\setminus \{1,2\}$, let
$\varLambda_{n}(t,W)\in \mathcal{M}(\OO_{n})$ be a cyclotomic model
set, and let $R>0$. Then, one has:
\begin{itemize}
\item[(a)]
The set $\mathcal{D}_{<R}(\varLambda_{n}(t,W))$ is determined by two $X$-rays in $\OO_{n}$-directions.
\item[(b)]
The set $\mathcal{D}_{<R}(\varLambda_{n}(t,W))$ is determined by
two projections on orthogonal complements of $1$-dimensional $\OO_{n}$-subspaces.
\end{itemize}
\end{coro}
\begin{proof}
This follows from Corollary~\ref{modelb} in conjunction
with Lemma~\ref{dirlem}. 
\end{proof}

\begin{rem}
In Section~\ref{sucdet}, we shall prove a result on 
the {\em successive} determination of the set
$\mathcal{F}(\varLambda_{n}(t,W))$ of all 
finite subsets of a fixed cyclotomic model set
$\varLambda_{n}(t,W)$ by $X$-rays or projections. There, we shall also
 give, independently of the investigations is this section, an alternative proof of Corollary~\ref{coromod}.
\end{rem}

\begin{coro}\label{bounded2}
Let $d\geq 2$, let $t\in\R^{d}$, and let $L\subset \R^{d}$ be a
$($full$)$ lattice. Furthermore, let $R>0$. Then, one has:
\begin{itemize}
\item[(a)]
The set $\mathcal{D}_{<R}(t+L)$
is determined by two $X$-rays in $L$-directions.
\item[(b)]
The set $\mathcal{D}_{<R}(t+L)$
is determined by two projections on orthogonal complements of
$1$-dimensional $L$-subspaces.
\end{itemize}
\end{coro}
\begin{proof}
This follows immediately from Corollary~\ref{modelb}, since translates
of lattices are model sets.
\end{proof}

\begin{rem}
Clearly, Theorem~\ref{bounded} and the subsequent corollaries in this
section are
 best possible with respect to the number of $X$-rays (resp.,
 projections) used.
\end{rem}

\section{Determination of Convex Sets in Cyclotomic Model Sets by $X$-Rays}\label{convdet}

\begin{lem}\label{dim=2}
Let $n\,\in\,\mathbbm{N}\setminus \{1,2\}$ and let
$\varLambda_{n}(t,W)\in \mathcal{M}(\OO_{n})$ be a cyclotomic model set. Let $U\subset \mathbb{S}^{1}$ be a finite set of at least three pairwise non-parallel $\OO_{n}$-directions. Suppose the existence of $F, F'\in \mathcal{C}(\varLambda_{n}(t,W))$ such that $X_{u}F=X_{u}F'$ for all $u\in U$. Then, one has $$F\neq F' \;\Longrightarrow\; \operatorname{dim}(F)=\operatorname{dim}(F')=2\,.$$
\end{lem}
\begin{proof}
No changes needed in comparison with the proof of~\cite[Lemma
5.2]{GG}.
\end{proof}

\begin{rem}
In general, Lemma~\ref{dim=2} is false if one reduces the
number of pairwise non-parallel $\OO_{n}$-directions of $U$ to
two. 
\end{rem}

\begin{ex}
Consider the cyclotomic model set
$\varLambda_{\text{AB}}$ from above. Further, consider the $\OO_{8}$-directions $u:=u_{(2+\sqrt{2})+\zeta_{8}}\in
\mathbb{S}^{1}$ and $u':=u_{-\sqrt{2}+\zeta_{8}}\in
\mathbb{S}^{1}$. Then, the $2$-dimensional convex set
$F:=\{-1-\zeta_{8},-1,0,1,1+\zeta_{8}\}$ in $\varLambda_{\text{AB}}$
and the $1$-dimensional convex set
$F':=\{-1-\sqrt{2},-1,0,1,1+\sqrt{2}\}$ in $\varLambda_{\text{AB}}$
have the same $X$-rays in the directions $u$ and $u'$; see
Figure~\ref{fig:ab} and compare also~\cite[Example 5.3 and Figure 2]{GG}. 
\end{ex}

The proof of the following result is a modified version of that of~\mbox{\cite[Theorem 5.5]{GG}.} For convenience, we present the details.

\begin{theorem}\label{characun}
Let $n\,\in\,\mathbbm{N}\setminus \{1,2\}$ and let
$\varLambda_{n}(t,W)\in \mathcal{M}(\OO_{n})$ be a cyclotomic model set. Further, let $U\subset \mathbb{S}^{1}$ be a set of two or more pairwise non-parallel $\OO_{n}$-directions. The following statements are equivalent:
\begin{itemize}
\item[(i)]
$\mathcal{C}(\varLambda_{n}(t,W))$ is determined by the $X$-rays in the directions of $U$.
\item[(ii)]
There is no $U$-polygon in $\varLambda_{n}(t,W)$.
\end{itemize}
\end{theorem}
\begin{proof}
For (i) $\Rightarrow$ (ii), suppose the existence of a $U$-polygon $P$
in $\varLambda_{n}(t,W)$. Partition the vertices of $P$ into two
disjoint sets $V,V'$, where the elements of these sets alternate round
the boundary $\partial P$ of $P$. Since $P$ is a $U$-polygon, each line
in the plane parallel to some $u\in U$ that contains a point in $V$
also contains a point in $V'$. In particular, one sees that $\operatorname{card}(V)=\operatorname{card}(V')$. Set
$$
C:=(\varLambda_{n}(t,W)\cap P)\setminus (V\cup V')
$$
and further $F_{1}:=C\cup V$ and $F_{2}:=C\cup V'$. Then, $F_{1}$ and $F_{2}$ are different convex sets in $\varLambda_{n}(t,W)$ with the same $X$-rays in the directions of $U$.

For (ii) $\Rightarrow$ (i), suppose that $F_{1}$ and $F_{2}$ are different convex sets in $\varLambda_{n}(t,W)$ with the same $X$-rays in the directions of $U$. Set $$E:=\operatorname{conv}(F_{1})\cap \operatorname{conv}(F_{2})\,.$$ We may assume that $\operatorname{card}(U)\geq 4$, since Lemma~\ref{uleq3} provides a $U$-polygon in $\varLambda_{n}(t,W)$ whenever $\operatorname{card}(U)\leq 3$. By Lemma~\ref{dim=2}, we have $\operatorname{dim}(F_{1})=\operatorname{dim}(F_{2})=2$ and Lemma~\ref{cardinality}(b) shows that $F_{1}$ and $F_{2}$ have the same centroid. It follows that $E^{\circ}\neq \varnothing$.

Since $\operatorname{conv}(F_{1})$ and $\operatorname{conv}(F_{2})$ are convex polygons, one knows that $$(\operatorname{conv}(F_{1})\,\triangle \, \operatorname{conv}(F_{2}))^{\circ}$$ has finitely many components. By the assumption $F_{1}\neq F_{2}$, there is at least one component. Let these components be $C_{j}$, and call $C_{j}$ of type $r\in \{1,2\}$ if $C_{j}\subset (\operatorname{conv}(F_{r})\setminus E)^{\circ}$. Consider the set of type $1$ (resp., type $2$) components together with the equivalence relation generated by the reflexive and symmetric relation $R$ given by adjacency, i.e., $C\,R\, C'\Longleftrightarrow \overline{C}\cap \overline{C'}\neq \varnothing$. Let the set $\mathcal{D}_{1}$ (resp., $\mathcal{D}_{2}$) consist of all unions $\cup\, \mathcal{C}$, where $\mathcal{C}$ is an equivalence class of type $1$ (resp., type $2$) components. Let $\mathcal{D}:= \mathcal{D}_{1}\cup \mathcal{D}_{2}$. Note that the elements of $\mathcal{D}_{1}$ and $\mathcal{D}_{2}$ alternate round the boundary $\partial E$ of $E$. 

Suppose that $D\in \mathcal{D}_{1}$. The set $A:=(\overline{D}\setminus E)\cap \varLambda_{n}(t,W)$ is non-empty, finite and contained in $F_{1}\setminus E$. If $u\in U$ and $z\in A$, then, since $X_{u}F_{1}=X_{u}F_{2}$, there is an element $z' \in \varLambda_{n}(t,W)$ which satisfies
$$
z'\in (F_{2}\setminus E)\cap \ell_{u}^{z}\,.
$$
It follows that $\ell_{u}^{z}$ meets some element of $\mathcal{D}_{2}$. Let us denote this element by $D(u)$.

We first claim that $D(u)$ does not depend on the choice of $z\in A$. To see this, let $\tilde{z}\in A$ be another element of $A$ (i.e., $z\neq \tilde{z}\in A$) such that $\ell_{u}^{\tilde{z}}$ meets $\tilde{D}(u)\in \mathcal{D}_{2}$, where $\tilde{D}(u) \neq D(u)$. The latter inequality implies that $\tilde{D}(u)$ and $D(u)$ are disjoint and, moreover, we see that with respect to the clockwise ordering round $\partial E$ there exists an element $D'$ of $\mathcal{D}_{1}$ between $\tilde{D}(u)$ and $D(u)$. There follows the existence of an element $\hat{z}\in \varLambda_{n}(t,W)$ contained in the open strip bounded by $\ell_{u}^{z}$ and $\ell_{u}^{\tilde{z}}$ such that $\hat{z}\in \overline{C}\setminus E$, where $C$ is one of the type $1$ components contained in $D'$. Since $X_{u}F_{1}=X_{u}F_{2}$, there follows the existence of an element $\hat{z}'\in \varLambda_{n}(t,W) \cap \ell_{u}^{\tilde{z}}$ with $\hat{z}' \in \overline{C'}\setminus E$, where $C'$ is a type $2$ component. It follows that $C'\subset D$, a contradiction. This proves the claim.

The set $A(u):=(\overline{D(u)}\setminus E)\cap \varLambda_{n}(t,W)$ is finite and contained in $F_{2}\setminus E$. Moreover, since $X_{u}A(u)=X_{u}A$, we have $\operatorname{card}(A(u))=\operatorname{card}(A)$ by Lemma~\ref{cardinality}(a). In particular, we see that $A(u)$ is non-empty. 

By symmetry, one gets analogous results for any element $D\in \mathcal{D}_{2}$. 

Choose an arbitrary $D\in \mathcal{D}$ and define the subset
$$
\mathcal{D}':=\left\{((\dots(D(u'_{i_{1}}))\dots )(u'_{i_{k-1}}))(u'_{i_{k}})\, | \, k\in \mathbbm{N}, u'_{i_{j}}\in U \mbox{ for all }j\in\{1,\dots,k\}\right\}\,.
$$
of $\mathcal{D}$, obtained from $D$ by applying the above process through any finite sequence of directions from $U$. Let $\mathcal{D}'=\{D_{j}\,|\, j \in \{1,\dots,m\}\}$ and let $A_{j}:=(\overline{D_{j}}\setminus E)\cap \varLambda_{n}(t,W)$ be the non-empty set of elements of $\varLambda_{n}(t,W)$ corresponding to $D_{j}$, $j\in\{1,\dots,m\}$.

Let $c_{j}$ be the centroid of $A_{j}$, $j\in\{1,\dots,m\}$, and let $t_{j}$ be the line through the common endpoints of the two arcs, one in $\partial(\operatorname{conv}(F_{1}))$, the other in $\partial(\operatorname{conv}(F_{2}))$, which bound $D_{j}$. Then, $t_{j}$ separates $A_{j}$, and hence $c_{j}$, from the convex hull of the remaining centroids $c_{k}$, with $k\in \{1,\dots,m\}\setminus \{j\}$. It follows that the points $c_{j}$, $j\in\{1,\dots,m\}$, are the vertices of a convex polygon $P$. If $u\in U$ and $j\in\{1,\dots,m\}$, suppose that $A_{k}$ is the set arising from $u$ and $A_{j}$ by the process described above, i.e., $A_{k}=A_{j}(u)$. Then, by Lemma~\ref{cardinality}(b), the line $\ell_{u}^{c_{j}}$ also contains $c_{k}$. The points $c_{j}$ therefore pair off in this fashion, so $m$ is even, and since $\operatorname{card}(U)\geq 2$, we  have $m\geq 4$, and $P$ is non-degenerate. Hence, $P$ is a $U$-polygon.

Let $\operatorname{card}(A_{1})=\dots=\operatorname{card}(A_{m})=:s\in \mathbbm{N}$. Then, each vertex of $P$ belongs to $t+\frac{1}{s}\OO_{n}\subset t+\mathbbm{K}_{n}$. Hence, $P$ is a $U$-polygon in $t+\mathbbm{K}_{n}$. By Lemma~\ref{dilate}, there is a homothety $h\!:\, \C \longrightarrow \C$ such that $P':= h[P]$ is a polygon in $\varLambda_{n}(t,W)$. Since $P'$ is a $U$-polygon (see Lemma~\ref{homotu}(a)), $P'$ is a $U$-polygon in $\varLambda_{n}(t,W)$. 
\end{proof}

\begin{rem}
In the proof of Theorem~\ref{characun}, it is necessary to employ finite unions of components; compare~\cite[Remark 5.6]{GG}.   
\end{rem}

We are now able to prove the following result.

\begin{theorem}\label{main1}
Let $n\,\in\,\mathbbm{N}\setminus \{1,2\}$. Then, one has:
\begin{itemize}
\item[(a)]
There is a set $U\in \mathcal{U}_{4,\mathbbm{Q}}^{n}$ with
$\operatorname{card}(U)=4$ such that, for all cyclotomic model sets
$\varLambda_{n}(t,W)\in \mathcal{M}(\OO_{n})$, the set
$\mathcal{C}(\varLambda_{n}(t,W))$ is determined by the $X$-rays in
the directions of $U$. Furthermore, any set $U\in \mathcal{U}_{4,\mathbbm{Q}}^{n}$ with
$\operatorname{card}(U)=4$ having the property that the cross ratio of slopes of the
directions of $U$, arranged in order of increasing angle with the
positive real axis, is not an element of the set
$N_{1}$ $($as defined in
Theorem~$\ref{intersectq}$$)$ is suitable for this purpose.
\item[(b)]
Let $U\subset \mathbb{S}^{1}$ be any set of three or less pairwise
non-parallel $\OO_{n}$-directions. Then, for all cyclotomic model sets $\varLambda_{n}(t,W)\in \mathcal{M}(\OO_{n})$, the set $\mathcal{C}(\varLambda_{n}(t,W))$ is not determined by the $X$-rays in the directions of $U$.  
\item[(c)]
Let $U\in \mathcal{U}_{4,\mathbbm{Q}}^{n}$ with
$\operatorname{card}(U)=7$. Then, for all cyclotomic model sets $\varLambda_{n}(t,W)\in \mathcal{M}(\OO_{n})$, the set $\mathcal{C}(\varLambda_{n}(t,W))$ is determined by the $X$-rays in the directions of $U$. 
\item[(d)]
 There is a set $U\in \mathcal{U}_{4,\mathbbm{Q}}^{n}$ with
 $\operatorname{card}(U)=6$ such that, for all cyclotomic model sets $\varLambda_{n}(t,W)\in \mathcal{M}(\OO_{n})$, the set $\mathcal{C}(\varLambda_{n}(t,W))$ is not determined by the $X$-rays in the directions of $U$. 
\end{itemize}
\end{theorem}
\begin{proof}
Let us start with (a). In view of Theorem~\ref{finitesetcr} and
Theorem~\ref{characun}, the additional statement is immediate. So, it
remains to show existence. For example, the following sets of
$\OO_{n}$-directions have the property that the cross ratio of slopes of their
directions, arranged in order of increasing angle with the
positive real axis, is not an element of the set
$N_{1}$: $$U_{n}:=\{u_{1},u_{1+\zetan},u_{1+2\zetan},u_{1+5\zetan}\}\,,$$
$$U'_{n}:=\{u_{1},u_{2+\zetan},u_{\zetan},u_{-1+2\zetan}\}$$ and
$$U''_{n}:=\{u_{2+\zetan},u_{3+2\zetan},u_{1+\zetan},u_{2+3\zetan}\}\,;$$ compare~\cite[Remark
5.8]{GG} and Remark~\ref{r4}. For these sets, the cross ratio
of slopes of the
directions, arranged in order of increasing angle with the
 positive real axis, is, in order of their appearance above, $8/5$,
 $5/4$ and $5/4$ again. 

Assertions (b) and (c) are immediate consequences of Theorem~\ref{characun} in conjunction with Lemma~\ref{uleq3} and Theorem~\ref{finitesetcr}, respectively. 

For (d), note that there is a
$U$-polygon $P'$ in $\varLambda_{n}(t,W)$, where $U$
consists of the six pairwise non-parallel $\OO_{n}$-directions
 $u_{1}$, $u_{2+\zetan}$, $u_{1+\zetan}$, $u_{1+2\zetan}$, $u_{\zetan}$ and
$u_{-1+\zetan}$, respectively. To see this, let $P$ be the
non-degenerate convex dodecagon with vertices at $3+\zetan$,
$3+2\zetan$, $2+3\zetan$, $1+3\zetan$, $-1+2\zetan$, $-2+\zetan$, and
the reflections of these points in the origin $0$;
compare~\cite[Example 4.3 and Figure 1]{GG} and see
Remark~\ref{r4}. Then, one easily verifies that $P$ is a $U$-polygon
in $\OO_{n}$. Invoking Lemma~\ref{dilate} again, there is a homothety
$h\!:\, \C \longrightarrow \C$ such that
$P':= h[P]$ is a polygon in
$\varLambda_{n}(t,W)$. Since $P'$ is a $U$-polygon (see
Lemma~\ref{homotu}(a)), $P'$ is a $U$-polygon in
$\varLambda_{n}(t,W)$. The assertion now follows immediately from Theorem~\ref{characun}. 
\end{proof}

\begin{rem}
By assertions (b) and (d) of Theorem~\ref{main1}, assertions (a) and (c) of Theorem~\ref{main1} are optimal with respect to the number of directions used.
\end{rem}

As shown in
Lemma~\ref{u4nq346} for $n\in\{3,4,6\}$ (corresponding to periodic
cyclotomic model sets), the set
$\mathcal{U}_{4,\mathbbm{Q}}^{n}$ coincides with the natural set of all sets $U\subset
\mathbb{S}^{1}$ of four or more pairwise non-parallel
$\OO_{n}$-directions. Though, for $n\in\N\setminus\{3,4,6\}$, the
restriction of the set of all sets of four or more pairwise
non-parallel $\OO_{n}$-directions to the set
$\mathcal{U}_{4,\mathbbm{Q}}^{n}$ in Theorem~\ref{main1} is rather
artificial. Now, we use our results from Section~\ref{upol} in order to remove
this restriction and prove corresponding more results in larger generality. For the case $n\in\{5,8,10,12\}$ (corresponding to (aperiodic) cyclotomic model 
sets $\varLambda_{n}(t,W)\in \mathcal{M}(\OO_{n})$ with co-dimension
 two), one has the following generalization of Theorem~\ref{main1}(a).

\begin{theorem}\label{main2}
Let $n\in\{5,8,10,12\}$ and let $U\subset \mathbb{S}^1$ be any set of four pairwise non-parallel
$\OO_{n}$-directions having the property
\begin{itemize}
\item[(C2)]
The cross ratio of slopes of the
directions of $U$, arranged in order of increasing angle with the
positive real axis, does not map under the norm
$N_{\mathbbm{k}_{n}/\Q}$ to the set $N_{2}$ as defined in
Theorem~$\ref{intersectk8}$.
\end{itemize}
 Then, for all $($aperiodic$)$ cyclotomic model sets
 $\varLambda_{n}(t,W)\in \mathcal{M}(\OO_{n})$, the set\linebreak $\mathcal{C}(\varLambda_{n}(t,W))$ is determined by the $X$-rays in the directions of $U$.
\end{theorem}
\begin{proof}
This follows immediately from Theorem~\ref{characun} and
Theorem~\ref{finitesetncr}. 
\end{proof}

\begin{ex}\label{ex1}
As an example with $n=8$, the following set of
$\OO_{8}$-directions has Property~(C2):
$$
U_8:=\left\{u_{1+\zeta_{8}},u_{(-1+\sqrt{2})+\sqrt{2}\zeta_{8}},u_{(-1-\sqrt{2})+\zeta_{8}},u_{-2+(-1+\sqrt{2})\zeta_8}\right\}\,.
$$
Here, the cross ratio of slopes of the
elements of $U_8$, arranged in order of increasing angle with the
positive real axis, equals $12/7-3/7\sqrt{2}$, hence $$N_{\mathbbm{k}_{8}/\Q}\left(\frac{12}{7}-\frac{3}{7}\sqrt{2}\right)=\frac{18}{7}\notin N_{2}\,.$$
 Further, for $n\in\{5,10\}$, the following set of
$\OO_{5}$-directions has Property~(C2):
$$
U_{10}:=U_{5}:=\left\{u_{(1+\tau)+\zeta_{5}},u_{(\tau-1)+\zeta_{5}},u_{-\tau+\zeta_{5}},u_{2\tau-\zeta_{5}}\right\}\,.
$$
Here, the cross ratio of slopes of the
elements of $U_5$, arranged in order of increasing angle with the
positive real axis, equals $4/5+1/5\sqrt{5}$, hence $$N_{\mathbbm{k}_{5}/\Q}\left(\frac{4}{5}+\frac{1}{5}\sqrt{5}\right)=\frac{11}{25}\notin N_{2}\,.$$
 Finally, for $n=12$, the following set of
$\OO_{12}$-directions has Property~(C2):
$$
U_{12}:=\left\{u_{1},u_{2+\zeta_{12}},u_{\zeta_{12}},u_{\sqrt{3}-\zeta_{12}}\right\}\,.
$$
Here, the cross ratio of slopes of the
elements of $U_{12}$, arranged in order of increasing angle with the
positive real axis, equals $2+1/2\sqrt{3}$, hence $$N_{\mathbbm{k}_{12}/\Q}\left(2+\frac{1}{2}\sqrt{3}\right)=\frac{13}{4}\notin N_{2}\,.$$
 \end{ex}

\begin{rem}
Using the fact that the quadratic fields $\kk_{8},\kk_{5}$ and
$\kk_{12}$ have class number $1$, Pleasants~\cite{PABP2} showed that the sets of
$\OO_{n}$-directions $U_8,U_{5}$ and
$U_{12}$ in Example~\ref{ex1} are well suited in order to yield dense lines in the corresponding (aperiodic) cyclotomic model
sets. Note that precise statements on densities of lines in aperiodic cyclotomic
model sets depend on the shape of the window. It follows that,
 in the practice of quantitative HRTEM,
 the resolution coming from the above directions
 is likely to be rather high, which makes Theorem~\ref{main2} look promising.
\end{rem}

For the general case $n\in\N\setminus\{1,2\}$, one has the following
generalizations of Theorem~\ref{main1}(a) and Theorem~\ref{main2},
dealing with arbitrary (but fixed) cyclotomic model sets and arbitrary sets $U$ of four pairwise
non-parallel $\OO_{n}$-directions.

\begin{theorem}\label{main3}
For all $e\in\operatorname{Im}(\phi/2)$, there is a finite set $N_{e}\subset \Q$
such that, for all $n\in(\phi/2)^{-1}[\{e\}]$, and all sets $U\subset \mathbb{S}^1$ of four pairwise
non-parallel $\OO_{n}$-directions, one has the following:

If $\,U$ has the property that the cross ratio of slopes of the
directions of $\,U$, arranged in order of increasing angle with the
positive real axis, does not map under the norm
$N_{\mathbbm{k}_{n}/\Q}$ to $N_{e}$, then, for all cyclotomic model sets $\varLambda_{n}(t,W)\in \mathcal{M}(\OO_{n})$, the set $\mathcal{C}(\varLambda_{n}(t,W))$ is determined by the $X$-rays in the directions of $U$.
\end{theorem}
\begin{proof}
This follows immediately from Theorem~\ref{characun} and
Theorem~\ref{finitesetncrgeneral}.
\end{proof}

\begin{theorem}\label{main32}
For all $n\in\N\setminus\{1,2\}$ and all sets $U\subset \mathbb{S}^1$ of four pairwise
non-parallel $\OO_{n}$-directions, one has the following:

If $U$ has the property that the cross ratio of slopes of the
directions of $U$, arranged in order of increasing angle with the
positive real axis, does not map under the norm
$N_{\mathbbm{k}_{n}/\Q}$ to the set $$\pm\left(\{1\}\cup \big[\Q_{>1}^{\mathbbm{P}_{\leq 2}}\big]^{\pm 1}\right)\,,$$ 
then, for all cyclotomic model sets $\varLambda_{n}(t,W)\in
\mathcal{M}(\OO_{n})$, the set $\mathcal{C}(\varLambda_{n}(t,W))$ is
determined by the $X$-rays in the directions of $U$.
\end{theorem}
\begin{proof}
This follows immediately from Theorem~\ref{characun} and
Theorem~\ref{finitesetncrgeneral2}.
\end{proof}

For the case of periodic cyclotomic model sets, we are able to prove the first uniqueness result, which is in full accordance with
the setting of Section~\ref{secset}.

\begin{theorem}\label{main22}
Let $n\in\{3,4,6\}$ and let $U\subset \mathbb{S}^1$ be any set of four pairwise non-parallel
$\OO_{n}$-directions having the following properties.
\begin{itemize}
\item[(C1)]
The cross ratio of slopes of the
directions of $U$, arranged in order of increasing angle with the
positive real axis, does not map under the norm
$N_{\mathbbm{k}_{n}/\Q}$ to the set $N_{1}$ as defined in
Theorem~$\ref{intersectq}$.
\item[(E)]
$U$ contains two $\OO_{n}$-directions of the form $u_{o},u_{o'}$, where
$o,o'\in\OO_{n}\setminus\{0\}$ satisfy one of the equivalent conditions {\rm
  (i)}-{\rm (iii)} 
of Proposition~$\ref{proponeequiv}$.
\end{itemize}
 Then, the set 
$\cup_{t\in\R^2}\,\,\mathcal{C}(t+\OO_{n})$
 is determined by the $X$-rays in the directions of $U$.
\end{theorem}
\begin{proof}
Let $n\in\{3,4,6\}$ and let $U\subset \mathbb{S}^1$ be a set of four pairwise non-parallel
$\OO_{n}$-directions having the Properties (C1) and (E). Let $F,F'\in \cup_{t\in\R^2}\mathcal{C}(t+\OO_{n})$, say $F\in\mathcal{C}(t+\OO_{n})$ and $F'\in\mathcal{C}(t'+\OO_{n})$, where $t,t'\in\R^2$, and suppose that $F$ and $F'$ have the same $X$-rays in the directions of $U$. Then, by Lemma~\ref{fgrid} and Theorem~\ref{oneequiv} in conjunction with Property~(E), one obtains
\begin{equation}\label{fton}
F,F'\,\,\subset \,\,G_{F}^{U}\,\,\subset\,\,t+\OO_{n}\,.
\end{equation}
If $F=\varnothing$, then, by Lemma~\ref{cardinality}(a), one also gets
$F'=\varnothing$. It follows that one may assume, without loss of generality, that $F$ and $F'$ are
non-empty. Then, since $F'\subset t'+\OO_{n}$, it follows from
Equation~(\ref{fton}) that $t+\OO_{n}$ meets $t'+\OO_{n}$, the latter
being equivalent to the identity $t+\OO_{n}=t'+\OO_{n}$. Now, the assertion follows
immediately from Property~(C1) and Theorem~\ref{main1}(a) in
conjunction with Lemma~\ref{u4nq346}. 
\end{proof}

\begin{ex}\label{ex2}
It was shown in the proof of Theorem~\ref{main1}(a) that, for
$n\in\{3,4\}$, the sets of
$\OO_{n}$-directions $U_n,U'_n,U''_n$ (as defined in the proof of
Theorem~\ref{main1}(a)) have Property~(C1). One can easily see that these sets of
$\OO_{n}$-directions additionally have Property~(E).
\end{ex}

\begin{rem}
For $n\in\{3,4\}$, the sets of
$\OO_{n}$-directions parallel to the elements of the sets
$U_n,U'_n,U''_n$ obviously yield dense lines in the corresponding (periodic) cyclotomic model
sets. It follows that,
 in the practice of quantitative HRTEM,
 the resolution coming from these directions
 is rather high, so their might be a real application of Theorem~\ref{main22}.
\end{rem}

\begin{rem}
In an approximative sense, which will be made precise in the following, one is
also able to deal with the `non-anchored' case for regular aperiodic
cyclotomic model sets, i.e., the
determination of sets of the form $\cup_{\varLambda\in
   W^{\star_{n}}_{\mathcal{M}_{g}(\mathcal{O}_{n})}}\mathcal{C}(\varLambda)$, where $n\in\N\setminus\{1,2,3,4,6\}$, by the $X$-rays in four prescribed $\OO_{n}$-directions. Note first that, for $n\in\N\setminus\{1,2,3,4,6\}$,
Theorem~\ref{main3} shows that there is a finite set $N_{\phi(n)/2}\subset \Q$
such that, for all sets $U\subset \mathbb{S}^1$ of four pairwise
non-parallel $\OO_{n}$-directions, one has: 

If $U$ has the property that the cross ratio of slopes of the
directions of $U$, arranged in order of increasing angle with the
positive real axis, does not map under the norm
$N_{\mathbbm{k}_{n}/\Q}$ to $N_{\phi(n)/2}$, then, for all
cyclotomic model sets $\varLambda_{n}(t,W)\in \mathcal{M}(\OO_{n})$,
the set $\mathcal{C}(\varLambda_{n}(t,W))$ is determined by the
$X$-rays in the directions of $U$. 

Now let $U\subset \mathbb{S}^1$ be any set of four pairwise non-parallel
$\OO_{n}$-directions having the following properties.
\begin{itemize}
\item[(C)]
The cross ratio of slopes of the
directions of $U$, arranged in order of increasing angle with the
positive real axis, does not map under the norm
$N_{\mathbbm{k}_{n}/\Q}$ to the set $N_{\phi(n)/2}$ from above.
\item[(E)]
$U$ contains two $\OO_{n}$-directions of the form $u_{o},u_{o'}$, where
$o,o'\in\OO_{n}\setminus\{0\}$ satisfy one of the equivalent conditions {\rm
  (i)}-{\rm (iii)} 
of Proposition~$\ref{proponeequiv}$.
\end{itemize}
 We claim that, in a sense, for all fixed windows  $W\subset(\R^2)^{\phi(n)/2-1}$ with boundary $\partial W$ having Lebesgue
measure $0$ in $(\R^2)^{\phi(n)/2-1}$ and all fixed star maps $$.^{\star_{n}}\! : \, \mathcal{O}_{n}\longrightarrow
(\R^2)^{\frac{\phi(n)}{2}-1}$$ $($as described in
 Definition~$\ref{cyclodef}$$)$, the set 
$\cup_{\varLambda\in
   W^{\star_{n}}_{\mathcal{M}_{g}(\mathcal{O}_{n})}}\mathcal{C}(\varLambda)$
 is determined by the $X$-rays in the directions of $U$. To see this, let $F,F'\in \cup_{\varLambda\in
   W^{\star_{n}}_{\mathcal{M}_{g}(\mathcal{O}_{n})}}\mathcal{C}(\varLambda)$, say $F\in\mathcal{C}(\varLambda_{n}^{\star_{n}}(t,\tau+W))$ and $F'\in\mathcal{C}(\varLambda_{n}^{\star_{n}}(t',\tau'+W))$, where $t,t'\in\R^2$ and $\tau,\tau'\in(\R^2)^{\frac{\phi(n)}{2}-1}$, and suppose that $F$ and $F'$ have the same $X$-rays in the directions of $U$. Then, by Lemma~\ref{fgrid} and Theorem~\ref{oneequiv} in conjunction with Property~(E), one obtains
\begin{equation}\label{fton225}
F,F'\,\,\subset \,\,G_{F}^{U}\,\,\subset\,\,t+\OO_{n}\,.
\end{equation}
If $F=\varnothing$, then, by Lemma~\ref{cardinality}(a), one also gets
$F'=\varnothing$. One may thus assume, without loss of generality, that $F$ and $F'$ are
non-empty. Then, since $F'\subset t'+\OO_{n}$, 
Equation~(\ref{fton}) implies that $t+\OO_{n}$ meets $t'+\OO_{n}$, the latter
being equivalent to the identity $t+\OO_{n}=t'+\OO_{n}$. Moreover, 
 the identity $t+\OO_{n}=t'+\OO_{n}$ is equivalent to the relation $t'-t\in\OO_{n}$. Hence, one
has $$F-t\in\mathcal{C}\Big(\varLambda_{n}^{\star_{n}}\big(0,\tau+W\big)\Big)$$
and, since the equality
$\varLambda_{n}^{\star_{n}}(t'-t,\tau'+W)=\varLambda_{n}^{\star_{n}}(0,(\tau'+(t'-t)^{\star_{n}})+W)$
holds,
$$F'-t\in\mathcal{C}\Big(\varLambda_{n}^{\star_{n}}\big(t'-t,\tau'+W\big)\Big)=\mathcal{C}\Big(\varLambda_{n}^{\star_{n}}\big(0,(\tau'+(t'-t)^{\star_{n}})+W\big)\Big)\,.$$
Clearly, $F-t$ and $F'-t$ again have the same $X$-rays in the
directions of $U$. Hence, by Lemma~\ref{cardinality}(b), $F-t$ and $F'-t$
have the same centroid. Since the star map $.^{\star_{n}}$ is $\Q$-linear, it follows that the finite subsets
$[F-t]^{\star_{n}}$ and $[F'-t]^{\star_{n}}$ of
$(\R^2)^{\phi(n)/2-1}$ also have the same
centroid. Now, if one has $F-t=B_{R}(a)\cap
\varLambda_{n}^{\star_{n}}(t,\tau+W)$ and
$F'-t=B_{R'}(a')\cap
\varLambda_{n}^{\star_{n}}(0,(\tau'+(t'-t)^{\star_{n}})+W)$ for
suitable $a,a'\in\R^2$ and large $R,R'>0$ (which is rather natural in
practice), then Theorem~\ref{weyl} allows us to write
\begin{eqnarray*}
\frac{1}{\operatorname{vol}(W)}\int_{\tau +W}y\,{\rm
  d}\lambda(y)&\approx&
\frac{1}{\operatorname{card}\left(F-t\right)}\sum_{x\in
  F-t}x^{\star_{n}}\\ &=&
\frac{1}{\operatorname{card}\left(F'-t\right)}\sum_{x\in
  F'-t}x^{\star_{n}}\\
&\approx&\frac{1}{\operatorname{vol}(W)}\int_{(\tau'+(t'-t)^{\star_{n}})
  +W}y\,{\rm d}\lambda(y)\,.
\end{eqnarray*}
Consequently,  
$$
\tau+\int_{W}y\,{\rm
  d}\lambda(y)\approx(\tau'+(t'-t)^{\star_{n}})+\int_{W}y\,{\rm d}\lambda(y)\,,
$$
and hence $\tau\approx\tau'+(t'-t)^{\star_{n}}$. The latter means
that, approximately, both
$F-t$ and $F'-t$ are elements of the set
$\mathcal{C}(\varLambda_{n}^{\star_{n}}(0,\tau+W))$. Now, it follows in a loose
sense 
 from Property~(C) and Theorem~\ref{main3} that $F-t\approx F'-t$,
 and hence $F\approx F'$.
\end{rem}

\begin{ex}\label{ex4}
For $n\in\{5,8,10,12\}$, the set of
$\OO_{n}$-direction $U_n$ (as
defined in Example~\ref{ex1}) has
Property~(C) with $N_2$ as defined in 
Theorem~$\ref{intersectk8}$. One can easily see that these $\OO_{n}$-directions additionally have Property~(E).
\end{ex}

\section{Successive Determination of Finite Subsets by $X$-Rays or
  Projections -- A Class of Delone Sets and Cyclotomic Model Sets}\label{sucdet}

Here, we shall investigate the successive determination of finite
subsets (not necessarily convex sets) of both cyclotomic model
sets and arbitrary Delone sets living on $\OO_{n}$, where $n\in\N\setminus\{1,2\}$. Recall that the interactive technique of
successive determination allows us to use the information from previous
$X$-rays in deciding on the direction for the next $X$-ray.

\begin{defi}
For $n\,\in\,\mathbbm{N}\setminus \{1,2\}$ and
$j\in\{1,\dots,\phi(n)/2\}$, set
$b_{j}^{_{(n)}}:=(\zeta_{n}+\bar{\zeta}_{n})^{j-1}$
and $b_{\phi(n)/2+j}^{_{(n)}}:=(\zeta_{n}+\bar{\zeta}_{n})^{j-1}\zeta_{n}$.
Further, set $$B_{1}^{_{(n)}}:=\{b_1^{_{(n)}},\dots,b_{\phi(n)/2}^{_{(n)}}\}\,,$$
$$B_{2}^{_{(n)}}:=\{b_{\phi(n)/2+1}^{_{(n)}},\dots,b_{\phi(n)}^{_{(n)}}\}\,,$$ and, finally,
 $B_{\vphantom{1}}^{_{(n)}}:=B_{1}^{_{(n)}}\,\,\dot\cup \, \,B_{2}^{_{(n)}}$.
\end{defi}

\begin{lem}\label{sbasis}
Let $n\,\in\,\mathbbm{N}\setminus \{1,2\}$. Then,
$B_{\vphantom{1}}^{_{(n)}}$ is both a $\mathbbm{Q}$-basis of $\mathbbm{K}_{n}$ and a
$\mathbbm{Z}$-basis of $\OO_{n}$.
\end{lem}
\begin{proof}
The assertion follows immediately from Lemma~\ref{Oo}, Corollary~\ref{cr5} and the second part of
Remark~\ref{r1}.
\end{proof}

\begin{defi}
For $n\,\in\,\mathbbm{N}\setminus \{1,2\}$, let $\,\,.\,\widetilde{\hphantom{a}}^{_{n}}$ be an arbitrary, but fixed Minkowski embedding
 of $\mathbbm{K}_{n}$, i.e., a map 
$
\,\,.\,\widetilde{\hphantom{a}}^{_{n}}\,:\, \mathbbm{K}_{n}\longrightarrow
(\R^2)^{\phi(n)/2}\,,
$
given by
$$
z\longmapsto \big(\sigma_{1}(z),\sigma_{2}(z),\dots,\sigma_{\phi(n)/2}(z)\big)\,,
$$
where the set $\{\sigma_{1},\dots,\sigma_{\phi(n)/2}\}$ arises from
$G(\mathbbm{K}_{n}/ \mathbbm{Q})$ by choosing exactly one automorphism
from each pair of complex conjugate automorphisms. 
\end{defi}

\begin{rem}\label{latticerem}
For the notion of Minkowski embeddings, compare also
  Remark~\ref{mink} and references given there. Note that the map $\,\,.\,\widetilde{\hphantom{a}}^{_{n}}\,$ is $\Q$-linear and
 injective. Further, by Lemma~\ref{sbasis}, the set $[\OO_{n}]\widetilde{\hphantom{a}}^{_{n}}$ is a lattice in
$(\R^2)^{\phi(n)/2}$ with basis $[B_{\vphantom{1}}^{_{(n)}}]\widetilde{\hphantom{a}}^{_{n}}$;
see~\cite[Ch. 2, Sec. 3]{Bo}. Further, again by Lemma ~\ref{sbasis}
and the $\Q$-linearity of $\,\,.\,\widetilde{\hphantom{a}}^{_{n}}\,$, one has $[\mathbbm{K}_{n}]\widetilde{\hphantom{a}}^{_{n}}=\langle[B_{\vphantom{1}}^{_{(n)}}]\widetilde{\hphantom{a}}^{_{n}}\rangle_{\Q}$.
\end{rem}

\begin{defi}
For $n\,\in\,\mathbbm{N}\setminus \{1,2\}$, we let $\,\,.\,\widehat{\hphantom{a}}\,^{_{n}}$ be the
co-restriction of $\,\,.\,\widetilde{\hphantom{a}}^{_{n}}$ to its image
$[\mathbbm{K}_{n}]\widetilde{\hphantom{a}}^{_{n}}$.
\end{defi}

\begin{lem}\label{hatiso}
Let $n\,\in\,\mathbbm{N}\setminus \{1,2\}$. Then,
$\,\,.\,\widehat{\hphantom{a}}\,^{_{n}}$ is a $\Q$-linear isomorphism.
\end{lem}
\begin{proof}
This follows from the fact that the map $\,\,.\,\widetilde{\hphantom{a}}^{_{n}}\,$ is $\Q$-linear and
 injective. 
\end{proof}

\begin{defi}
For $n\,\in\,\N$ and $\kappa\in
\mathbbm{K}_{n}$, we let $m_{\kappa}^{_{(n)}}\!:\,\mathbbm{K}_{n} \longrightarrow
\mathbbm{K}_{n}$ be the map that is given by multiplication by $\kappa$, i.e., $z\longmapsto \kappa z$.
\end{defi}

\begin{rem}\label{maut}
Let $n\,\in\,\N$. Note that $m_{\kappa}^{_{(n)}}$ is a
$\mathbbm{Q}$-linear endomorphism of
the $\phi(n)$-dimensional vector space $\mathbbm{K}_{n}$ over
$\mathbbm{Q}$; cf. Proposition~\ref{gau}. Further, $m_{\kappa}^{_{(n)}}$ is a $\mathbbm{Q}$-linear automorphism if and only if $\kappa\neq
0$. In particular, if $n\,\in\,\mathbbm{N}\setminus
\{1,2,3,4,6\}$, then  $m_{(\zeta_{n}+\bar{\zeta}_{n})}^{_{(n)}}$ is a
$\Q$-linear automorphism, since then $\zeta_{n}+\bar{\zeta}_{n}\neq
0$. In fact, this
 restriction means that $\phi(n)/2\geq 2$ from which follows
 that $\zeta_{n}+\bar{\zeta}_{n}\notin \Q$; cf. Corollary~\ref{cr5}.

The
 $\Q$-linear endomorphism $m_{\kappa}^{_{(n)}}$ (resp., automorphism, if $\kappa\neq 0$) corresponds via the $\Q$-linear isomorphism $\widehat{\hphantom{a}}\,^{_{n}}$ to a
 $\mathbbm{Q}$-linear endomorphism (resp., automorphism), say
 $(m_{\kappa}^{_{(n)}})\widehat{\hphantom{a}}\,^{_{n}}$, of $[\mathbbm{K}_{n}]\widetilde{\hphantom{a}}^{_{n}}$, i.e.,
 $$(m_{\kappa}^{_{(n)}})\widehat{\hphantom{a}}\,^{_{n}}=\,\,.\,\widehat{\hphantom{a}}\,^{_{n}}\;\circ
 m_{\kappa}^{_{(n)}}\circ (\,\,.\,\widehat{\hphantom{a}}\,^{_{n}})^{-1}\,.$$
Note further that $(m_{\kappa}^{_{(n)}})\widehat{\hphantom{a}}\,^{_{n}}$ extends uniquely
to an $\R$-linear endomorphism (resp., automorphism), say $(m_{\kappa}^{_{(n)}})\widetilde{\hphantom{a}}^{_{n}}$, of $(\R^2)^{\phi(n)/2}$; cf. Remark~\ref{latticerem}.
\end{rem}

The following result characterizes, for $n\,\in\,\mathbbm{N}\setminus
\{1,2\}$, the intersections of $1$-dimensional $\OO_{n}$-subspaces in the Euclidean plane with the $n$th
cyclotomic field $\mathbbm{K}_{n}$ in terms of existence of certain $\Q$-bases.

\begin{lem}\label{lines}
Let $n\,\in\,\mathbbm{N}\setminus \{1,2\}$ and let
$o\in\OO_{n}\setminus \{0\}$. Then, one has
 $$\mathbbm{K}_{n}\cap (\R o)=\mathbbm{k}_{n}o\,.$$ Moreover,
 $\mathbbm{K}_{n}\cap (\R o)$ is a $(\phi(n)/2)$-dimensional $\Q$-linear
subspace of $\mathbbm{K}_{n}$ with $\Q$-basis $B_{1}^{_{(n)}}o$.
\end{lem}
\begin{proof}
If $z\in \mathbbm{K}_{n}\cap (\R o)$, one has $z=\lambda o$
with $\lambda \in \mathbbm{k}_{n}$ by Proposition~\ref{p1}. The assertion follows immediately from the last observation,
Corollary~\ref{cr5}, and our assumption $o\neq 0$.
\end{proof}

\begin{prop}
Let $n\,\in\,\mathbbm{N}\setminus \{1,2\}$ and let
$o\in\OO_{n}\setminus \{0\}$. Then, one has:
\begin{itemize}
\item[(a)]$\langle[\OO_{n}\cap
    (\R o)]\widetilde{\hphantom{a}}^{_{n}} \rangle_{\R}$
is a $(\phi(n)/2)$-dimensional $[\OO_{n}]\widetilde{\hphantom{a}}^{_{n}}$-subspace of
$(\R^2)^{\phi(n)/2}$ with $\R$-basis $$[B_{1}^{_{(n)}}o]\widetilde{\hphantom{a}}^{_{n}}=\left\{o\widetilde{\hphantom{a}}^{_{n}},(m_{(\zeta_{n}+\bar{\zeta}_{n})}^{_{(n)}})\widehat{\hphantom{a}}\,^{_{n}}(o\widetilde{\hphantom{a}}^{_{n}}),\dots,\big((m_{(\zeta_{n}+\bar{\zeta}_{n})}^{_{(n)}})\widehat{\hphantom{a}}\,^{_{n}}\big)^{\frac{\phi(n)}{2}-1}(o\widetilde{\hphantom{a}}^{_{n}})\right\}\,.$$ 
\smallskip
\vspace{-2mm}
\item[(b)]
$\langle[\OO_{n}\cap
    (\R o)]\widetilde{\hphantom{a}}^{_{n}} \rangle_{\R}\cap
    [\mathbbm{K}_{n}]\widetilde{\hphantom{a}}^{_{n}}=\left[\mathbbm{k}_{n}o\right]\widetilde{\hphantom{a}}^{_{n}}$ and
    $[B_{1}^{_{(n)}}o]\widetilde{\hphantom{a}}^{_{n}}$ is a $\Q$-basis of $$\left\langle\big[\OO_{n}\cap
    (\R o)\big]\widetilde{\hphantom{a}}^{_{n}}\right \rangle_{\R}\cap
    \big[\mathbbm{K}_{n}\big]\widetilde{\hphantom{a}}^{_{n}}\,.$$
\end{itemize}
\end{prop}
\begin{proof}
First, observe that 
\begin{eqnarray*}
[B_{1}^{_{(n)}}o]\widetilde{\hphantom{a}}^{_{n}}&=&\left\{(b_{1}^{_{(n)}}o)\widetilde{\hphantom{a}}^{_{n}},(b_{2}^{_{(n)}}o)\widetilde{\hphantom{a}}^{_{n}},\dots,(b_{\frac{\phi(n)}{2}}^{_{(n)}}o)\widetilde{\hphantom{a}}^{_{n}}\right\}\\
&=&\left\{o\widetilde{\hphantom{a}}^{_{n}},(m_{(\zeta_{n}+\bar{\zeta}_{n})}^{_{(n)}})\widehat{\hphantom{a}}\,^{_{n}}(o\widetilde{\hphantom{a}}^{_{n}}),\dots,\big((m_{(\zeta_{n}+\bar{\zeta}_{n})}^{_{(n)}})\widehat{\hphantom{a}}\,^{_{n}}\big)^{\frac{\phi(n)}{2}-1}(o\widetilde{\hphantom{a}}^{_{n}})\right\}
\end{eqnarray*}
and, trivially,
$$
[B_{1}^{_{(n)}}o]\widetilde{\hphantom{a}}^{_{n}}\subset\big[\OO_{n}\cap (\R o)\big]\widetilde{\hphantom{a}}^{_{n}}\,.
$$
Secondly, one can easily see that, since the elements
$(b_{j}^{_{(n)}})\widetilde{\hphantom{a}}^{_{n}}\in (\R^2)^{\phi(n)/2}$,
$j\in\{1,\dots,\phi(n)/2\}$, are $\R$-linearly independent,
the elements
$(b_{j}^{_{(n)}}o)\widetilde{\hphantom{a}}^{_{n}}\in (\R^2)^{\phi(n)/2}$,
$j\in\{1,\dots,\phi(n)/2\}$, are $\R$-linearly
independent as well. This shows that one has
$$\operatorname{dim}\Big(\big\langle[\OO_{n}\cap (\R o)]\widetilde{\hphantom{a}}^{_{n}} \big\rangle_{\R}\Big)\geq\frac{\phi(n)}{2}\,.$$ Assume that one has $\operatorname{dim}(\langle[\OO_{n}\cap (\R o)]\widetilde{\hphantom{a}}^{_{n}}\rangle_{\R})> \phi(n)/2$, i.e., assume the existence of more
than $\phi(n)/2$ elements of $[\OO_{n}\cap(\R
  o)]\widetilde{\hphantom{a}}^{_{n}}$ that are $\R$-linearly
  independent.
The inclusion 
$$
\big[\OO_{n}\cap(\R o)\big]\widetilde{\hphantom{a}}^{_{n}}\subset\big[\mathbbm{K}_{n}\cap(\R o)\big]\widetilde{\hphantom{a}}^{_{n}}\,,
$$
then implies the existence of more
than $\phi(n)/2$ elements of $[\mathbbm{K}_{n}\cap(\R
  o)]\widetilde{\hphantom{a}}^{_{n}}$ that are linearly independent over $\R$ and hence, linearly
independent over
$\Q$. This is a contradiction since, by Lemma~\ref{hatiso} and
Lemma~\ref{lines}, one has
$\operatorname{dim}_{\Q}([\mathbbm{K}_{n}\cap(\R
o)]\widetilde{\hphantom{a}}^{_{n}})=\phi(n)/2$. Part (a) follows.

For part (b), we are done if we can show that the inclusion
\begin{eqnarray*}
\left[\mathbbm{k}_{n}o\right]\widetilde{\hphantom{a}}^{_{n}}=\big[\mathbbm{K}_{n}\cap(\R
  o)\big]\widetilde{\hphantom{a}}^{_{n}}&\subset &\Big\langle\big[\OO_{n}\cap(\R o)\big]\widetilde{\hphantom{a}}^{_{n}}\Big\rangle_{\R}\cap[\mathbbm{K}_{n}]\widetilde{\hphantom{a}}^{_{n}}\\
&\subset & \Big\langle\big[\mathbbm{K}_{n}\cap(\R
  o)\big]\widetilde{\hphantom{a}}^{_{n}}\Big\rangle_{\R}\cap[\mathbbm{K}_{n}]\widetilde{\hphantom{a}}^{_{n}}\\ &=&\big\langle[\mathbbm{k}_{n}o]\widetilde{\hphantom{a}}^{_{n}}\big\rangle_{\R}\cap[\mathbbm{K}_{n}]\widetilde{\hphantom{a}}^{_{n}}
\end{eqnarray*}
holds and, moreover, that the left-hand side and the
right-hand side of this inclusion are equal. Since, by
Lemma~\ref{lines}, $\mathbbm{K}_{n}\cap (\R o)=\mathbbm{k}_{n}o$, the
only non-trivial part of the inclusion is $$\big[\mathbbm{K}_{n}\cap(\R
  o)\big]\widetilde{\hphantom{a}}^{_{n}}\subset
\Big\langle\big[\OO_{n}\cap(\R
  o)\big]\widetilde{\hphantom{a}}^{_{n}}\Big\rangle_{\R}\cap[\mathbbm{K}_{n}]\widetilde{\hphantom{a}}^{_{n}}\,.$$Let $z\in \mathbbm{K}_{n}\cap (\R o)$. It follows that $z=\lambda o$ with $\lambda\in\mathbbm{k}_{n}$ by Proposition~\ref{p1}. We
consider the $\mathbbm{Q}$-coordinates of $\lambda$ with
respect to the $\mathbbm{Q}$-basis
$$\{1,(\zetan+\bar{\zeta}_{n}),(\zetan+\bar{\zeta}_{n})^2, \dots
,(\zetan+\bar{\zeta}_{n})^{\frac{\phi(n)}{2}-1}\}$$ of $\mathbbm{k}_{n}$
 (cf. Corollary~\ref{cr5}) and let $k\in \mathbbm{N}$ be the least
common multiple of all their denominators. Then, by Remark~\ref{r1}, we
get $k\lambda \in \oo_{n}$ and $k\lambda o \in
\OO_{n}\cap(\R o)$. Since $$(k\lambda
o)\widetilde{\hphantom{a}}^{_{n}}=k(\lambda
o)\widetilde{\hphantom{a}}^{_{n}}\in \big[\OO_{n}\cap(\R o)\big]\widetilde{\hphantom{a}}^{_{n}}\,,$$ one gets $z\widetilde{\hphantom{a}}^{_{n}}=(\lambda
o)\widetilde{\hphantom{a}}^{_{n}}\in \left\langle[\OO_{n}\cap(\R
  o)]\widetilde{\hphantom{a}}^{_{n}}\right\rangle_{\R}$. This proves
the inclusion. In order to prove equality of the left-hand side and the
right-hand side of this inclusion, consider the $\Q$-linear
automorphism
$(m_{1/o}^{_{(n)}})\widehat{\hphantom{a}}\,^{_{n}}$, of
$[\mathbbm{K}_{n}]\widetilde{\hphantom{a}}^{_{n}}$ together with its unique
extension to an $\R$-linear automorphism, say
$(m_{1/o}^{_{(n)}})\widetilde{\hphantom{a}}^{_{n}}$, of
$(\R^2)^{\phi(n)/2}$; cf. Remark~\ref{maut}. For the
left-hand side, one has
$$
\big(m_{\frac{1}{o}}^{_{(n)}}\big)\widetilde{\hphantom{a}}^{_{n}}\big[
[\mathbbm{k}_{n}o]\widetilde{\hphantom{a}}^{_{n}}\big]=[\mathbbm{k}_{n}]\widetilde{\hphantom{a}}^{_{n}}\,.
$$
Furthermore, for the
right-hand side, one has
\begin{eqnarray*}
\big(m_{\frac{1}{o}}^{_{(n)}}\big)\widetilde{\hphantom{a}}^{_{n}}\Big[
\big\langle[\mathbbm{k}_{n}o]\widetilde{\hphantom{a}}^{_{n}}\big\rangle_{\R}\cap[\mathbbm{K}_{n}]\widetilde{\hphantom{a}}^{_{n}}\Big]&=&\big(m_{\frac{1}{o}}^{_{(n)}}\big)\widetilde{\hphantom{a}}^{_{n}}\Big[\big\langle[\mathbbm{k}_{n}o]\widetilde{\hphantom{a}}^{_{n}}\big\rangle_{\R}\Big]\cap
\big(m_{\frac{1}{o}}^{_{(n)}}\big)\widetilde{\hphantom{a}}^{_{n}}\Big[[\mathbbm{K}_{n}]\widetilde{\hphantom{a}}^{_{n}}\Big]\\
&=&\left\langle\big(m_{\frac{1}{o}}^{_{(n)}}\big)\widetilde{\hphantom{a}}^{_{n}}\big[[\mathbbm{k}_{n}o]\widetilde{\hphantom{a}}^{_{n}}\big]\right\rangle_{\R}\cap
[\mathbbm{K}_{n}]\widetilde{\hphantom{a}}^{_{n}}\\
&=&\big\langle
  [\mathbbm{k}_{n}]\widetilde{\hphantom{a}}^{_{n}}\big\rangle_{\R}\cap [\mathbbm{K}_{n}]\widetilde{\hphantom{a}}^{_{n}}\,.
\end{eqnarray*}
It remains to prove that
$$
\big\langle
  [\mathbbm{k}_{n}]\widetilde{\hphantom{a}}^{_{n}}\big\rangle_{\R}\cap [\mathbbm{K}_{n}]\widetilde{\hphantom{a}}^{_{n}}=[\mathbbm{k}_{n}]\widetilde{\hphantom{a}}^{_{n}}\,.
$$
Clearly, one has $
\left\langle
  [\mathbbm{k}_{n}]\widetilde{\hphantom{a}}^{_{n}}\right\rangle_{\R}\cap [\mathbbm{K}_{n}]\widetilde{\hphantom{a}}^{_{n}}\supset[\mathbbm{k}_{n}]\widetilde{\hphantom{a}}^{_{n}}
$, so we are only concerned with the other inclusion. Let $x\in\left\langle
  [\mathbbm{k}_{n}]\widetilde{\hphantom{a}}^{_{n}}\right\rangle_{\R}\cap [\mathbbm{K}_{n}]\widetilde{\hphantom{a}}^{_{n}}$. Since $$\big\langle
  [\mathbbm{k}_{n}]\widetilde{\hphantom{a}}^{_{n}}\big\rangle_{\R}=\left\langle\big\langle [B_{1}^{_{(n)}}]\widetilde{\hphantom{a}}^{_{n}}\big\rangle_{\Q}\right\rangle_{\R}=\big\langle [B_{1}^{_{(n)}}]\widetilde{\hphantom{a}}^{_{n}}\big\rangle_{\R}$$ and $[\mathbbm{K}_{n}]\widetilde{\hphantom{a}}^{_{n}}=\left\langle [B_{\vphantom{1}}^{_{(n)}}]\widetilde{\hphantom{a}}^{_{n}}\right\rangle_{\Q}$, there are $\mu_1,\dots,\mu_{\phi(n)/2}\in\R$ and $q_1,\dots,q_{\phi(n)}\in\Q$ such that
$$
x=\sum_{j=1}^{\frac{\phi(n)}{2}}\mu_j
(b_j^{_{(n)}})\widetilde{\hphantom{a}}^{_{n}}=\sum_{j=1}^{\phi(n)}q_j (b_j^{_{(n)}})\widetilde{\hphantom{a}}^{_{n}}\,.
$$Since $[B_{\vphantom{1}}^{_{(n)}}]\widetilde{\hphantom{a}}^{_{n}}$
is an $\R$-basis of $(\R^2)^{\phi(n)/2}$, one obtains
$\mu_j=q_j$ and $q_{\phi(n)/2+j}=0$ for all
$j\in\{1,\dots,\phi(n)/2\}$. It follows that
$x\in\left\langle
  [B_{1}^{_{(n)}}]\widetilde{\hphantom{a}}^{_{n}}\right\rangle_{\Q}=[\mathbbm{k}_{n}]\widetilde{\hphantom{a}}^{_{n}}$. This completes the proof.
\end{proof}

From now on, we consider $(\R^2)^{\phi(n)/2}$ as an
Euclidean vector space, equipped with the inner product
$\langle \,.\, ,\, .\,\rangle$, defined by taking the basis
$[B_{\vphantom{1}}^{_{(n)}}]\widetilde{\hphantom{a}}^{_{n}}$ of the lattice $[\OO_{n}]\widetilde{\hphantom{a}}^{_{n}}$ to be an
orthonormal basis. In particular, the orthogonal complement
$T^{\perp}$ of a linear subspace $T$ of $(\R^2)^{\phi(n)/2}$
is understood with respect to this inner product, and the norm
$\Arrowvert x\Arrowvert$ of an element
$x\in(\R^2)^{\phi(n)/2}$ is given by
\begin{equation}\label{normo}
\Arrowvert x\Arrowvert=\sqrt{\sum_{j=1}^{\phi(n)}\lambda_j^2}\,,
\end{equation} where the $\lambda_j$ are
the coordinates of $x$ with respect to the basis
$[B_{\vphantom{1}}^{_{(n)}}]\widetilde{\hphantom{a}}^{_{n}}$. We wish to emphasize that the latter should not be confused with
  the use of orthogonal projections in
  the notion of determination of finite sets by projections, where we
  always use the {\em canonical} inner product on $\R^d$, $d\geq 2$.

\begin{rem}
Note that if $o\in \OO_{n}$, then, by Proposition~\ref{p1}(a),
one also has $\bar{o}\in \OO_{n}$, where $\bar{o}$ denotes the complex conjugate of $o$.
\end{rem}

\begin{lem}\label{transfer}
Let $n\,\in\,\mathbbm{N}\setminus \{1,2\}$, let $F\in \mathcal{F}(\OO_{n})$ and let $u\in\mathbb{S}^{1}$ be an
$\mathcal{O}_{n}$-direction, say parallel to $o\in
\OO_{n}\setminus\{0\}$. Then, the set
$$
G_{F}^{\{u\}}\,\,\cap\,\, \OO_{n}\,.
$$
corresponds via the $\Q$-linear isomorphism $$\,\,.\,\widehat{\hphantom{a}}\,^{_{n}}\;\circ\;m_{\bar{o}}^{_{(n)}}\!:\,\mathbbm{K}_{n}\longrightarrow [\mathbbm{K}_{n}]\widetilde{\hphantom{a}}^{_{n}}$$
to a subset of the lattice $[\OO_{n}]\widetilde{\hphantom{a}}^{_{n}}$ which is contained in a finite union of
translates of the form $y+[\oo_{n}]\widetilde{\hphantom{a}}^{_{n}}$, where $y\in [\OO_{n}]\widetilde{\hphantom{a}}^{_{n}}$.
\end{lem}
\begin{proof}
First, note that $m_{\bar{o}}$ is indeed a $\Q$-linear automorphism,
since $\bar{o}\neq 0$; see Remark~\ref{maut}. Let $\ell\in \operatorname{supp}(X_{u}F)$ and consider an element
$f+\lambda o \in \ell \cap \OO_{n}$, where $f\in F\subset \OO_{n}$ and
$\lambda \in \R$. One sees that $\lambda o \in \OO_{n}$ and, further,
$\lambda \in \mathbbm{k}_{n}$; recall that $\mathbbm{k}_{n}$ is the maximal real subfield of $\mathbbm{K}_{n}$. Since $\bar{o}\in\OO_{n}$, one has 
$(\lambda o)\bar{o} \in \oo_{n}$; cf. Proposition~\ref{p1}. This shows
that $m_{\bar{o}}^{_{(n)}}(f+\lambda
o)=f\bar{o}+(\lambda o)\bar{o}\in f\bar{o}+\oo_{n}$. Hence, one has
$m_{\bar{o}}^{_{(n)}}[\ell \cap \OO_{n}]\in f\bar{o}+\oo_{n}$. Consequently, $m_{\bar{o}}^{_{(n)}}[\,G_{F}^{\{u\}}\cap\OO_{n}\,]$ is contained in a finite union of
translates of the form $a+\oo_{n}$, where $a\in \OO_{n}$. The assertion follows immediately.
\end{proof}

\begin{rem}
Note that
 $\langle[\oo_{n}]\widetilde{\hphantom{a}}^{_{n}}\rangle_{\R}$ and $(\langle[\oo_{n}]\widetilde{\hphantom{a}}^{_{n}}\rangle_{\R})^{\perp}$ are
  $(\phi(n)/2)$-dimensional $[\OO_{n}]\widetilde{\hphantom{a}}^{_{n}}$-sub\-spaces, since one has
 $\langle[\oo_{n}]\widetilde{\hphantom{a}}^{_{n}}\rangle_{\R}=\langle[\mathcal{B}_{1}^{_{(n)}}]\widetilde{\hphantom{a}}^{_{n}}\rangle_{\R}$ and $(\langle[\oo_{n}]\widetilde{\hphantom{a}}^{_{n}}\rangle_{\R})^{\perp}=\langle[\mathcal{B}_{2}^{_{(n)}}]\widetilde{\hphantom{a}}^{_{n}}\rangle_{\R}$.
\end{rem}

In the next lemma, $(\langle[\oo_{n}]\widetilde{\hphantom{a}}^{_{n}}\rangle_{\R})^{\perp}$
  is seen as a metric space in the
  canonical way.

\begin{lem}\label{suitable}
Let $n\,\in\,\mathbbm{N}\setminus \{1,2\}$. Furthermore, let $B$ be an open ball  of
positive radius $r$ in the $(\phi(n)/2)$-dimensional
$[\OO_{n}]\widetilde{\hphantom{a}}^{_{n}}$-subspace
 $(\langle[\oo_{n}]\widetilde{\hphantom{a}}^{_{n}}\rangle_{\R})^{\perp}$. Then, there is a $(\phi(n)/2)$-dimensional
 $[\OO_{n}]\widetilde{\hphantom{a}}^{_{n}}$-subspace $S$ of
$(\R^2)^{\phi(n)/2}$ with the following properties:
\begin{itemize}
\item[(i)]
For each $y\in
[\OO_{n}]\widetilde{\hphantom{a}}^{_{n}}$, the translate $y+S$ meets at most one point of $$[\OO_{n}]\widetilde{\hphantom{a}}^{_{n}}\cap\left(B+\big\langle[\oo_{n}]\widetilde{\hphantom{a}}^{_{n}}\big\rangle_{\R}\right)\,.$$
\item[(ii)]
$S\cap [\mathbbm{K}_{n}]\widetilde{\hphantom{a}}^{_{n}}$ has a $\Q$-basis of 
the form 
$$
\left\{x,(m_{(\zeta_{n}+\bar{\zeta}_{n})}^{_{(n)}})\widehat{\hphantom{a}}\,^{_{n}}(x),\dots,\big((m_{(\zeta_{n}+\bar{\zeta}_{n})}^{_{(n)}})\widehat{\hphantom{a}}\,^{_{n}}\big)^{\frac{\phi(n)}{2}-1}(x) \right\}\,,
$$
where $x\in[\OO_{n}]\widetilde{\hphantom{a}}^{_{n}}\setminus\{0\}$. 
\end{itemize}
\end{lem}
\begin{proof}
For $\varepsilon \in \Q$ positive, set
$$S_{\varepsilon}:=\left\{x\in
(\R^2)^{\frac{\phi(n)}{2}}\,\left |\,\left\langle(b_{j}^{_{(n)}})\widetilde{\hphantom{a}}^{_{n}}+\varepsilon(b_{\frac{\phi(n)}{2}+j}^{_{(n)}})\widetilde{\hphantom{a}}^{_{n}},x\right\rangle=0
  \,\,\,\forall j\in\big\{1,\dots,\phi(n)/2\big\} \right.\right\}\,.$$ Then, $S_{\varepsilon}$ is a
 $[\OO_{n}]\widetilde{\hphantom{a}}^{_{n}}$-subspace of dimension $\phi(n)/2$. Further,
$[\mathbbm{K}_{n}]\widetilde{\hphantom{a}}^{_{n}}\cap S_{\varepsilon}$ is a vector space of dimension
$\phi(n)/2$ over $\Q$. Moreover, an
$\R$-basis of $S_{\varepsilon}$ which simultaneously is a $\Q$-basis of
$[\mathbbm{K}_{n}]\widetilde{\hphantom{a}}^{_{n}}\cap S_{\varepsilon}$ is given by
\begin{eqnarray*}
&&\left \{-\varepsilon (b_{j}^{_{(n)}})\widetilde{\hphantom{a}}^{_{n}} 
+(b_{\frac{\phi(n)}{2}+j}^{_{(n)}})\widetilde{\hphantom{a}}^{_{n}}\,\left
  |\, j\in \big \{1,\dots,\phi(n)/2 \big\}
 \right. \right \}\\
&=&\left \{\big(m_{(\zeta_{n}+\bar{\zeta}_{n})}^{_{(n)}})\widehat{\hphantom{a}}\,^{_{n}}\big)^{j}\big(-\varepsilon (b_{1}^{_{(n)}})\widetilde{\hphantom{a}}^{_{n}} 
+(b_{\frac{\phi(n)}{2}+1}^{_{(n)}})\widetilde{\hphantom{a}}^{_{n}}\big)\,\left
  |\, j\in \big \{0,\dots,\phi(n)/2-1 \big\}
 \right. \right \}.
\end{eqnarray*}
Note that there is even a suitable
$x_{\varepsilon}\in[\OO_{n}]\widetilde{\hphantom{a}}^{_{n}}\setminus\{0\}$ such that 
$$\Big\{\big(m_{(\zeta_{n}+\bar{\zeta}_{n})}^{_{(n)}})\widehat{\hphantom{a}}\,^{_{n}}\big)^{j}(x_{\varepsilon})\,|\,j\in\big\{0,\dots,\phi(n)/2-1\big\}\Big\}$$
is an $\R$-basis of $S_{\varepsilon}$ as well as a $\Q$-basis of
$[\mathbbm{K}_{n}]\widetilde{\hphantom{a}}^{_{n}}\cap S_{\varepsilon}$. For example, if $q_{\varepsilon}$ is the
denominator of $\varepsilon$, then 
$$
x_{\varepsilon}:= q_{\varepsilon} \Big (-\varepsilon (b_{1}^{_{(n)}})\widetilde{\hphantom{a}}^{_{n}} 
+\big(b_{\frac{\phi(n)}{2}+1}^{_{(n)}}\big)\widetilde{\hphantom{a}}^{_{n}}\Big)
$$
has this property.

Suppose that $v\in
([\OO_{n}]\widetilde{\hphantom{a}}^{_{n}}\setminus\{0\})\cap S_{\varepsilon}$, say
$v=\sum_{j=1}^{\phi(n)}\lambda_{j}(b_{j}^{_{(n)}})\widetilde{\hphantom{a}}^{_{n}}$ with uniquely determined
 $\lambda_{j}\in \Z$; see Remark~\ref{latticerem}. Hence, one has
$$
\lambda_{j}=-\varepsilon \lambda_{\frac{\phi(n)}{2}+j}
$$ 
for $j\in\{1,\dots,\phi(n)/2\}$. Further, Since $v\neq 0$, one
has $\lambda_{j}\neq 0$ for at least one $j\in\{1,\dots,\phi(n)\}$.
 If $v\notin (\langle[\oo_{n}]\widetilde{\hphantom{a}}^{_{n}}\rangle_{\R})^{\perp}$, there is a
 $j_{0}\in\{1,\dots,\phi(n)/2\}$ such that $\lambda_{j_{0}}\neq 0$. It follows that
$$
1\leq \lvert \lambda_{j_{0}} \rvert =\varepsilon \big\lvert
  \lambda_{\frac{\phi(n)}{2}+j_{0}}\big \rvert
 \leq \varepsilon \Arrowvert v \Arrowvert.
$$
If $v\in
(\langle[\oo_{n}]\widetilde{\hphantom{a}}^{_{n}}\rangle_{\R})^{\perp}$, then $\lambda_{j}=0$ for
$j\in\{1,\dots,\phi(n)/2\}$. This means that there is a
$j_{0}\in\{1,\dots,\phi(n)/2\}$ with 
$\lambda_{\frac{\phi(n)}{2}+j_{0}}\neq 0$. Then, one has
$$
0=-\varepsilon \lambda_{\frac{\phi(n)}{2}+j_{0}}\neq 0\,,
$$ 
a contradiction. 

One sees
that every $v\in
([\OO_{n}]\widetilde{\hphantom{a}}^{_{n}}\setminus\{0\})\cap S_{\varepsilon}$ satisfies $\Arrowvert v \Arrowvert\geq 1/\varepsilon$. Further, one
has
$$
\Big\Arrowvert v\big|\big\langle[\oo_{n}]\widetilde{\hphantom{a}}^{_{n}}\big\rangle_{\R} \Big\Arrowvert = \varepsilon \sqrt{\Big(\lambda_{\frac{\phi(n)}{2}+1}\Big)^{2}+\dots+\Big(\lambda_{\phi(n)}\Big)^{2}} \leq \varepsilon
\Arrowvert v \Arrowvert\,.
$$ 
It follows that, if $\varepsilon \leq 1/2$, then
$$
\Big\Arrowvert v\big|\big(\langle[\oo_{n}]\widetilde{\hphantom{a}}^{_{n}}\rangle_{\R}\big)^{\perp} \Big\Arrowvert \geq \Arrowvert v \Arrowvert -
\Big\Arrowvert v\big|\big\langle[\oo_{n}]\widetilde{\hphantom{a}}^{_{n}}\big\rangle_{\R} \Big\Arrowvert \geq \Arrowvert v \Arrowvert - \varepsilon \Arrowvert v
\Arrowvert = (1-\varepsilon) \Arrowvert v \Arrowvert \geq \frac{1}{2\varepsilon}\,.
$$
This implies that, if $v_{1},v_{2}\in [\OO_{n}]\widetilde{\hphantom{a}}^{_{n}}$ are distinct
lattice points in a translate $y+S_{\varepsilon}$ of
$S_{\varepsilon}$, where $y\in [\OO_{n}]\widetilde{\hphantom{a}}^{_{n}}$, and if $\varepsilon \leq
1/2$, then the distance $$\big\Arrowvert v_{1}|(\langle[\oo_{n}]\widetilde{\hphantom{a}}^{_{n}}\rangle_{\R})^{\perp} - v_{2}|(\langle[\oo_{n}]\widetilde{\hphantom{a}}^{_{n}}\rangle_{\R})^{\perp}\big\Arrowvert=\big\Arrowvert(v_{1}- v_{2})|(\langle[\oo_{n}]\widetilde{\hphantom{a}}^{_{n}}\rangle_{\R})^{\perp}\big\Arrowvert$$
between $v_{1}|(\langle[\oo_{n}]\widetilde{\hphantom{a}}^{_{n}}\rangle_{\R})^{\perp}$ and $v_{2}|(\langle[\oo_{n}]\widetilde{\hphantom{a}}^{_{n}}\rangle_{\R})^{\perp}$ is at least
$1/(2\varepsilon)$.

Choose an element $\varepsilon_{0}
\in \Q$ such that
$$
0<\varepsilon_{0}\leq \operatorname{min}\left
  \{\frac{1}{2},\frac{1}{4r} \right \}\,.
$$ 
Then, for each $y\in
[\OO_{n}]\widetilde{\hphantom{a}}^{_{n}}$, the translate $y+S_{\varepsilon_{0}}$ of
$S_{\varepsilon_{0}}$ meets at most one point of
$$[\OO_{n}]\widetilde{\hphantom{a}}^{_{n}}\cap \Big(B+\big\langle[\oo_{n}]\widetilde{\hphantom{a}}^{_{n}}\big\rangle_{\R}\Big)\,.$$ To see this, assume the existence of an element $y\in
[\OO_{n}]\widetilde{\hphantom{a}}^{_{n}}$ and assume the existence of two distinct
points $v_{1},v_{2}$ in $(y+S_{\varepsilon_{0}})\cap([\OO_{n}]\widetilde{\hphantom{a}}^{_{n}}\cap
(B+\langle[\oo_{n}]\widetilde{\hphantom{a}}^{_{n}}\rangle_{\R}))$. Then, one has
\begin{eqnarray*}
2r &>&\big\Arrowvert v_{1}|(\langle[\oo_{n}]\widetilde{\hphantom{a}}^{_{n}}\rangle_{\R})^{\perp} -
v_{2}|(\langle[\oo_{n}]\widetilde{\hphantom{a}}^{_{n}}\rangle_{\R})^{\perp}\big\Arrowvert\\&=&\big\Arrowvert(v_{1}- v_{2})|(\langle[\oo_{n}]\widetilde{\hphantom{a}}^{_{n}}\rangle_{\R})^{\perp}\big\Arrowvert
\,\,\geq\,\, \frac{1}{2\varepsilon_{0}}\,\,\geq\,\, 2r\,,
\end{eqnarray*}
a contradiction. We also saw above that a $\Q$-basis of $[\mathbbm{K}_{n}]\widetilde{\hphantom{a}}^{_{n}}\cap S_{\varepsilon_{0}}$
is given by 
$$\left\{\left.\big(m_{(\zeta_{n}+\bar{\zeta}_{n})}^{_{(n)}})\widehat{\hphantom{a}}\,^{_{n}}\big)^{j}(x_{\varepsilon_{0}})\,\right
  |\,j\in\big\{0,\dots,\phi(n)/2-1\big\}
  \right\}\,.$$
This completes the proof.
\end{proof}

Now we are able to prove a first result of this section on successive determination.

\begin{theorem}\label{mainsucdet}
Let $n\,\in\,\mathbbm{N}\setminus \{1,2\}$. Then, one has:
\begin{itemize}
\item[(a)]
The set $\mathcal{F}(\OO_{n})$ is successively determined by two $X$-rays in $\mathcal{O}_{n}$-directions.
\item[(b)]
The set
$\mathcal{F}(\OO_{n})$ is successively determined by two
projections on orthogonal complements of
$1$-dimensional $\OO_{n}$-subspaces.
\end{itemize}
\end{theorem}
\begin{proof}
Let us first prove part (a). Let $F\in \mathcal{F}(\OO_{n})$ and let $u\in\mathbb{S}^{1}$ be an
$\mathcal{O}_{n}$-direction. Let $o\in \OO_{n}\setminus\{0\}$ be
parallel to $u$. Suppose that $F'\in \mathcal{F}(\OO_{n})$ satisfies
$X_{u}F=X_{u}F'$. Then, by Lemma~\ref{fgrid}, one has
$$
F,F'\,\,\subset\,\, G_{F}^{\{u\}}\,\,\cap\,\, \OO_{n}\,.
$$
By Lemma~\ref{transfer}, the set $G_{F}^{\{u\}}\cap \OO_{n}$ corresponds via the $\Q$-linear isomorphism $$\,\,.\,\widehat{\hphantom{a}}\,^{_{n}}\;\circ\;m_{\bar{o}}^{_{(n)}}\!:\,\mathbbm{K}_{n}\longrightarrow [\mathbbm{K}_{n}]\widetilde{\hphantom{a}}^{_{n}}$$
to a subset of the lattice $[\OO_{n}]\widetilde{\hphantom{a}}^{_{n}}$ which is contained in a finite union of
translates of the form $y+[\oo_{n}]\widetilde{\hphantom{a}}^{_{n}}$, where $y\in [\OO_{n}]\widetilde{\hphantom{a}}^{_{n}}$. It follows that
$$\big(m_{\bar{o}}^{_{(n)}}\big)\widehat{\hphantom{a}}\,^{_{n}}\big[\,G_{F}^{\{u\}}\cap \OO_{n}\,\big]\big|\big(\langle[\oo_{n}]\widetilde{\hphantom{a}}^{_{n}}\rangle_{\R}\big)^{\perp}$$ is a finite set. Hence, it is
contained in an open ball $B$ of
positive radius $r$ in $(\langle[\oo_{n}]\widetilde{\hphantom{a}}^{_{n}}\rangle_{\R})^{\perp}$. By
Lemma~\ref{suitable}, there
is a $(\phi(n)/2)$-dimensional $[\OO_{n}]\widetilde{\hphantom{a}}^{_{n}}$-subspace $S$ of
$(\R^2)^{\phi(n)/2}$ with the following properties:
\begin{itemize}
\item[(i)]
For each $y\in
[\OO_{n}]\widetilde{\hphantom{a}}^{_{n}}$, the translate $y+S$ meets at most one point of $$[\OO_{n}]\widetilde{\hphantom{a}}^{_{n}}\cap\left(B+\big\langle[\oo_{n}]\widetilde{\hphantom{a}}^{_{n}}\big\rangle_{\R}\right)\,.$$
\item[(ii)]
$S\cap [\mathbbm{K}_{n}]\widetilde{\hphantom{a}}^{_{n}}$ has a $\Q$-basis of 
the form 
$$
\left\{x,(m_{(\zeta_{n}+\bar{\zeta}_{n})}^{_{(n)}})\widehat{\hphantom{a}}\,^{_{n}}(x),\dots,\big((m_{(\zeta_{n}+\bar{\zeta}_{n})}^{_{(n)}})\widehat{\hphantom{a}}\,^{_{n}}\big)^{\frac{\phi(n)}{2}-1}(x) \right\}\,,
$$
where $x\in[\OO_{n}]\widetilde{\hphantom{a}}^{_{n}}\setminus\{0\}$. 
\end{itemize}

In particular, it follows that $S\cap [\mathbbm{K}_{n}]\widetilde{\hphantom{a}}^{_{n}}$ corresponds via
$(\,\,.\,\widehat{\hphantom{a}}\,^{_{n}})^{-1}$ to a $\Q$-linear subspace of $\mathbbm{K}_{n}$ of dimension
$\phi(n)/2$, say $L$, with $\Q$-basis
$$
\left\{\left.\Big(m_{(\zeta_{n}+\bar{\zeta}_{n})}^{_{(n)}}\Big)^{j}(o')\,\right
  |\,j\in\big\{0,\dots,\phi(n)/2-1\big\}
  \right\}\,,
$$
where $o'$ is the unique element of $\OO_{n}$
satisfying $(o')\widetilde{\hphantom{a}}^{_{n}}=x$. Clearly, one has $o'\in \OO_{n}\setminus
\{0\}$. Now, Lemma~\ref{lines} immediately implies that
$$L=\mathbbm{K}_{n}\cap (\R o')\,.$$ 
Consider $(m_{\bar{o}}^{_{(n)}})^{-1}(o')=o'/\bar{o}\in \mathbbm{K}_{n}\setminus
\{0\}$. By Proposition~\ref{gau} and Remark~\ref{r1},
$(m_{\bar{o}}^{_{(n)}})^{-1}(o')$ is parallel to a non-zero
element of $\OO_{n}$. Let $u'\in\mathbb{S}^{1}$ be an $\OO_{n}$-direction parallel
to $(m_{\bar{o}}^{_{(n)}})^{-1}(o')$, e.g.,
$$u':=u_{(m_{\bar{o}}^{_{(n)}})^{-1}(o')}$$ (as defined in Definition~\ref{u4nq}(b)). 
We claim that $X_{u'}F=X_{u'}F'$ implies that $F=F'$. In order to
prove our claim, we
shall actually show that any line in the Euclidean plane of the form
$\ell_{u'}^{o''}$ (cf. Definition~\ref{xray..}), where $o''\in \OO_{n}$,
meets at most one point of the set $G_{F}^{\{u\}}\cap \OO_{n}$ defined above. To see this,
assume the existence of an element $o''\in \OO_{n}$, and assume the
existence of two distinct points $g$ and $g'$ in $\ell_{u'}^{o''}\cap
(G_{F}^{\{u\}}\cap \OO_{n})$. We claim that $m_{\bar{o}}^{_{(n)}}(g)$ and $m_{\bar{o}}^{_{(n)}}(g')$ are two distinct points in
$$(o''\bar{o}+L) \cap
m_{\bar{o}}^{_{(n)}}\left[G_{F}^{\{u\}}\cap \OO_{n}\right]\,.$$ 
To see this, let $h\in \ell_{u'}^{o''}\cap
(G_{F}^{\{u\}}\cap \OO_{n})$. It follows that there is a suitable $\lambda\in \R$ such that
$$h\,=\,o''+\lambda \big(m_{\bar{o}}^{_{(n)}}\big)^{-1}(o')\in G_{F}^{\{u\}}\cap \OO_{n}\,\subset\, \OO_{n}\,.$$
Proposition~\ref{p1} implies that $\lambda \in \mathbbm{k}_{n}$ and,
moreover, one gets
$$
m_{\bar{o}}^{_{(n)}}(h)\,=\,m_{\bar{o}}^{_{(n)}}\big(o''+\lambda
(m_{\bar{o}}^{_{(n)}})^{-1}(o')\big)\,=\,o''\bar{o}+\lambda
\frac{o'}{\bar{o}}\bar{o}\,=\,o''\bar{o}+\lambda
o'\,\in\, o''\bar{o}+L\,. 
$$
This proves the claim. Finally,
$(m_{\bar{o}}^{_{(n)}}(g))\widetilde{\hphantom{a}}^{_{n}}$ and $(m_{\bar{o}}^{_{(n)}}(g'))\widetilde{\hphantom{a}}^{_{n}}$ are two distinct points in 
\begin{eqnarray*}
&& \big((o''\bar{o})\widetilde{\hphantom{a}}^{_{n}}+[L]\widetilde{\hphantom{a}}^{_{n}}\big) \cap
\big[m_{\bar{o}}^{_{(n)}}[G_{F}^{\{u\}}\cap
\OO_{n}]\big]\widetilde{\hphantom{a}}^{_{n}}\\
&=& \big
((o''\bar{o})\widetilde{\hphantom{a}}^{_{n}}+(S\cap[\mathbbm{K}_{n}]\widetilde{\hphantom{a}}^{_{n}})\big ) \cap
\big[m_{\bar{o}}^{_{(n)}}[G_{F}^{\{u\}}\cap \OO_{n}]\big]\widetilde{\hphantom{a}}^{_{n}}\\
&\subset&  \big((o''\bar{o})\widetilde{\hphantom{a}}^{_{n}}+ S\big) \cap
\big[m_{\bar{o}}^{_{(n)}}[G_{F}^{\{u\}}\cap \OO_{n}]\big]\widetilde{\hphantom{a}}^{_{n}}\,,
\end{eqnarray*}
which contradicts Property~(i) above, since
$(o'\bar{o})\widetilde{\hphantom{a}}^{_{n}}\in
[\OO_{n}]\widetilde{\hphantom{a}}^{_{n}}$ and since
$$\big[m_{\bar{o}}^{_{(n)}}[G_{F}^{\{u\}}\cap \OO_{n}]\big]\widetilde{\hphantom{a}}^{_{n}}$$ is a subset of
$[\OO_{n}]\widetilde{\hphantom{a}}^{_{n}}\cap(B+\langle[\oo_{n}]\widetilde{\hphantom{a}}^{_{n}}\rangle_{\R})$.
This proves part (a). Part (b) follows immediately from an analysis of the proof of part (a).
\end{proof}

\begin{coro}\label{mainsucdetcoro}
Let $n\,\in\,\mathbbm{N}\setminus \{1,2\}$. Then, one has:
\begin{itemize}
\item[(a)]
The set $\cup_{t\in\R^2} \mathcal{F}(t+\OO_{n})$ is successively determined by three $X$-rays in $\mathcal{O}_{n}$-directions.
\item[(b)]
The set
$\cup_{t\in\R^2} \mathcal{F}(t+\OO_{n})$ is successively determined by three
projections on orthogonal complements of
$1$-dimensional $\OO_{n}$-subspaces.
\end{itemize}
\end{coro}
\begin{proof}
Let us first prove part (a). Let $U\subset \mathbb{S}^1$ be any set of two non-parallel $\OO_{n}$-directions, say $u,u'\in\mathbb{S}^1$, having the property
\begin{itemize}
\item[(E)]
There are
$o,o'\in\OO_{n}\setminus\{0\}$ with $u_o=u$ and $u_{o'}=u'$, and satisfying one of the equivalent conditions {\rm
  (i)}-{\rm (iii)} 
of Proposition~$\ref{proponeequiv}$.
\end{itemize}
Clearly, there are sets $U$ having property (E), e.g.,
 $U:=\{1,\zeta_{n}\}$ has this property. Let $F,F'\in \cup_{t\in\R^2}
\mathcal{F}(t+\OO_{n})$, say $F\in\mathcal{F}(t+\OO_{n})$ and
$F'\in\mathcal{F}(t'+\OO_{n})$, where $t,t'\in \R^2$, and suppose that
$F$ and $F'$ have the same $X$-rays in the directions of $U$. Then, by
Lemma~\ref{fgrid} and Theorem~\ref{oneequiv} in conjunction with
Property~(E), one obtains
\begin{equation}\label{fton1}
F,F'\,\,\subset \,\,G_{F}^{U}\,\,\subset\,\,t+\OO_{n}\,.
\end{equation}
If $F=\varnothing$, then, by Lemma~\ref{cardinality}(a), one also gets
$F'=\varnothing$. It follows that one may assume, without loss of generality, that $F$ and $F'$ are
non-empty. Then, since $F'\subset t'+\OO_{n}$, it follows from
Equation~(\ref{fton}) that $t+\OO_{n}$ meets $t'+\OO_{n}$, the latter
being equivalent to the identity $t+\OO_{n}=t'+\OO_{n}$. Hence, one
has $F-t,F'-t\in\mathcal{F}(\OO_{n})$, and, moreover, $F-t$ and
$F'-t$ have the same $X$-rays in the directions of $U$. In particular,
one has $X_{u}F=X_{u}F'$ for $u\in U$. Now, beginning with this direction $u$, one obviously can proceed as in the
proof of Theorem~\ref{mainsucdet}(a). Part (b) follows immediately from an analysis of the proof of part (a).
\end{proof}

With having the modelling of atomic constellations in mind, we conclude as follows.

\begin{coro}\label{cth6}
Let $n\,\in\,\mathbbm{N}\setminus \{1,2\}$ and let
$\varLambda\subset \R^2$ be a Delone set living on $\OO_{n}$. Then,
one has:
\begin{itemize}
\item[(a)]
The set $\mathcal{F}(\varLambda)$ $($resp.,
$\cup_{t\in\R^2}\mathcal{F}(t+\varLambda)$$)$ is successively determined
by two $($resp., three$)$ $X$-rays in $\mathcal{O}_{n}$-directions.
\item[(b)]
The set $\mathcal{F}(\varLambda)$ $($resp.,
$\cup_{t\in\R^2}\mathcal{F}(t+\varLambda)$$)$ is successively determined
by two $($resp., three$)$ 
projections on orthogonal complements of $1$-dimensional $\OO_{n}$-subspaces.
\end{itemize}
\end{coro}
\begin{proof}
This follows immediately from Theorem~\ref{mainsucdet} (resp., Corollary~\ref{mainsucdetcoro}).
\end{proof}

\begin{coro}\label{ccth6a}
Let $n\,\in\,\mathbbm{N}\setminus \{1,2\}$ and let
$\varLambda_{n}(t,W)\in\mathcal{M}(\OO_{n})$ be a cyclotomic model
set. Then, one has:
\begin{itemize}
\item[(a)]
The set $\mathcal{F}(\varLambda_{n}(t,W))$ is successively
determined by two  $X$-rays in $\mathcal{O}_{n}$- directions.
\item[(b)]
The set $\mathcal{F}(\varLambda_{n}(t,W))$ is successively
determined by two  projections on orthogonal complements of
$1$-dimensional $\OO_{n}$-subspaces.
\end{itemize}
\end{coro}
\begin{proof}
This follows immediately from Theorem~\ref{mainsucdet}.
\end{proof}

\begin{coro}\label{ccth6}
Let $n\,\in\,\mathbbm{N}\setminus \{1,2\}$. Then, for all windows  $W\subset
 (\R^2)^{\phi(n)/2-1}$, for all star maps $.^{\star_{n}}\! : \, \mathcal{O}_{n}\longrightarrow
(\R^2)^{\phi(n)/2-1}$ $($as described in
 Definition~$\ref{cyclodef}$$)$ and for all $R>0$, one has:
\begin{itemize}
\item[(a)]
The set $\cup_{\varLambda\in
   W^{\star_{n}}_{\mathcal{M}_{g}(\mathcal{O}_{n})}}\mathcal{F}(\varLambda)$ is successively
determined by three $X$-rays in $\mathcal{O}_{n}$-directions.
\item[(b)]
The set $\cup_{\varLambda\in
   W^{\star_{n}}_{\mathcal{M}_{g}(\mathcal{O}_{n})}}\mathcal{F}(\varLambda)$ is successively
determined by three projections on orthogonal complements of
$1$-dimensional $\OO_{n}$-subspaces.
\end{itemize}
\end{coro}
\begin{proof}
This follows immediately from Corollary~\ref{mainsucdetcoro}.
\end{proof}

\begin{rem}
Clearly, Theorem~\ref{mainsucdet} and its immediate implications
contained in Corollary~\ref{cth6} and Corollary~\ref{ccth6a} are
 optimal with respect to the number of $X$-rays (resp.,
 projections) used. Note that $\mathcal{F}(\R^{2})$ needs at least three $X$-rays
(resp., projections) for its successive determination; see~\cite[Corollary
  7.5]{GG}. Consequently, Corollary~\ref{cth6} and
  Corollary~\ref{ccth6a} show that cyclotomic model sets and even
  general Delone sets living on $\OO_{n}$ are closer to
lattices than to general point sets as far as successive determination is
concerned; see~\cite[Corollary 7.3]{GG}.
\end{rem}

We are now able to give an alternative proof of Corollary~\ref{coromod}.

\begin{coro}[cf. Corollary~\ref{coromod}]\label{cth6d}
Let $n\,\in\,\mathbbm{N}\setminus \{1,2\}$, let
$\varLambda_{n}(t,W)\in \mathcal{M}(\OO_{n})$ be a cyclotomic model
set, and let $R>0$. Then,
one has:
\begin{itemize}
\item[(a)]
The set $\mathcal{D}_{<R}(\varLambda_{n}(t,W))$ is determined by two $X$-rays in $\mathcal{O}_{n}$-directions.
\item[(b)]
The set $\mathcal{D}_{<R}(\varLambda_{n}(t,W))$ is determined by two projections on orthogonal complements of
$1$-dimensional $\OO_{n}$-subspaces.
\end{itemize}
\end{coro}
\begin{proof}[Second proof]
We may assume, without loss of generality, that $t=0$ and hence
$\varLambda_{n}(t,W)\subset \OO_{n}$. Let $\,\,.\,\widetilde{\hphantom{a}}^{_{n}}$ be the Minkowski embedding
 of $\OO_{n}$ that is used in the construction of $\varLambda_{n}(t,W)$, i.e., a map 
$$
\,\,.\,\widetilde{\hphantom{a}}^{_{n}}\,:\, \OO_{n}\longrightarrow
\R^2\times(\R^2)^{\frac{\phi(n)}{2}-1}\,,
$$
given by
$$
z\longmapsto \big(z,z^{\star_{n}}\big)\,,
$$where
$$z^{\star_{n}}=
\left(\sigma_2(z),\dots,\sigma_{\frac{\phi(n)}{2}}(z)\right)\,;$$ see
Section~\ref{deficyclo}. The assertion will follow from the ensuing analysis of the proof of
Theorem~\ref{mainsucdet}. Let us first prove (a). Let $F\in
\mathcal{D}_{<R}(\varLambda_{n}(t,W))$ and choose the $\mathcal{O}_{n}$-direction
$u:=1\in\mathbb{S}^{1}\cap \OO_{n}$. Suppose that $F'\in \mathcal{D}_{<R}(\varLambda_{n}(t,W))$ satisfies
$X_{1}F=X_{1}F'$. Then, by Lemma~\ref{fgrid}, one has
$$
F,F'\,\,\subset\,\, G_{F}^{\{1\}}\,\,\cap\,\, \varLambda_{n}(t,W)\,\,\subset\,\,\OO_{n}\,.
$$
An analysis of the proof of Lemma~\ref{transfer} shows that the set $G_{F}^{\{1\}}\cap \varLambda_{n}(t,W)$ maps via $\,\,.\,\widetilde{\hphantom{a}}^{_{n}}$
to a subset of the lattice $[\OO_{n}]\widetilde{\hphantom{a}}^{_{n}}$ which is contained in a finite union of
translates of the form
$f\widetilde{\hphantom{a}}^{_{n}}+[\oo_{n}]\widetilde{\hphantom{a}}^{_{n}}$,
where $f\in F$. In particular, one has
$$
\big[G_{F}^{\{1\}}\cap \varLambda_{n}(t,W)\big]\widetilde{\hphantom{a}}^{_{n}}\,\,\subset\,\, [F]\widetilde{\hphantom{a}}^{_{n}}+[\oo_{n}]\widetilde{\hphantom{a}}^{_{n}}
$$
and, further,
\begin{equation}\label{fgeq}
\big[\,G_{F}^{\{1\}}\cap \varLambda_{n}(t,W)\,\big]\widetilde{\hphantom{a}}^{_{n}}\big|\big(\langle[\oo_{n}]\widetilde{\hphantom{a}}^{_{n}}\rangle_{\R}\big)^{\perp}\,\,\subset \,\,[\,F\,]\widetilde{\hphantom{a}}^{_{n}}\big|\big(\langle[\oo_{n}]\widetilde{\hphantom{a}}^{_{n}}\rangle_{\R}\big)^{\perp}\,.
\end{equation}
Clearly,
$[\,F\,]\widetilde{\hphantom{a}}^{_{n}}|(\langle[\oo_{n}]\widetilde{\hphantom{a}}^{_{n}}\rangle_{\R})^{\perp}$
is a finite set. Further, by definition of $\varLambda_{n}(t,W)$, one has
$$
[\,F\,]\widetilde{\hphantom{a}}^{_{n}}\,=\,\{(f,f^{\star_{n}})\,|\,f\in
F\}\,\subset\, F\times [F]^{\star_{n}}\,\subset\, F\times W\,.
$$
Since $F\in
\mathcal{D}_{<R}(\varLambda_{n}(t,W))$ by assumption and since $W$ is
bounded (recall that $\overline{W}$ is compact), one sees that there
is a positive $D\in\R$ such that, for all
$F\in\mathcal{D}_{<R}(\varLambda_{n}(t,W))$, the set
$[\,F\,]\widetilde{\hphantom{a}}^{_{n}}$ has diameter at most $D$ with respect to the maximum norm 
 on $\R^2\times(\R^2)^{\phi(n)/2-1}$, defined by the Euclidean norms on $\R^2$ and
$(\R^2)^{\phi(n)/2-1}\cong \R^{\phi(n)-2}$, respectively. Since
all norms on $\R^2\times(\R^2)^{\phi(n)/2-1}$ are equivalent,
there
is a positive $D'\in\R$ such that, for all
$F\in\mathcal{D}_{<R}(\varLambda_{n}(t,W))$, the set $
[\,F\,]\widetilde{\hphantom{a}}^{_{n}}$ has diameter at most $D'$ with respect to the norm on
$\R^2\times(\R^2)^{\phi(n)/2-1}$ defined by the
lattice $[\OO_{n}]\widetilde{\hphantom{a}}^{_{n}}$;
cf. Equation~(\ref{normo}). It follows from~(\ref{fgeq}) that, for all
$F\in\mathcal{D}_{<R}(\varLambda_{n}(t,W))$, the set
$[\,G_{F}^{\{1\}}\cap \varLambda_{n}(t,W)\,]\widetilde{\hphantom{a}}^{_{n}}|(\langle[\oo_{n}]\widetilde{\hphantom{a}}^{_{n}}\rangle_{\R})^{\perp}$
has diameter at most $D'$ with respect to the induced norm on $(\langle[\oo_{n}]\widetilde{\hphantom{a}}^{_{n}}\rangle_{\R})^{\perp}$.
Hence, there
is a positive $r\in\R$ such that, for all
$F\in\mathcal{D}_{<R}(\varLambda_{n}(t,W))$, the set
$[\,G_{F}^{\{1\}}\cap \varLambda_{n}(t,W)\,]\widetilde{\hphantom{a}}^{_{n}}|(\langle[\oo_{n}]\widetilde{\hphantom{a}}^{_{n}}\rangle_{\R})^{\perp}$ is
contained in a suitable open ball $B_{F}$ of
 radius $r$ in
 $(\langle[\oo_{n}]\widetilde{\hphantom{a}}^{_{n}}\rangle_{\R})^{\perp}$. Observing that the $(\phi(n)/2)$-dimensional
 $[\OO_{n}]\widetilde{\hphantom{a}}^{_{n}}$-subspace $S$ of
$(\R^2)^{\phi(n)/2}$ in Lemma~\ref{suitable} does only depend
on the radius $r$ of the open ball, one sees that it is possible here
to choose the second direction $u'$ independently from $F$ and
$F'$. This proves part (a). Part (b) again follows immediately from an analysis of the proof of part (a).
\end{proof}

\begin{rem}
Similarly, one can modify~\cite[Lemma 7.1]{GG} in order to obtain an alternative
 proof of Corollary~\ref{bounded2}. Clearly, Corollary~\ref{cth6d} is 
 best possible with respect to the number of $X$-rays (resp.,
 projections) used.
\end{rem}

\begin{theorem}\label{cth6d2}
Let $n\in\N\setminus\{1,2\}$. Then, for all windows  $W\subset
 (\R^2)^{\phi(n)/2-1}$, for all star maps $.^{\star_{n}}\! : \, \mathcal{O}_{n}\longrightarrow
(\R^2)^{\phi(n)/2-1}$ $($as described in
 Definition~$\ref{cyclodef}$$)$ and for all $R>0$, one has:
\begin{itemize}
\item[(a)]
The set $\cup_{\varLambda\in W^{\star_{n}}_{\mathcal{M}_{g}(\mathcal{O}_{n})}}\mathcal{D}_{<R}(\varLambda)$ is 
determined by three $X$-rays in $\mathcal{O}_{n}$-directions.
\item[(b)]
The set $\cup_{\varLambda\in W^{\star_{n}}_{\mathcal{M}_{g}(\mathcal{O}_{n})}}\mathcal{D}_{<R}(\varLambda)$ is 
determined by three projections on orthogonal complements of
$1$-dimensional $\OO_{n}$-subspaces.
\end{itemize}
 \end{theorem}
\begin{proof}
To prove part (a), let $U\subset \mathbb{S}^1$ be a set of two non-parallel
$\OO_{n}$-directions containing the $\OO_{n}$-direction $1$, say
$U=\{1,u\}$. Suppose that $u\in \mathbb{S}^1$ has the property
\begin{itemize}
\item[(E')]
There is an element 
 $o\in\OO_{n}\setminus\{0\}$ with $u_o=u$ such that $\{1,o\}$ satisfy one of the equivalent conditions {\rm
  (i)}-{\rm (iii)} 
of Proposition~$\ref{proponeequiv}$.
\end{itemize}
Clearly, there are elements $u\in \mathbb{S}^1$ having
Property~(E'), such as $u:=\zeta_{n}\in\mathbb{S}^1$. Let $F,F'\in \cup_{\varLambda\in
   W^{\star_{n}}_{\mathcal{M}_{g}(\mathcal{O}_{n})}}\mathcal{D}_{<R}(\varLambda)$, say $F\in\mathcal{D}_{<R}(\varLambda_{n}^{\star_{n}}(t,\tau+W))$ and $F'\in\mathcal{D}_{<R}(\varLambda_{n}^{\star_{n}}(t',\tau'+W))$, where $t,t'\in\R^2$ and $\tau,\tau'\in(\R^2)^{\phi(n)/2-1}$, and suppose that $F$ and $F'$ have the same $X$-rays in the directions of $U$. Then, by Lemma~\ref{fgrid} and Theorem~\ref{oneequiv} in conjunction with Property~(E'), one obtains
\begin{equation}\label{ftonfinal}
F,F'\,\,\subset \,\,G_{F}^{U}\,\,\subset\,\,t+\OO_{n}\,.
\end{equation}
If $F=\varnothing$, then, by Lemma~\ref{cardinality}(a), one also gets
$F'=\varnothing$. It follows that one may assume, without loss of generality, that $F$ and $F'$ are
non-empty. Then, since $F'\subset t'+\OO_{n}$, it follows from
Equation~(\ref{fton}) that $t+\OO_{n}$ meets $t'+\OO_{n}$, the latter
being equivalent to the identity $t+\OO_{n}=t'+\OO_{n}$. Moreover, 
 the identity $t+\OO_{n}=t'+\OO_{n}$ is equivalent to the relation $t'-t\in\OO_{n}$. Hence, one
has $$F-t\in\mathcal{D}_{<R}\Big(\varLambda_{n}^{\star_{n}}\big(0,\tau+W\big)\Big)$$
and, since the equality
$\varLambda_{n}^{\star_{n}}(t'-t,\tau'+W)=\varLambda_{n}^{\star_{n}}(0,(\tau'+(t'-t)^{\star_{n}})+W)$
holds,
$$F'-t\in\mathcal{D}_{<R}\Big(\varLambda_{n}^{\star_{n}}\big(t'-t,\tau'+W\big)\Big)=\mathcal{D}_{<R}\Big(\varLambda_{n}^{\star_{n}}\big(0,(\tau'+(t'-t)^{\star_{n}})+W\big)\Big)\,.$$
Clearly, $F-t$ and $F'-t$ again have the same $X$-rays in the
directions of $U$. In particular, one has
$X_{1}(F-t)=X_{1}(F'-t)$. Now, proceeding as in the proof of
Corollary~\ref{cth6d}, one sees that there is an $\OO_{n}$-direction
 $u'$ such that the set 
$\cup_{\varLambda\in
   W^{\star_{n}}_{\mathcal{M}_{g}(\mathcal{O}_{n})}}\mathcal{D}_{<R}(\varLambda)$
 is determined by the $X$-rays in the directions of the set 
 $U':=U\cup\{u'\}$. In particular, the above set is determined by
 three $X$-rays in $\OO_{n}$-directions. Part (b) follows immediately from an analysis of the proof of part (a).
\end{proof}

\section*{Final Remark}
Note that many of the above results also hold for model
sets associated with the well-known Penrose tiling of the plane; cf.~\cite{BH,H}. 

\section*{Towards Three-Dimensional Model Sets}
Though we have focussed on the planar case here, many results can be
lifted to higher dimensions. This is due to the dimensional hierarchy
mentioned in the introduction; 
see~\cite{PABP2}. However, model sets in dimension
$d\geq 3$ do not possess such a systematic description as the planar
ones do on the basis of cyclotomic fields. Also, particular symmetries
play an important role, such as icosahedral symmetries in
$3$-space. This will require some more detailed attention to the
special cases at hand. We hope to report on progress in this direction
in the near future; see~\cite{H}.

\section*{Acknowledgements}
It is a pleasure to thank Uwe
Grimm, Peter Gritzmann, Barbara Langfeld and Bernd Sing for helpful comments and
suggestions. The author is most grateful to Michael Baake, Richard J. Gardner  and Peter
A. B. Pleasants for valuable discussions and suggested
simplifications.
The author was supported by the German Research Council
  (Deutsche Forschungsgemeinschaft), within the Collaborative Research
  Centre (Sonderforschungsbereich) 701.

\end{document}